\pdfoutput=1

\documentclass[letterpaper, preprint, paper,10pt]{AAS}

\usepackage[utf8]{inputenc} %
\usepackage[T1]{fontenc}    %
\usepackage{kotex}
\usepackage{amsmath}
\usepackage{amsfonts}
\usepackage{amssymb}
\usepackage{mathtools}
\usepackage{bm}

\usepackage{graphicx}
\graphicspath{{media/}}     %

\usepackage{subcaption}
\newif\ifdraft
\drafttrue

\ifdraft
  \usepackage[colorlinks=true, pdfstartview=FitV, linkcolor=blue, citecolor=teal, urlcolor=magenta]{hyperref}
\else
  \usepackage[colorlinks=true, pdfstartview=FitV, linkcolor=black, citecolor=black, urlcolor=black]{hyperref}
\fi
\usepackage{xurl}
\usepackage{url}
\usepackage{doi}

\usepackage{booktabs}       %
\usepackage{tabularx}       %
\usepackage{multirow}       %
\usepackage{makecell}       %
\usepackage{nicefrac}       %
\usepackage{microtype}      %
\usepackage{soul}
\usepackage[dvipsnames]{xcolor}
\usepackage{tikz}
\usetikzlibrary{shapes.geometric}
\usepackage{tikz-3dplot}
\usepackage{circuitikz}
\usepackage{tablefootnote} 
\usepackage{enumitem}

\usepackage{pifont}

\usepackage[square,numbers]{natbib}

\usepackage[acronym]{glossaries}
\newacronym{cr3bp}{CR3BP}{Circular Restricted Three-Body Problem}
\newacronym{er3bp}{ER3BP}{Elliptic Restricted Three-Body Problem}
\newacronym{bcr4bp}{BCR4BP}{Bi-Circular Restricted Four-Body Problem}
\newacronym{hr3bp}{HR3BP}{Hill Restricted Three-Body Problem}
\newacronym{hr4bp}{HR4BP}{Hill Restricted Four-Body Problem}
\newacronym{qbcp}{QBCP}{Quasi Bi-Circular Problem}
\newacronym{qhr4bp}{QHR4BP}{Quasi-Hill Restricted Four-Body Problem}
\newacronym{qhr4bps}{QHR4BPs}{Quasi-Hill Restricted Four-Body Problems}

\newacronym{iqhr4bp}{I-QHR4BP}{In-plane Quasi-Hill Restricted Four-Body Problem}
\newacronym{oqhr4bp}{O-QHR4BP}{Out-of-plane Quasi-Hill Restricted Four-Body Problem}

\newacronym{hfem}{HFEM}{Higher-Fidelity Ephemeris Model}

\newacronym{nrho}{NRHO}{Near Rectilinear Halo Orbit}

\newacronym{qpo}{QPO}{Quasi-Periodic Orbit}
\newacronym{qpos}{QPOs}{Quasi-Periodic Orbits}

\newacronym{myt}{MYT}{Multi-Year Trajectory}
\newacronym{myts}{MYTs}{Multi-Year Trajectories}

\newacronym{gmos}{GMOS}{G{\'o}mez Mondelo Olikara Scheeres}

\newacronym{po}{PO}{Periodic Orbit}
\newacronym{pos}{POs}{Periodic Orbits}

\newacronym{dft}{DFT}{Discrete Fourier Transform}
\newacronym{fft}{FFT}{Fast Fourier Transform}
\newacronym{cft}{CFT}{Continuous Fourier Transform}

\newacronym{dro}{DRO}{Distant Retrograde Orbit}
\newacronym{stm}{STM}{State Transition Matrix}

\newacronym{qpts}{QPTs}{Quasi-Periodic Trajectories}

\newacronym{elfos}{ELFOs}{Elliptical Lunar Frozen Orbits}
\newacronym{elfo}{ELFO}{Elliptical Lunar Frozen Orbit}
\newacronym{lsp}{LSP}{Lunar South Pole}
\newacronym{gdop}{GDOP}{Geometric Dilution of Precision}
\newacronym{lcrns}{LCRNS}{Lunar Communications Relay and Navigation Systems}

\newacronym{pso}{PSO}{Particle Swarm Optimizer}
\newacronym{ga}{GA}{Genetic Algorithm}
\newacronym{de}{DE}{Differential Evolution}
\newacronym{aco}{ACO}{Ant Colony Optimizer}

\newacronym{laskar}{L-NAFF}{Laskar-Numerical Analysis of Fundamental Frequency}
\newacronym{gomez}{GMS-C}{G{\'o}mez-Mondelo-Sim{\'o}-Collocation}

\newacronym{fddc}{FDDC}{Frequency-Domain Differential Corrector}
\newacronym{dadm}{DADM}{Doubly Averaged Dynamical Model}

\newacronym{brf}{BRF}{Barycentric Rotating Frame}
\newacronym{mrf}{MRF}{Moon-Centered Rotating Frame}
\newacronym{eof}{EOF}{Moon-Centered Earth Orbit Frame}
\newacronym{mci}{MCI}{Moon-Centered Inertial}
\glsdisablehyper

\ifdraft
  \newcommand{\us}[1]{{ }(\textcolor{blue}{#1})}
  \newcommand{\usc}[1]{{ }(\textcolor{orange}{[C] #1})}
  \newcommand{\uso}[1]{{ }(\textcolor{Purple}{[O] #1})}
\else
  \newcommand{\us}[1]{}
  \newcommand{\usc}[1]{}
  \newcommand{\uso}[1]{}
\fi

\newcommand{\freeVE}{\ensuremath{\bm{X}}}

\newcommand{\constraintV}[1]{\ensuremath{\bm{F}_{#1}}}

\newcommand{\stm}[4]{\mathbf{\Phi}_{#1}^{#2}(#3; #4)}

\newcommand{\mb}{\begin{bmatrix}}
\newcommand{\me}{\end{bmatrix}}

\newcommand{\redpentagon}{\protect\tikz[baseline=-0.4ex]\protect\node[regular polygon, regular polygon sides=5, fill=red, draw=red, minimum size=1.1ex, inner sep=0pt]{};}

\definecolor{peakA}{HTML}{E41A1C}
\definecolor{peakB}{HTML}{377EB8}
\definecolor{peakC}{HTML}{4DAF4A}
\definecolor{peakD}{HTML}{984EA3}
\definecolor{peakE}{HTML}{FF7F00}
\definecolor{peakF}{HTML}{A65628}
\newcommand{\pkmark}[1]{\protect\tikz[baseline=-0.4ex]\protect\fill[#1] circle(0.55ex);}

\PaperNumber{26-738}

\newif\ifarxiv
\arxivtrue

\ifarxiv
  \usepackage{etoolbox}
  \makeatletter
  \patchcmd{\@maketitle}{(Preprint) AAS~\AA@papernumber}{}{}%
           {\typeout{AASPATCH WARNING: preprint stamp not removed}}
  \patchcmd{\@maketitle}{AAS~\AA@papernumber}{}{}%
           {\typeout{AASPATCH WARNING: plain AAS stamp not removed}}
  \makeatother
\fi

\title{Elliptical Lunar Frozen Orbit Constellations: \\ Torus-Based Design and Analysis}

\author{%
Beom Park\thanks{Apollo 11 Postdoctoral Fellow, School of Aeronautics and Astronautics, Purdue University, West Lafayette, IN 47907, USA.},
Kathleen C. Howell\thanks{Hsu Lo Distinguished Professor, School of Aeronautics and Astronautics, Purdue University, West Lafayette, IN 47907, USA.},
Daniel Brack\thanks{Lunar Data Network Trajectory Lead, Intuitive Machines, 13467 Columbia Shuttle Street, Houston, TX 77059, USA.}
}

\begin{document}

\setcounter{secnumdepth}{3}
\ifdraft\setcounter{tocdepth}{3}\fi   %

\maketitle{}

\renewcommand{\thefootnote}{\arabic{footnote}}

\begin{abstract}
Elliptical Lunar Frozen Orbits (ELFOs) are attractive candidates for lunar satellite constellations supporting lunar south pole exploration. Building on prior frequency-based orbit analysis methods, this investigation develops a torus-based, frequency-domain framework for constellation-level ELFO design and analysis. Within a doubly-averaged dynamical model, each ELFO is characterized by three frequency-angle pairs; a systematic comparison across progressively higher-fidelity models identifies the additional frequency components introduced at each level. The angular parametrization admits multiple symmetries that substantially reduce the design space, enabling a rapid, epoch-free global survey of constellation configurations. A higher-fidelity refinement and analysis layer then combines a (1) \acrfull{fddc} that targets desired frequency properties and suppresses undesired long-period oscillations, (2) Fourier surrogate re-optimization for inter-satellite phasing without propagation in the loop, and (3) spectral attribution that separates design-refinable coverage losses from those inherent to the dynamics. Together, these elements supply a fidelity-bridging and diagnostic framework for \acrshort{elfo} constellation design and analysis.
\end{abstract}

\ifarxiv\else
  \ifdraft
  {
  \hypersetup{linkcolor=black}
  \tableofcontents
  }
  \clearpage
  \fi
\fi

\section{Introduction}

\acrfull{elfos} are nominally bounded trajectories in the vicinity of the Moon that display quasi-constant eccentricity, inclination, and argument of perilune~\cite{von1910application,lidov1962evolution,kozai1959effects}. The frozen characteristic of these highly inclined, eccentric orbits arises primarily from the Earth's gravitational perturbation, and the resulting geometry supplies consistent line-of-sight to the \acrfull{lsp} and its vicinity~\cite{ely2005stable,ely2006constellations,folta2006lunar}. With these characteristics, multiple space agencies are envisioning lunar satellite constellations that incorporate one or more \acrshort{elfos} into the design architecture. Notable examples include LunaNet~\cite{israel2020lunanet}, a NASA-led initiative pursued with international (ESA, JAXA) and commercial partners such as Intuitive Machines to establish a ``lunar internet'' for navigation and communication services, and the Queqiao relay satellite constellation~\cite{jones2024queqiao} developed along a comparable concept by CNSA.

As lunar satellite constellations advance from concept to implementation, the strategic placement of the constituent satellites becomes a key design challenge. Existing constellation analyses exploit closed-form features of \acrshort{elfos}, including the frozen condition and lunar period resonances, to reduce the effective design space. Candidate configurations are refined through finite-horizon optimization, with the propagation fidelity ranging from short-arc Keplerian evaluations~\cite{ceresoli2025design,zanotti2024high} to higher-fidelity propagators in the optimization loop~\cite{brack2025}. Specific degradation mechanisms are also identified and corrected at the case level. \citet{ely2006constellations} identifies the relative mean-anomaly drift between satellites as a primary source of phasing degradation, mitigated through a careful selection of initial semi-major axes. The same mechanism is operationalized by \citet{brack2025} via mean-to-osculating refinement of the initial conditions through a particle-swarm optimization, and by \citet{ceresoli2025design} via a tailored analytical semi-major-axis correction. These efforts deliver mission-ready point designs and isolate the dominant mechanism for each reference configuration.

The current investigation extends this body of work through a torus-based, frequency-domain route, synthesizing frequency-based ingredients from the authors' prior work into a constellation-level design and diagnosis framework. The contributions of this integration are three-fold. First, the frequency structure of \acrshort{elfos} is examined across a hierarchy of dynamical models~\cite{park2025numerical}; each orbit is characterized by three frequency-angle pairs within the simplest model, and a systematic comparison with higher-fidelity dynamics classifies the additional frequency components introduced at each level as controllable through initial conditions or unremovable, a distinction not made explicit in prior treatments. Second, the angular parametrization yields multiple symmetries~\cite{park2026bridging} that reduce the effective design space, enabling a rapid, epoch-free global survey that reveals structural features of the constellation design landscape. Third, a higher-fidelity refinement-and-analysis layer bridges the lower-fidelity survey to the target ephemeris dynamics. The layer combines a \acrfull{fddc} refinement~\cite{park2025frequency} of targeted frequency properties and undesired long-period oscillations with a Fourier surrogate that re-optimizes the inter-satellite phasing without propagation in the loop. A frequency-domain attribution then separates design-refinable coverage losses from forced components that are inherent to the dynamics. 

The remainder of this document is organized as follows. Section~\ref{sec:prelim} supplies background on the reference frames as well as dynamical models. Section~\ref{sec:freq-structure} examines the frequency structure of the \acrshort{elfos} across dynamical models. Section~\ref{sec:design-exploration} leverages this structure for a systematic exploration of the constellation design space within the lower-fidelity model. Section~\ref{sec:fddc} pursues refinement and analysis in a higher-fidelity model. Lastly, Section~\ref{sec:conclusions} offers concluding remarks.

\section{Preliminaries}
\label{sec:prelim}
This section establishes the reference frames and dynamical models employed throughout the analysis. Three models of increasing fidelity are introduced: the \acrfull{dadm}, the \acrfull{cr3bp}, and the \acrfull{hfem}. The hierarchical relationship among these models enables the systematic frequency-based analysis of \acrshort{elfo} evolution developed in subsequent sections. The conventions for the frames as well as the dynamical models follow the authors' previous treatment~\cite{park2026bridging}.

\subsection{Reference Frames}

The dimensional gravitational parameters for the Earth and Moon are $\tilde{\mu}_E$ and $\tilde{\mu}_M$, respectively, and the mass ratio is $\mu = \tilde{\mu}_M / (\tilde{\mu}_E + \tilde{\mu}_M) \approx 0.01215$. A reference Earth-Moon distance $l_* \approx 3.847479920112920\cdot 10^5$ km defines the characteristic time $t_* = \sqrt{l_*^3 / (\tilde{\mu}_E + \tilde{\mu}_M)}\approx 3.756998590849907\cdot 10^5$ s; these characteristic quantities are employed to nondimensionalize all distance and time variables, respectively. The nondimensional (nd) time is denoted $t$, and an overdot indicates differentiation with respect to $t$. Four reference frames are employed as follows:
\begin{itemize}[topsep=0pt, partopsep=0pt]
    \item \ul{Moon-Centered Rotating Frame} (\acrshort{mrf}): An orthogonal basis $\hat{\bm{x}}$-$\hat{\bm{y}}$-$\hat{\bm{z}}$ is constructed with $\hat{\bm{x}}$ directed from the Earth to the Moon, $\hat{\bm{z}}$ along the Earth-Moon orbital angular momentum vector, and $\hat{\bm{y}}$ completing the right-handed triad. The origin is at the Moon. The nd spacecraft position is $\bm{r} = x\hat{\bm{x}} + y\hat{\bm{y}} + z\hat{\bm{z}}$, nondimensionalized\footnote{An alternative, ``pulsating'' convention nondimensionalizes the spacecraft position by the instantaneous Earth-Moon distance $l$ \cite{park2024assessment}, applicable to the \acrshort{hfem} dynamics. While this pulsating formulation embeds the Earth-Moon distance variation and is useful in describing dynamics near libration points (where orbital geometry relative to the Earth-Moon is relevant), the \acrshort{elfos} are near-Keplerian with respect to the Moon, forgoing the need for the time-varying distance in the scaling process.} by $l_*$, and the velocity is $\bm{v} = \dot{\bm{r}}$. In the \acrshort{cr3bp} and \acrshort{dadm}, the basis is fixed by the constant Earth-Moon geometry; in the \acrshort{hfem}, the basis directions track the instantaneous Earth-Moon geometry supplied at each epoch via the ephemerides. The \acrshort{mrf} is primarily adopted for trajectory visualization throughout this work; dimensional position vectors are recovered as $l_*\bm{r}$ and reported in km.

    \item \ul{Barycentric Rotating Frame} (\acrshort{brf}): The \acrshort{brf} shares the $\hat{\bm{x}}$-$\hat{\bm{y}}$-$\hat{\bm{z}}$ basis of the \acrshort{mrf} but is centered at the Earth-Moon barycenter. It is adopted for the \acrshort{cr3bp} formulation, where the Earth-Moon distance is constant and equals $l_*$; the Earth and Moon are located at $\bm{r}_E = -\mu\hat{\bm{x}}$ and $\bm{r}_M = (1-\mu)\hat{\bm{x}}$, respectively. The nd spacecraft state is $(\bm{r}_B,\, \bm{v}_B)$ with $\bm{r}_B = x_B\hat{\bm{x}} + y_B\hat{\bm{y}} + z_B\hat{\bm{z}}$ and $\bm{v}_B = \dot{\bm{r}}_B$, related to the \acrshort{mrf} state via $\bm{r} = \bm{r}_B - \bm{r}_M$ and $\bm{v} = \bm{v}_B$.

    \item \ul{Moon-Centered Earth Orbit Frame} (\acrshort{eof}): A basis $\hat{\bm{X}}$-$\hat{\bm{Y}}$-$\hat{\bm{Z}}$ is defined with $\hat{\bm{Z}} = \hat{\bm{z}}$ and $\hat{\bm{X}}$-$\hat{\bm{Y}}$ oriented such that the rotating basis $\hat{\bm{x}}$-$\hat{\bm{y}}$ revolves about $\hat{\bm{Z}}$ at the nd rate of unity relative to $\hat{\bm{X}}$-$\hat{\bm{Y}}$. The origin is at the Moon. The nd spacecraft position in this frame relates to the rotating-frame state via,
    \begin{align}
        \label{eq:eof_mrf}
        \bm{R} = \bm{C}_{\hat{\bm{z}}} \bm{r} = X\hat{\bm{X}} + Y\hat{\bm{Y}} + Z\hat{\bm{Z}},
        \quad
        \bm{C}_{\hat{\bm{z}}}(t) =
        \begin{bmatrix}
        \cos (t + t_{0}) & -\sin (t + t_{0}) & 0 \\
        \sin (t + t_{0}) & \cos (t + t_{0}) & 0 \\
        0 & 0 & 1
        \end{bmatrix},
    \end{align}
    The angle $t_0$ specifies the orientation of $\hat{\bm{X}}$ within the $xy$ plane at a reference epoch. This parameter may be selected to align $\hat{\bm{X}}$ with a physically defined direction, such as the intersection of the lunar equatorial plane and the Earth orbit plane adopted in \citet{ely2005stable}; however, any value is admissible without fundamentally altering the dynamics. The velocity within the frame is evaluated as $\bm{V} = d\bm{R}/dt$ and connects to the \acrshort{mrf} velocity through the derivative of $\bm{R}$ in Eq.~\eqref{eq:eof_mrf}. Within the \acrshort{dadm} and \acrshort{cr3bp}, the \acrshort{eof} is inertial. Within the \acrshort{hfem}, the basis vectors evolve with the ephemerides and the frame is not inertially fixed; integrations are, therefore, not simulated in the \acrshort{eof}, and the \acrshort{eof} is reserved for orbital-element specification at each epoch.

    \item \ul{J2000 Moon-Centered Inertial Frame} (\acrshort{mci}): This inertial frame is centered at the Moon with its unit vectors $\hat{\bm{X}}_J$-$\hat{\bm{Y}}_J$-$\hat{\bm{Z}}_J$ aligned with the J2000 ecliptic and equinox, i.e., the \texttt{ECLIPJ2000} frame. The ecliptic, rather than the Earth mean equator, is adopted as the fundamental plane since the lunar ascending node regresses at a nearly uniform rate only when referred to the ecliptic, a property leveraged in Sec.~\ref{subsec:freq-across-models}. The nd spacecraft position vector in this frame is~
    \begin{equation}
        \label{eq:mci_mrf} \bm{R}_J = \bm{C} \bm{r} = X_J\hat{\bm{X}}_J + Y_J\hat{\bm{Y}}_J+ Z_J\hat{\bm{Z}}_J, \quad \bm{C}(t) = \begin{bmatrix} \hat{\bm{x}} & \hat{\bm{y}} & \hat{\bm{z}} \end{bmatrix},
    \end{equation}
    The matrix $\bm{C}$ is the direction cosine matrix whose columns are the \acrshort{mrf} unit vectors expressed in the \acrshort{mci}, evaluated at each epoch using the JPL DE440 ephemerides~\cite{park2021jpl}. The velocity in the frame is evaluated as $\bm{V}_J = d\bm{R}_J/dt$.
\end{itemize}
Classical Keplerian elements, namely the semi-major axis $a$, eccentricity $e$, inclination $i$, argument of periapsis $\omega$, right ascension of the ascending node $\Omega$, and mean anomaly $M$, are employed throughout this work and are defined with respect to the \acrshort{eof}, unless otherwise noted. Standard two-body transformations~\cite{vallado2001fundamentals} allow shifts between osculating elements and the Cartesian state $(\bm{R},\, \bm{V})$ in the \acrshort{eof}. Equation~\eqref{eq:eof_mrf} supplies the corresponding states for trajectory analysis in the \acrshort{mrf} and \acrshort{brf} within the \acrshort{dadm} and \acrshort{cr3bp}. 

Transformation from the \acrshort{eof} to the \acrshort{mci} requires more caution in terms of its velocity components due to the derivatives of the elements of $\bm{C}$. Incorporating the $\dot{\bm{C}}$ terms in Eqs.~\eqref{eq:eof_mrf}-\eqref{eq:mci_mrf}, although mathematically consistent, causes changes in shape-related Keplerian elements (e.g., $a, e$) across the \acrshort{eof} and \acrshort{mci}. To eliminate this drift, a direct transformation between the two frames is introduced at a reference epoch, i.e.,
\begin{align}
    \label{eq:eof_mci_snapshot}
    \bm{R}_J = \bm{C}_{\!\text{eff}}\,\bm{R},
    \qquad
    \bm{V}_J = \bm{C}_{\!\text{eff}}\,\bm{V},
    \qquad
    \bm{C}_{\!\text{eff}} \equiv \bm{C}\,\bm{C}_{\hat{\bm{z}}}^T(t_0),
\end{align}
evaluated at the selected epoch via the JPL DE440 ephemerides. The orthogonality of $\bm{C}_{\!\text{eff}}$ preserves $|\bm{R}|$, $|\bm{V}|$, and $\bm{R}\cdot\bm{V}$, such that the Keplerian shape parameters $a$, $e$ are invariant between the two frames. This transformation effectively assumes the \acrshort{eof} reference frame to be frozen at its epoch and is consistent with previous work \cite{ely2005stable}. While this nuance from the transformation is noted for contextual clarity, the primary findings in the current analysis are valid regardless of the specific transformation scheme between \acrshort{eof} and \acrshort{mci}. 

\subsection{Dynamical Models}

Three dynamical models spanning different fidelity levels are employed: (1) The \acrshort{dadm}, (2) \acrshort{cr3bp}, and (3) \acrshort{hfem}.

\subsubsection{\acrshort{dadm}}

The \acrshort{dadm} governs the mean orbital element evolution as viewed in the \acrshort{eof} under the following assumptions: (1) the Earth is the sole perturbing body, (2) both the Earth and Moon are point masses in mutual circular orbits, and (3) the perturbing potential of the Earth is expanded in Legendre polynomials and truncated at second order~\cite{longuski2022introduction}. Averaging the truncated potential over both the spacecraft Keplerian period around the Moon and the Earth orbital period around the Moon yields the doubly-averaged perturbing potential\footnote{Consistent with Broucke~\cite{broucke2003long}, the full factor $(1-\mu)$ is retained throughout the expressions that follow, whereas the numerical evaluations reported in this work adopt the approximation $1-\mu \approx 1$ employed elsewhere~\cite{folta2006lunar}.},
\begin{align}
    \label{eq:perturbing_potential} U_{E} = \frac{1}{32}(1-\mu)n_E^2 a^2 \left[ (1+3\cos 2i)(2 + 3e^2) + 30 e^2 \sin^2 i \cos 2\omega \right],
\end{align}
The dimensional mean angular rate for the Earth in its orbit with respect to the Moon is $n_E = 1/t_*$. The orbital elements in Eq.~\eqref{eq:perturbing_potential} are mean elements defined within the \acrshort{eof}. Lagrange's planetary equations yield the following evolution rates with respect to the nd time $t$,
\begin{align}
    \label{eq:dadt} & \dot{a} = 0, \\
    \label{eq:dedt}& \dot{e} = \frac{15}{8}\frac{(1-\mu)n_E}{n}\,e\sqrt{1-e^2}\,\sin^2 i\,\sin 2\omega, \\
    \label{eq:didt}& \dot{i} = - \frac{15}{16}\frac{(1-\mu)n_E}{n}\,\frac{e^2}{\sqrt{1-e^2}}\,\sin 2i\, \sin 2\omega, \\
    \label{eq:domegadt}& \dot{\omega} = \frac{3}{16}\frac{(1-\mu)n_E}{n}\,\frac{1}{\sqrt{1-e^2}} \left[ (3 + 2e^2 + 5\cos 2i) + 5(1-2e^2 - \cos 2i) \cos 2\omega \right], \\
   \label{eq:dOmegadt}   & \dot{\Omega} = \frac{3}{8}\frac{(1-\mu)n_E}{n}\,\frac{1}{\sqrt{1-e^2}}\, (5e^2 \cos 2\omega - 3e^2 - 2) \cos i, \\
   \label{eq:dthetadt}  & \dot{M} = nt_* = \sqrt{\frac{\tilde{\mu}_M}{a^3}}\,t_*,
\end{align}
The parameter $n$ denotes the dimensional mean motion of the spacecraft with respect to the Moon as determined by $a$. In addition to $a$, the system possesses two additional integrals of motion \cite{broucke2003long},
\begin{align}
    \label{eq:C1_def} C_1 &= (1 - e^2)\cos^2 i, \\
    \label{eq:C2_def} C_2 &= e^2\!\left(\tfrac{2}{5} - \sin^2 i\,\sin^2\omega\right),
\end{align}
These integrals are related to the $\hat{\bm{z}}$-directional angular momentum vectors and the total energy for the spacecraft. Since the evolution of $e$, $i$, and $\omega$ depends only on these two integrals (with $n=n(a)$ fixed from Eq.~\eqref{eq:dadt}), the \acrshort{dadm} reduces to a single-degree-of-freedom system for given $(C_1, C_2)$. The resulting dynamics in $(e, i, \omega)$ are hereafter denoted as the \emph{reduced dynamics} and are visualized in the $(e\cos\omega,\, e\sin\omega)$ phase portrait. The equilibria in this reduced system satisfy $\dot{e} = \dot{i} =  \dot{\omega} = 0$, yielding the frozen-orbit condition,\footnote{Another equilibrium type, denoted ``circular,'' exists at $e = 0$~\cite{nie2018lunar}, but remains out of scope.}
\begin{align}
    \label{eq:frozen} \omega = \frac{\pi}{2}, \qquad e = \sqrt{1 - \tfrac{5}{3}\cos^2 i}.
\end{align}
These equilibria are elliptic, and nearby trajectories librate around them in the $(e\cos\omega,\, e\sin\omega)$ plane (a detailed interpretation of the phase portrait is deferred to Sec.~\ref{subsec:freq-modes}). An \acrshort{elfo} in the current analysis refers to either a frozen equilibrium or a librating trajectory in its vicinity.

\subsubsection{\acrshort{cr3bp}}

The \acrshort{cr3bp} is formulated in the \acrshort{brf}. The equations of motion are expressed as
\begin{align}
    \label{eq:cr3bp}
    \dot{\bm{v}}_B = -2\hat{\bm{z}} \times \bm{v}_B + \nabla U_C,
\end{align}
The pseudo-potential is $U_C = \tfrac{1}{2}(x_B^2 + y_B^2) + (1-\mu)/d + \mu/r$, with $d = \|\bm{r}_B - \bm{r}_E\|_2$ and $r = \|\bm{r}_B - \bm{r}_M\|_2$, and $\times$ denotes the vector cross product. As in the \acrshort{dadm}, the Earth and Moon are modeled as point masses on mutually circular orbits about their barycenter. In contrast, the gravitational influence of the Earth is retained in full (but as a point mass), without Legendre truncation or averaging.

\subsubsection{\acrshort{hfem}}

The \acrshort{hfem} incorporates point-mass gravity from three bodies: the Sun ($S$), Earth ($E$), and Moon ($M$). Their instantaneous positions are supplied from the JPL DE440 ephemerides~\cite{park2021jpl}. With the Moon as the central body, the equations of motion in the \acrshort{mci} are
\begin{align}
    \label{eq:hfem} \dot{\bm{V}}_J = -\frac{\mu}{R_J^3}\bm{R}_J + \sum_{j \in \{E,\, S\}} \left[ \frac{\mu_j}{R_{jM}^3}\bm{R}_{jM} - \frac{\mu_j}{R_{jc}^3} \bm{R}_{jc} \right],
\end{align}
The vectors $\bm{R}_{jc}$ and $\bm{R}_{jM}$ are directed from body $j$ to the spacecraft and the Moon, respectively, and $\mu_j = \tilde{\mu}_j/(\tilde{\mu}_E + \tilde{\mu}_M)$ is the nd gravitational parameter of body $j$. Although the \acrshort{hfem} may be extended to include additional perturbations such as solar radiation pressure and the irregular gravitational field of the Moon, the current scope is limited to point-mass gravity to isolate the largest perturbation for the \acrshort{elfos}, i.e., the realistic motion of the Earth with respect to the Moon \cite{ely2005stable,nie2018lunar}.

\section{ELFOs as Invariant Tori: Frequency Structures}
\label{sec:freq-structure}
The frequency structure for the \acrshort{elfos} underpins the constellation design and analysis framework developed in subsequent sections. This section first reviews the intrinsic frequency modes and their coupling to the orbital elements within the baseline model (\acrshort{dadm}), and then examines the evolution of the frequency content across the other two models, i.e., the \acrshort{cr3bp} and \acrshort{hfem}. Finally, lunar obliquity and the satellite-user geometry are examined in the frequency domain.

\subsection{Three Frequency Modes within the DADM}
\label{subsec:freq-modes}
\begin{figure}[h!]
    \centering
    \begin{subfigure}[b]{0.48\textwidth}
        \centering
        \includegraphics[width=\linewidth]{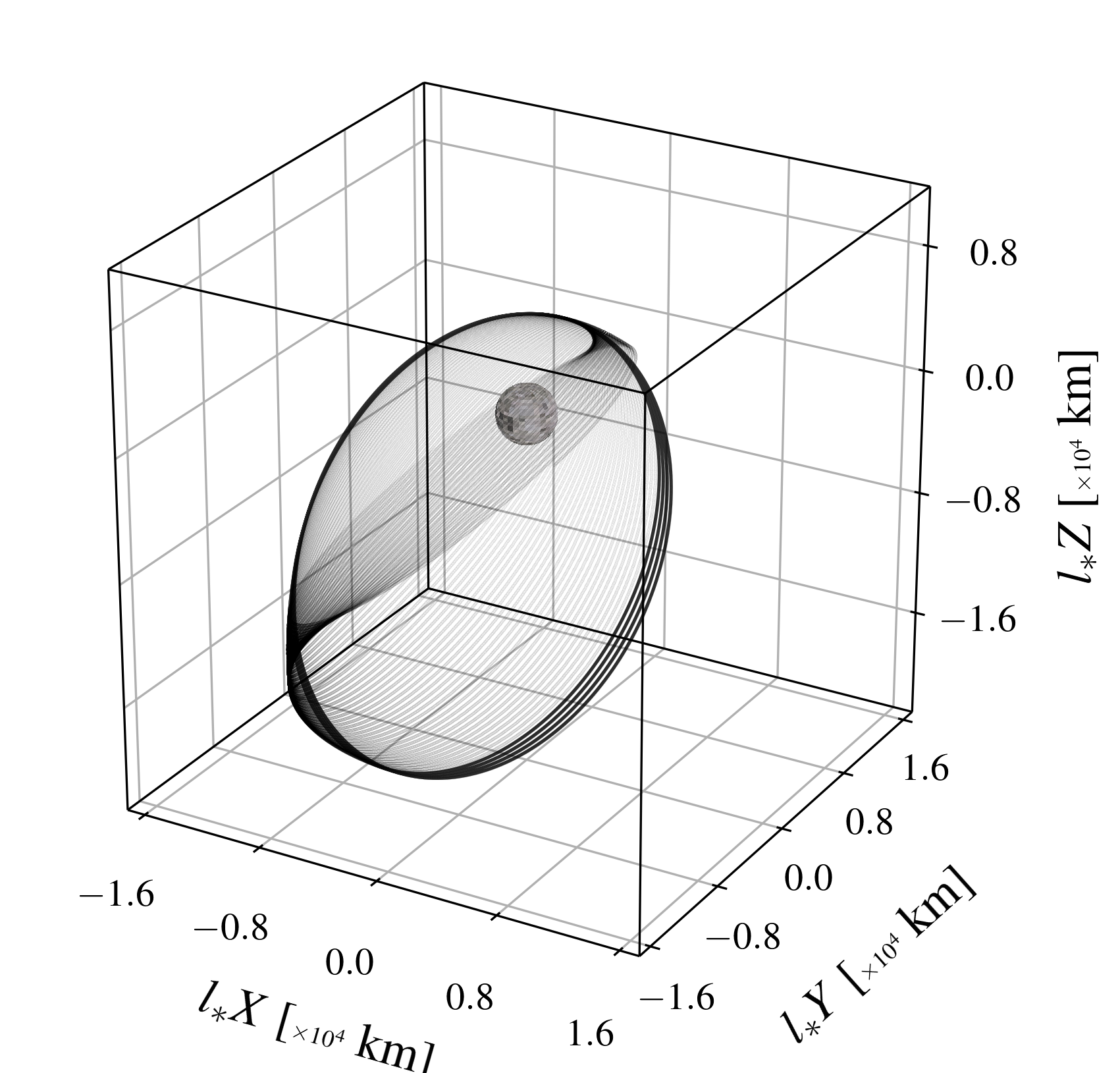}
        \caption{An \acrshort{elfo} at $A_L = 0$ in the \acrshort{eof}.}\label{fig:freq_inertial}
    \end{subfigure}
    \hfill
    \begin{subfigure}[b]{0.48\textwidth}
        \centering
        \includegraphics[width=\linewidth]{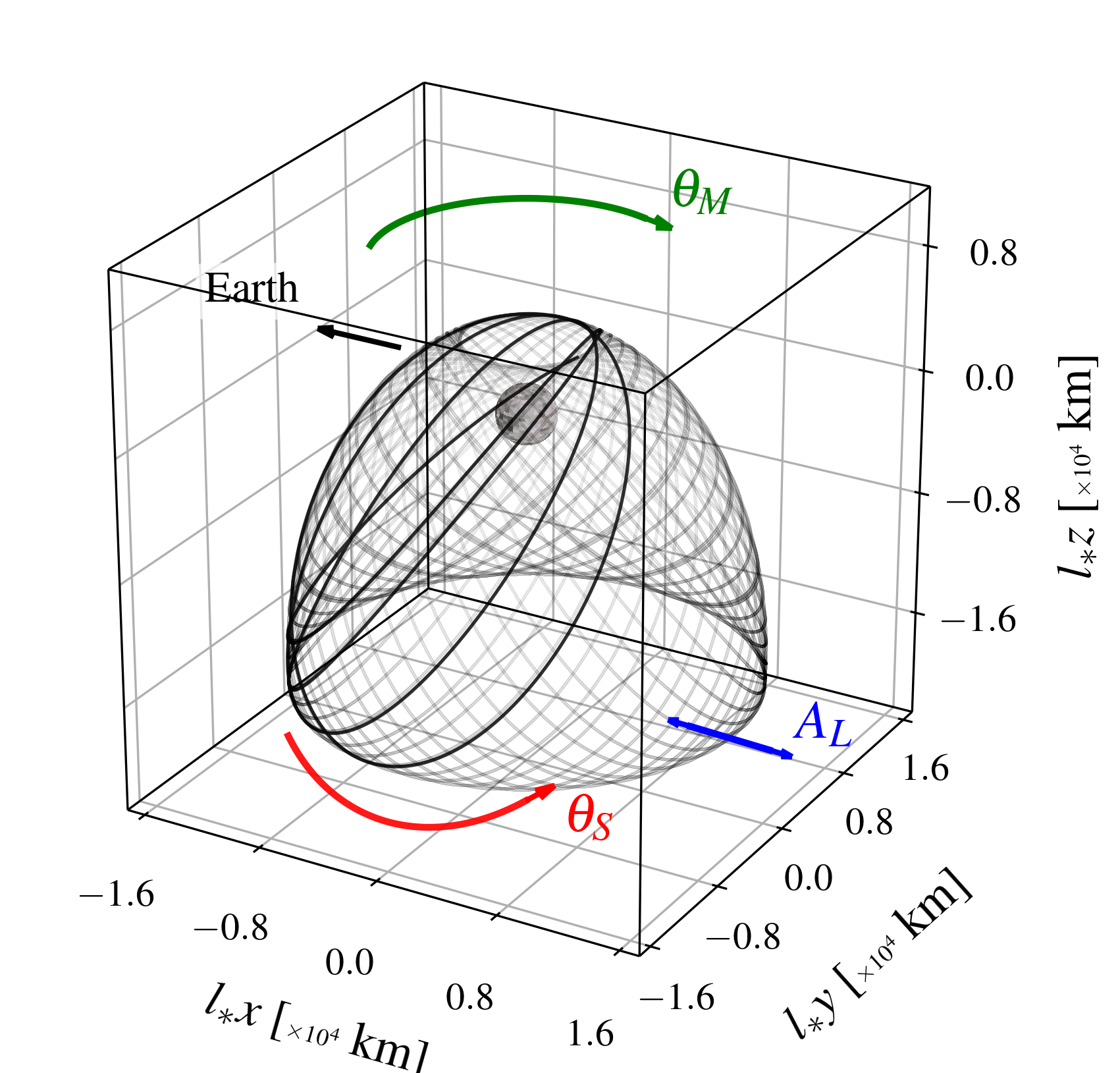}
        \caption{An \acrshort{elfo} at $A_L = 0$ in the \acrshort{mrf}.}\label{fig:freq_structure}
    \end{subfigure}
    \\[6pt]
    \begin{subfigure}[b]{0.48\textwidth}
        \centering
        \includegraphics[width=\linewidth]{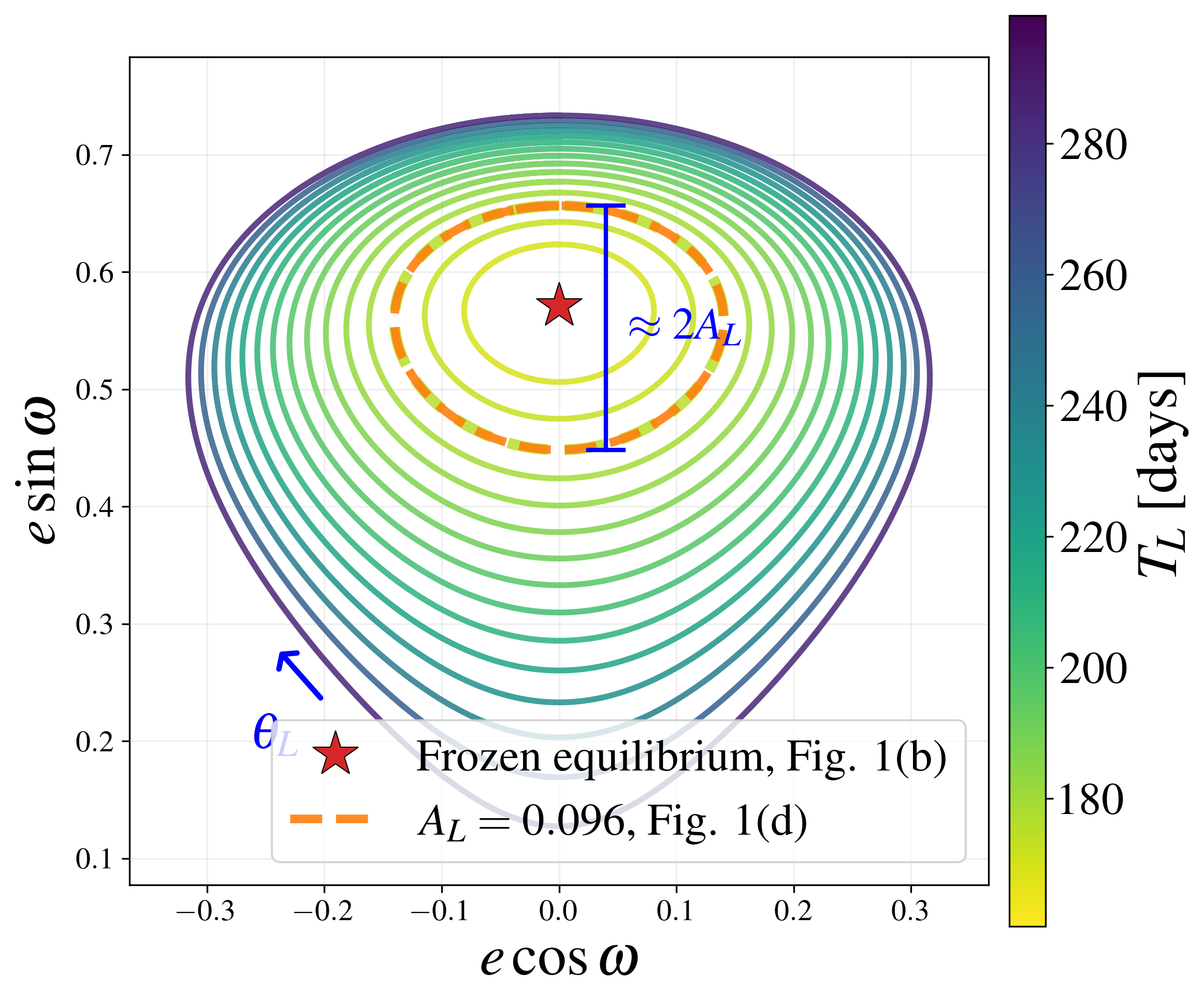}
        \caption{Phase portrait in the $(e\cos\omega,\, e\sin\omega)$ plane.}\label{fig:e_w_plane}
    \end{subfigure}
    \hfill
    \begin{subfigure}[b]{0.48\textwidth}
        \centering
        \includegraphics[width=\linewidth]{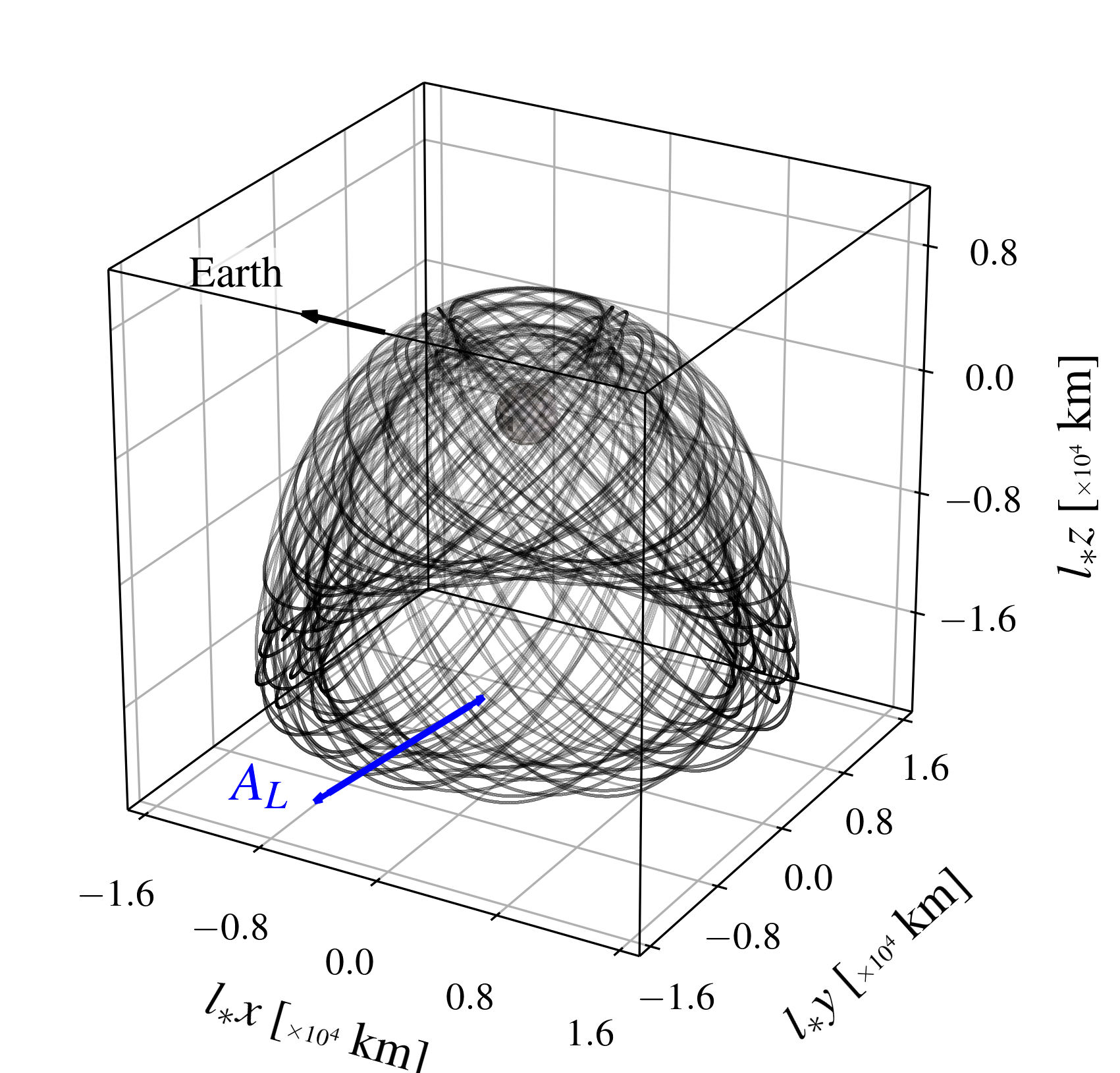}
        \caption{An \acrshort{elfo} at $A_L = 0.096$ in the \acrshort{mrf}.}\label{fig:freq_nonfrozen}
    \end{subfigure}
    \caption{Frequency structures of representative \acrshort{elfos} ($a = 14{,}200$~km, $i = 50.5^\circ$).}\label{fig:freq_graphics}
\end{figure}

The \acrshort{dadm} globally admits a quasi-periodic solution structure due to its integrability. With three integrals of motion ($a, C_1, C_2$), the solutions evolve on tori of dimension three at most. This dimensionality is readily understood through the reduced dynamics (Sec.~\ref{sec:prelim}). The frozen equilibria correspond to fixed points in the reduced phase space with $e= e(i)$ and $\omega =\pm \pi/2$. As \acrshort{lsp} missions motivate the present application, $\omega = \pi/2$ is assumed in the current analysis. One such frozen configuration is plotted in Fig.~\ref{fig:e_w_plane} (red star) as a fixed point within the reduced phase space. These frozen equilibria and their associated libration trajectories, collectively denoted \acrshort{elfos} in the current analysis, correspond to two-dimensional (2D) and three-dimensional (3D) invariant tori within the full phase space, respectively~\cite{park2026bridging}. Accordingly, three frequency and phase-angle pairs, denoted $(\bm{\nu}, \bm{\theta}) = (\nu_S, \nu_M, \nu_L, \theta_S, \theta_M, \theta_L)$, fully determine an \acrshort{elfo} and the specific location along it~\cite{park2026bridging}. Each frequency-phase pair admits a physical interpretation, summarized in Table~\ref{tab:freq_summary}; the table lists the physical origin of each mode, its typical period, and representative scalar signals where the corresponding mode is clearly detected.

The three frequency modes are illustrated through a representative \acrshort{elfo} with $a = 14{,}200$~km and $i = 50.5^\circ$ at the frozen equilibrium condition ($e \approx 0.5707$, $\omega = \pi/2$), as specified in Table~\ref{tab:baseline_orbit}. The Keplerian orbital period is $T_S \approx 1.76$~days, corresponding to a mean motion $\nu_S = 2\pi / T_S \approx 15.55$~rad/nd with $T_S$ expressed in nd time. A three-month trajectory in the \acrshort{eof} (Fig.~\ref{fig:freq_inertial}) reveals a fixed orbital geometry with a slowly precessing $\Omega$ ($\dot{\Omega} <0$ from Eq.~\eqref{eq:dOmegadt}). The same trajectory viewed in the \acrshort{mrf} (Fig.~\ref{fig:freq_structure}) separates the short-period orbital revolution from the medium-period nodal precession more clearly: the first four revolutions (black traces) are distinguishable, and the directions of $\theta_S$ (red) and $\theta_M$ (green) are annotated. The short-period phase $\theta_S$ is approximately equal to the osculating mean anomaly $M$, and $\nu_S$ is analogous to the mean motion. Although $\theta_S = \nu_S t$ is strictly linear by definition, the osculating $M$ varies nonlinearly due to short-period perturbations within any unaveraged models, e.g., the \acrshort{cr3bp}; nevertheless, $M$ serves as a practical observable for extracting $\theta_S$. Similarly, the medium-period phase satisfies $\theta_M \approx -\Omega_R$, where $\Omega_R$ is the right ascension of the ascending node in \acrshort{mrf}; the negative sign ensures a positive $\nu_M$. For the sample orbit, the period for $\nu_M$ is $T_M \approx 25$~days (Table~\ref{tab:baseline_orbit}). Although $\Omega_R$ is an osculating quantity, $\theta_M = \nu_M t$ supplies its averaged counterpart beyond the \acrshort{dadm}.

\begin{figure}[h!]
    \centering
    \begin{subfigure}[b]{0.48\textwidth}
        \centering
        \includegraphics[width=\linewidth]{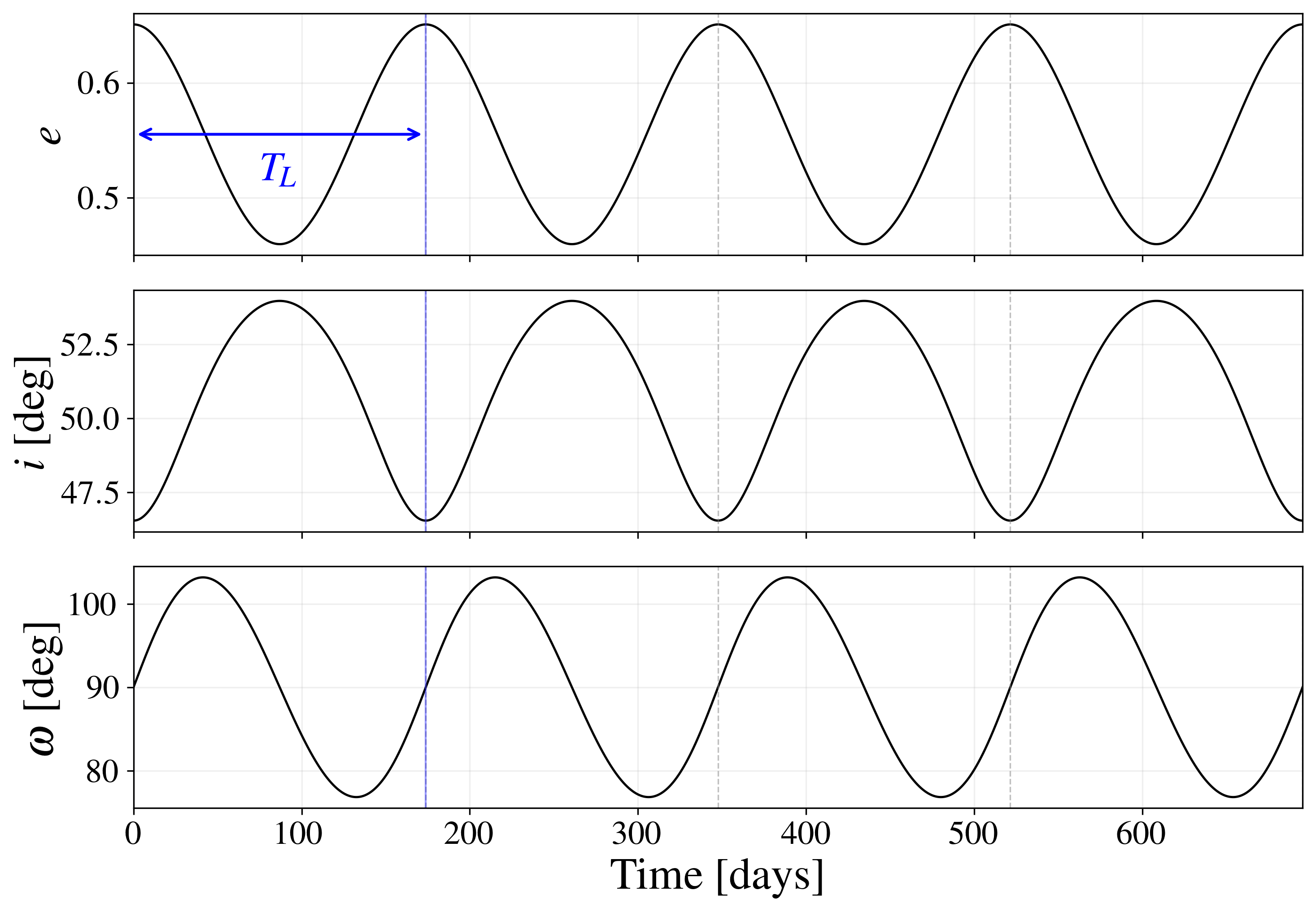}
        \caption{Time history of $e$, $i$, $\omega$ within the \acrshort{eof}.}\label{fig:coupling_time}
    \end{subfigure}
    \hfill
    \begin{subfigure}[b]{0.48\textwidth}
        \centering
        \includegraphics[width=\linewidth]{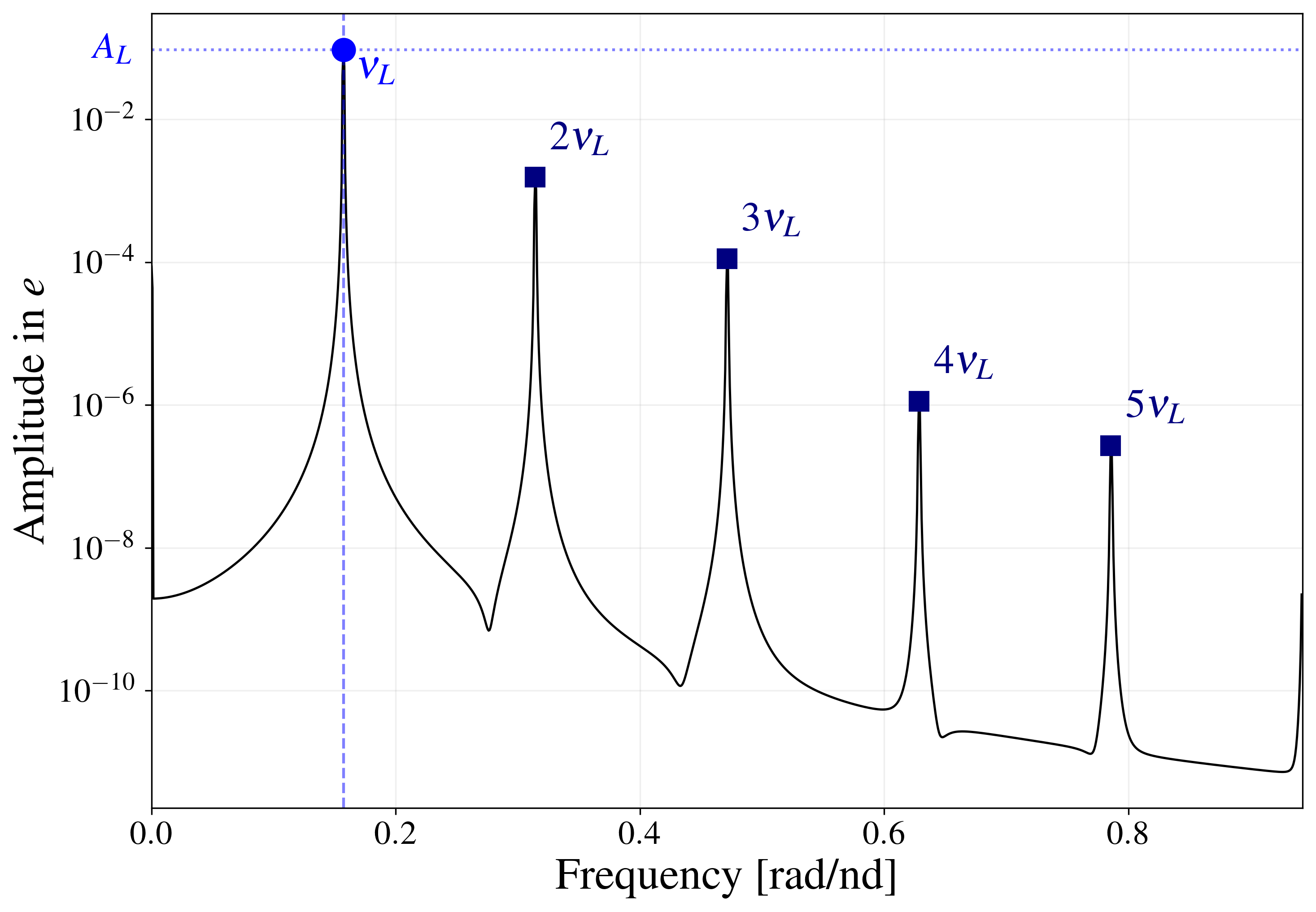}
        \caption{Spectral decomposition of $e$.}\label{fig:coupling_fft}
    \end{subfigure}
    \caption{Long-period coupling of $(e, i, \omega)$ within the \acrshort{dadm} ($A_L \neq 0$ configuration from Table~\ref{tab:baseline_orbit}).}\label{fig:coupling}
\end{figure}

The third frequency component, the long-period pair ($\nu_L, \theta_L$), is illustrated within the reduced domain (Fig.~\ref{fig:e_w_plane}). Since the frozen equilibrium corresponds to a fixed point within the reduced space, $\nu_L$ is degenerate at the equilibrium. Each closed contour in Fig.~\ref{fig:e_w_plane} represents a different value of $C_2$ at the same $C_1$, colored by its corresponding long-period $T_L = 2\pi/\nu_L$\footnote{For a given pair of integrals $(C_1, C_2)$, the libration period $T_L$ is evaluated as a nonsingular quadrature over the eccentricity range $[e_{\min}, e_{\max}]$, motivated by Russell and Brinckerhoff \cite{russell2009circulating}; see Appendix~\ref{app:lp-quadrature} for the derivation.}. Orbits closer to the equilibrium are associated with smaller $T_L$, whereas those farther away span wider excursions in $(e, i, \omega)$ space with longer libration periods. Due to the degeneracy of $\nu_L$ at the fixed point, $A_L = A_L(\nu_L)$ defines a more tractable quantity that measures deviation from the frozen condition, with $A_L = 0$ corresponding to the equilibrium. In the current investigation, $A_L$ is defined as the first harmonic amplitude for the eccentricity, $e$. For illustration, the dashed orange curve in Fig.~\ref{fig:e_w_plane} represents an $A_L \neq 0$ configuration with the initial condition supplied in Table~\ref{tab:baseline_orbit}. The corresponding trajectory in \acrshort{mrf} is illustrated in Fig.~\ref{fig:freq_nonfrozen} for a year. As this configuration does not satisfy the frozen condition, $e$, $i$, and $\omega$ oscillate with the common frequency $\nu_L$ in a coupled fashion, as illustrated in Fig.~\ref{fig:coupling_time}. The amplitude $A_L$ is then measured from the spectrum of $e$ as the amplitude at $\nu_L$\footnote{While the other two variables $i, \omega$ may be leveraged to define $A_L$ as well, $e$ is most easily accessible as it is independent of the frame, i.e., $e$ in \acrshort{mci} is identical to $e$ within the \acrshort{eof} \cite{ceresoli2025design}.} in Fig.~\ref{fig:coupling_fft}. Since the harmonics decay in amplitude, $A_L$ is approximately supplied with the maximum and minimum value of $e$ as $A_L \approx (e_{\max} - e_{\min})/2$ within the reduced phase space, with the approximation becoming exact as $A_L \rightarrow 0$. Within the \acrshort{mrf}, $A_L$ manifests as oscillations in the orbital plane, together with the coupled behavior in $e$ and $\omega$, as depicted in Fig.~\ref{fig:freq_nonfrozen} in contrast to the frozen equilibrium case (Fig.~\ref{fig:freq_structure}).

Since the amplitude $A_L$ quantifies the deviation from the often desired frozen configuration rather than a frequency, it serves as a key design variable alongside $(\nu_S, \nu_M)$. Existing constellation design often imposes the frozen condition \emph{a priori} in lower-fidelity analysis \cite{brack2025,ceresoli2025design, zanotti2024high}, while nonzero $A_L$ is absorbed implicitly into high-dimensional numerical searches within realistic dynamical models; in either case, the long-period mode is seldom treated as an explicit design variable.

\begin{table}[h!]
\centering
\caption{Dominant frequency-phase components and representative scalar signals.}\label{tab:freq_summary}
\setlength{\tabcolsep}{5pt}
\renewcommand{\arraystretch}{1.15}
\footnotesize
\begin{tabularx}{\textwidth}{@{} l l l l l X @{}}
\toprule
Period & Component & Physical meaning & Typical period & Fig.~\ref{fig:freq_structure} color & Representative scalars \\
\midrule
Short & $(\nu_S,\theta_S)$ & Orbital revolution around the Moon & $\sim$1~day &
Red &
$M$; $z$ \\

Medium & $(\nu_M,\theta_M)$ & Nodal precession in rotating frame & $\sim$1~month &
Green &
$-\Omega_R$; $x, y$ \\

Long & $(\nu_L,\theta_L)$ & Libration about frozen equilibrium & months to years &
Blue &
$e$; $i$; $\omega$ \\

\bottomrule
\end{tabularx}
\end{table}

\begin{table}[h!]
\centering
\caption{Initial conditions for the representative \acrshort{elfos} in Fig.~\ref{fig:freq_graphics} within the \acrshort{dadm}. Frequencies in rad/nd.}\label{tab:baseline_orbit}
\setlength{\tabcolsep}{5pt}
\renewcommand{\arraystretch}{1.2}
\footnotesize
\begin{tabular}{@{} l c c c c c c c c c c @{}}
\toprule
& $a$ [km] & $e$ & $i$ [deg] & $\omega$ [deg] & $C_1$ & $C_2$ & $\nu_S$ ($T_S$) & $\nu_M$ ($T_M$) & $\nu_L$ ($T_L$) & $A_L$ \\
\midrule
$A_L = 0$ (Fig.~\ref{fig:freq_structure})
& 14{,}200 & 0.5707 & 50.5 & 90.0 & 0.2728 & $-0.0636$
& \multirow{2}{*}{\makecell{15.55 \\ (1.76~d)}}
& \makecell{1.086 \\ (25.2~d)} & - & 0 \\
$A_L \neq 0$ (Fig.~\ref{fig:freq_nonfrozen})
& 14{,}200 & 0.6507 & 46.5 & 90.0 & 0.2728 & $-0.0537$
&
& \makecell{1.083 \\ (25.2~d)} & \makecell{0.157 \\ (173.8~d)} & 0.096 \\
\bottomrule
\end{tabular}
\end{table}

\subsection{Frequency Content across Dynamical Models}
\label{subsec:freq-across-models}

The intrinsic \acrshort{dadm} frequency structure is examined under progressively more realistic dynamical models. A hierarchical progression from the \acrshort{dadm} to the \acrshort{cr3bp} and then to the \acrshort{hfem} suggests empirical robustness of the quasi-periodic structure for the investigated \acrshort{elfo} while identifying additional frequency components introduced by perturbations.

\subsubsection{\acrshort{cr3bp}} The \acrshort{cr3bp} retains the autonomous, i.e., time-independent, character of the system while introducing unaveraged gravitational forcing from the Earth. The model is generally non-integrable, and its phase space may admit hyperbolic manifolds beyond invariant tori, especially as semi-major axis ($a$) increases and the orbit is subject to larger perturbations from the Earth. The current work assumes the existence of dense invariant tori within the \acrshort{cr3bp}; resonances and diffusive mechanisms may erode tori at sufficiently large semi-major axis ($a$), but the configurations here exhibit numerically stable quasi-periodic behavior over the propagation horizons employed. A detailed characterization of the torus boundaries is deferred to future work.

\begin{figure}[h!]
    \centering
    \begin{subfigure}[b]{0.48\textwidth}
        \centering
        \includegraphics[width=\linewidth]{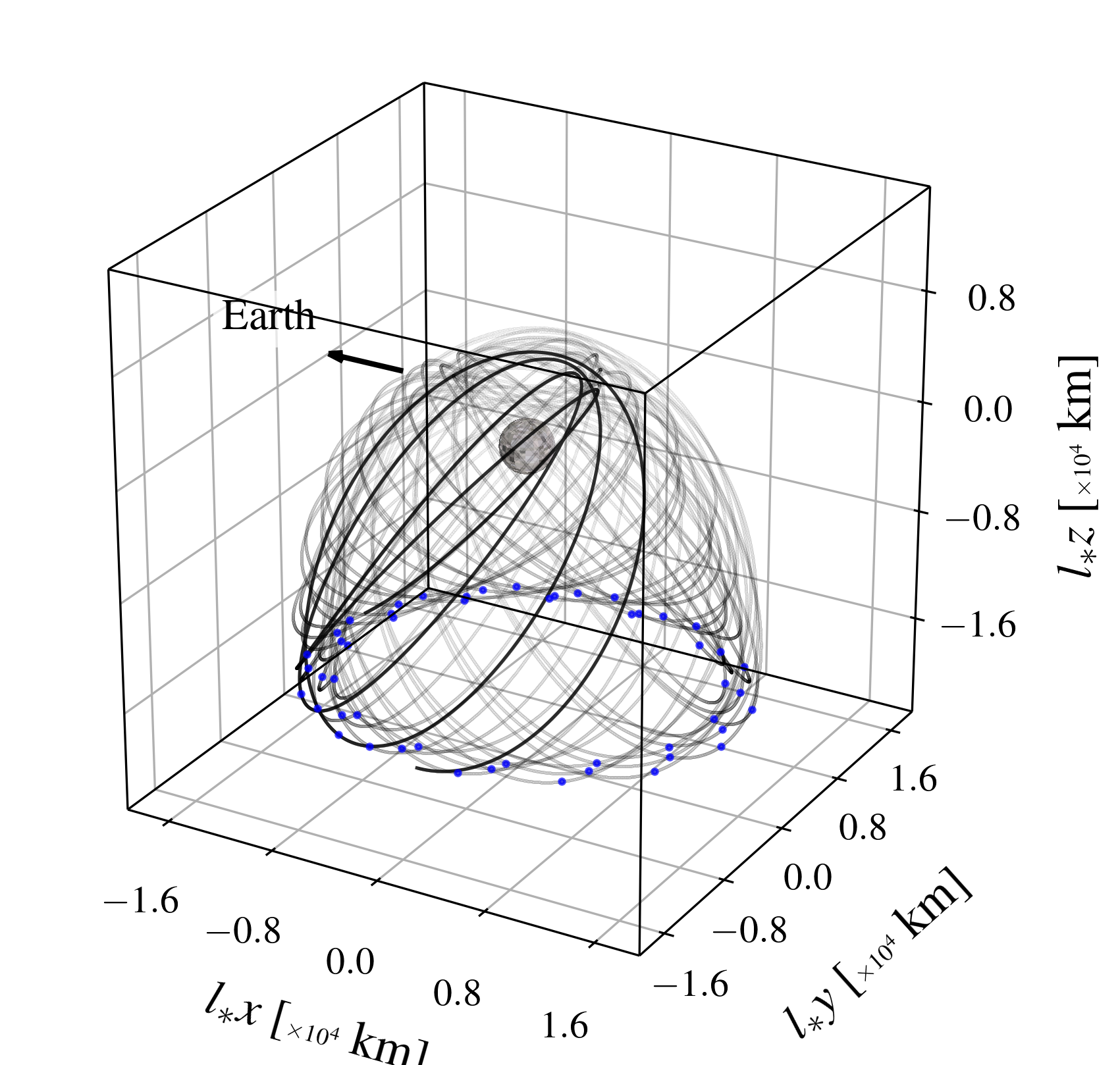}
        \caption{Three-month trajectory (\acrshort{mrf}).}\label{fig:cr3bp_geom}
    \end{subfigure}
    \hfill
    \begin{subfigure}[b]{0.48\textwidth}
        \centering
        \includegraphics[width=\linewidth]{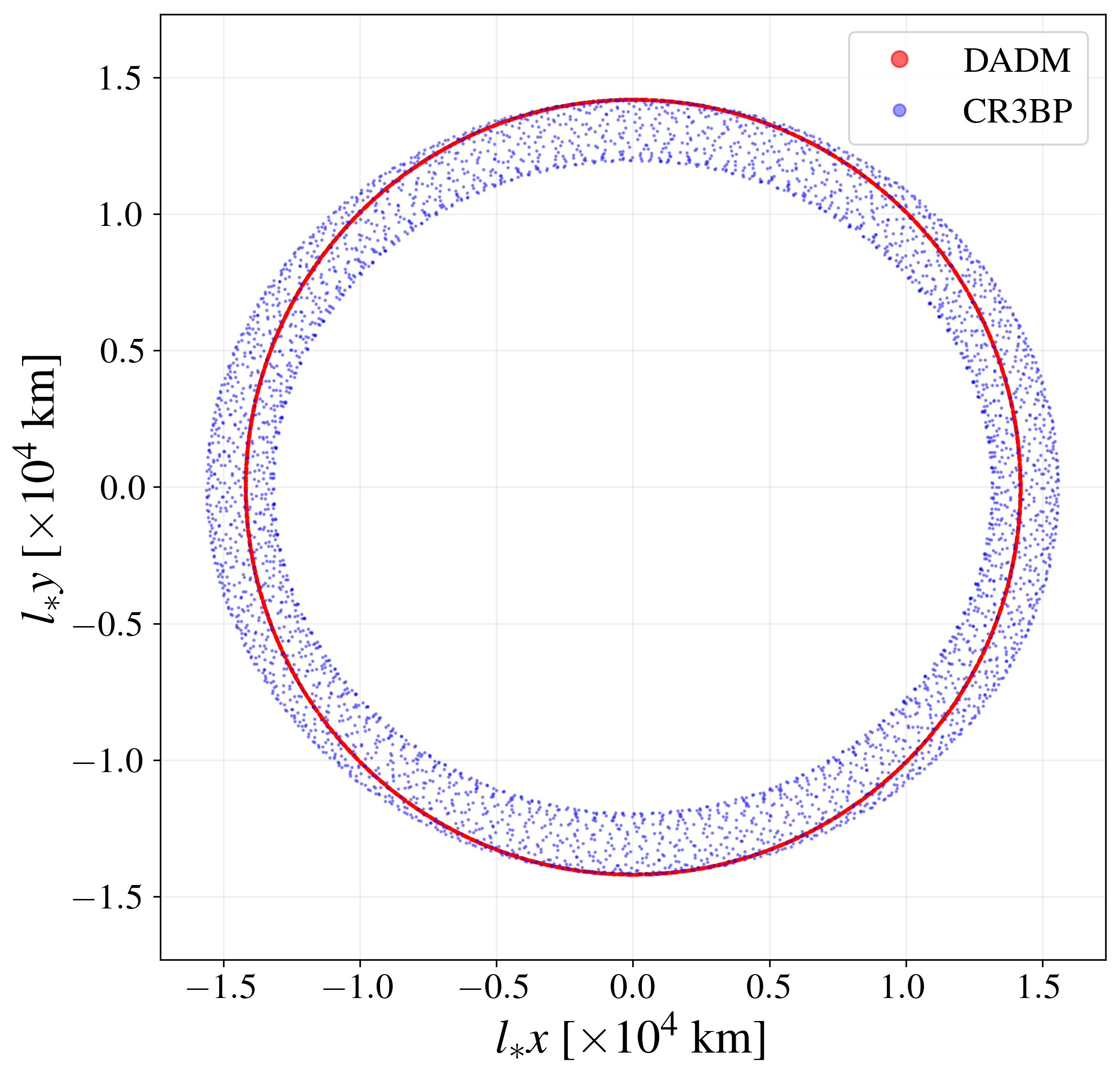}
        \caption{Apolune states (20 years), \acrshort{dadm} vs \acrshort{cr3bp}.}\label{fig:strobo_cr3bp}
    \end{subfigure}
    \\[6pt]
    \begin{subfigure}[b]{0.48\textwidth}
        \centering
        \includegraphics[width=\linewidth]{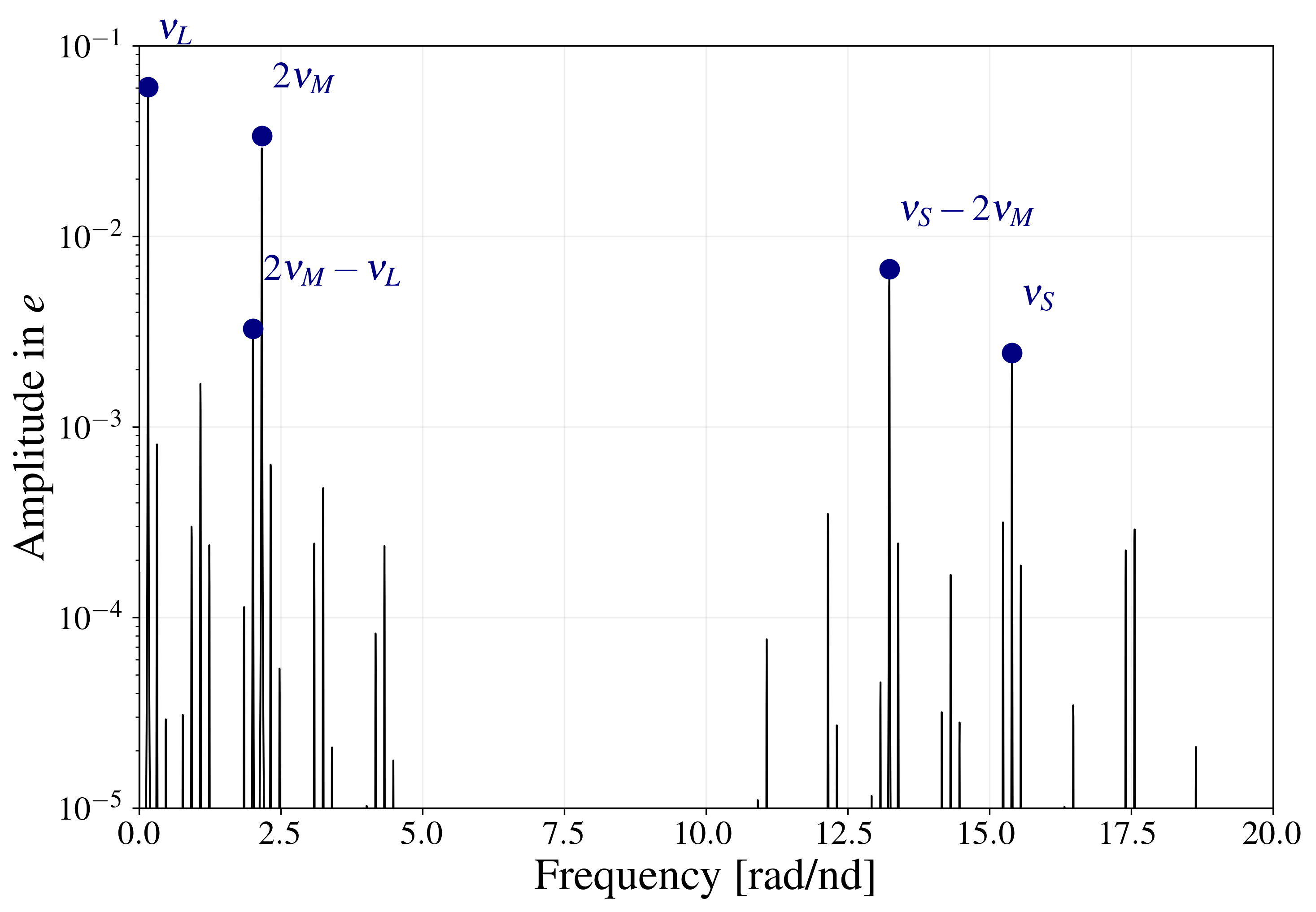}
        \caption{Spectral decomposition of $e$ with the five dominant frequencies annotated.}\label{fig:cr3bp_fft_e}
    \end{subfigure}
    \hfill
    \begin{subfigure}[b]{0.48\textwidth}
        \centering
        \includegraphics[width=\linewidth]{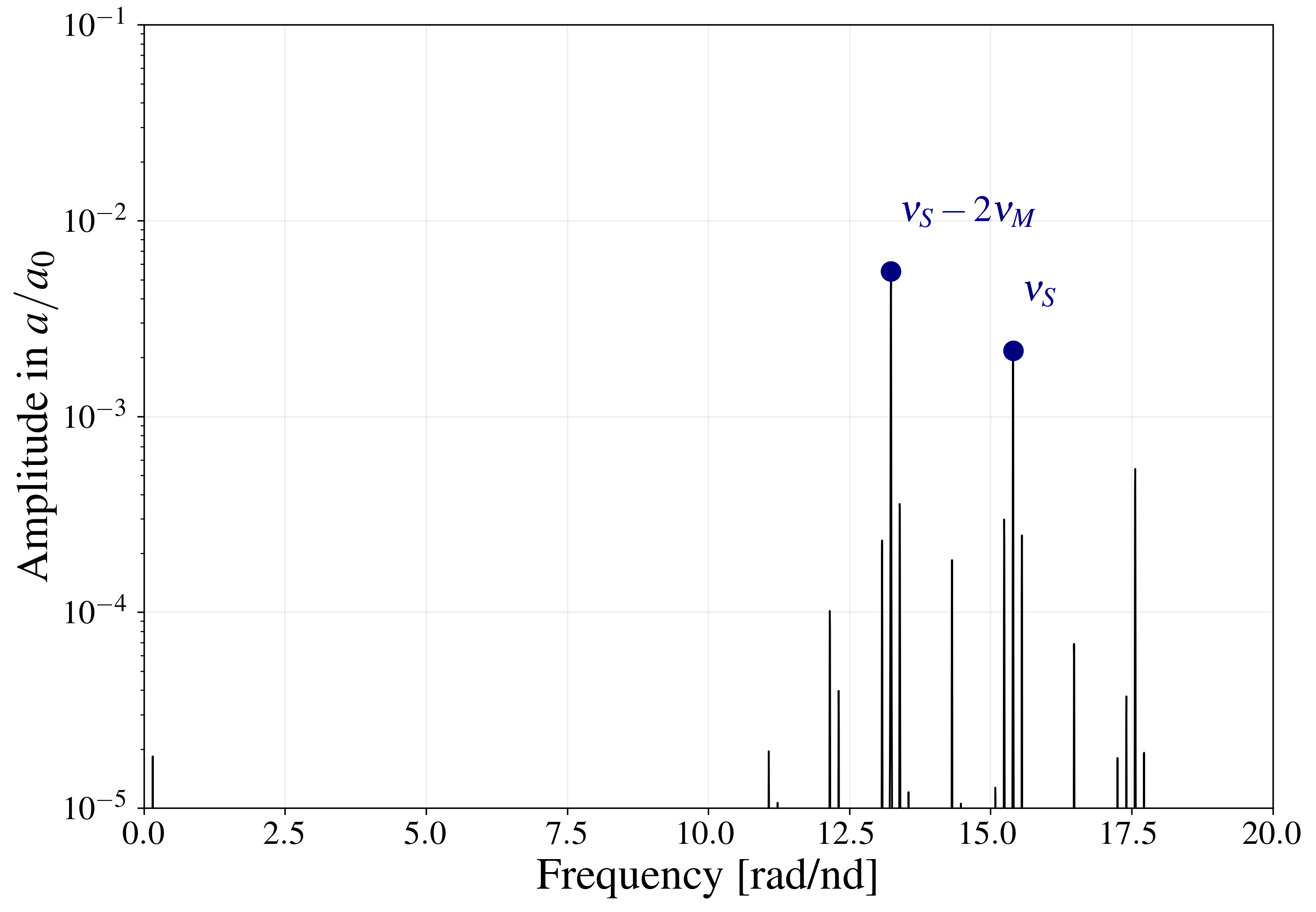}
        \caption{Spectral decomposition of $a/a_0$ with the two dominant frequencies annotated.}\label{fig:cr3bp_fft_a}
    \end{subfigure}
    \caption{Geometry and frequency content of a representative \acrshort{elfo} propagated in the \acrshort{cr3bp} ($A_L = 0$ configuration from Table~\ref{tab:baseline_orbit}; 20-year propagation for the apolune and spectral analyses).}\label{fig:cr3bp_summary}
\end{figure}

Assuming the fundamental quasi-periodic structure is preserved as 2D tori ($A_L = 0$) and 3D tori ($A_L \neq 0$) within the \acrshort{cr3bp}, the design variables $(\nu_S, \nu_M, A_L, \theta_S, \theta_M, \theta_L)$ remain well-defined and are extracted from the trajectory. However, the unaveraged forcing in the \acrshort{cr3bp} deforms the geometry: the orbital elements $a, e, i, \omega$ are no longer separated from the short- and medium-period variations and acquire additional oscillations. For illustration, the initial (mean) \acrshort{eof}-orbital elements from the \acrshort{dadm} frozen equilibrium in Table~\ref{tab:baseline_orbit} are employed as osculating elements and mapped to the \acrshort{brf} via the inverse of Eq.~\eqref{eq:eof_mrf}. The angles $\Omega = t_0 = 0, M = \pi$ are leveraged for the initial configuration. Such a transformation supplies an initial condition for a numerical propagation within the \acrshort{cr3bp}. Figure~\ref{fig:cr3bp_summary}\subref{fig:cr3bp_geom} depicts the \acrshort{cr3bp} trajectory in the rotating frame over the same three-month interval as in Fig.~\ref{fig:freq_structure}, and Fig.~\ref{fig:cr3bp_summary}\subref{fig:strobo_cr3bp} compares the states at the apolunes along each revolution against the \acrshort{dadm} over 20 years. Whereas the \acrshort{dadm} frozen orbit traces a circle of constant radius at apolunes, the \acrshort{cr3bp} orbit deforms into an ellipse-like annulus. This deformation is further illustrated with a frequency-domain analysis.

Two representative scalars are examined under the spectral decomposition: osculating $e$ and $a$ that determine the geometry of the deformed \acrshort{elfo} within the \acrshort{cr3bp}. The other parameters $i, \omega$ evolve in a manner coupled to $e$ and typically do not constitute independent shape parameters, while the angles $\Omega, M$ mainly determine the phase along the \acrshort{elfo}, rather than the geometry itself. The spectra for $e$ and $a/a_0$ with $a_0 = 14{,}200$~km (reference value) are illustrated in Figs.~\ref{fig:cr3bp_fft_e} and \ref{fig:cr3bp_fft_a}, respectively. For the range of oscillatory amplitudes examined ($10^{-5}$ to $10^{-1}$), the eccentricity responds to the unaveraged perturbation from the \acrshort{cr3bp} more sensitively, resulting in peaks with larger magnitudes in the lower-frequency domain that are absent for $a$. Five and two dominant peaks for $e$ and $a$ are annotated in the figures. The corresponding frequencies and amplitudes are detected with a strategy devised by~\citet{gomez2010collocation} and included in Table~\ref{tab:cr3bp_freq}. These peak frequencies are then identified as linear combinations of the frequencies $\nu_S, \nu_M, \nu_L$. The first frequency in $e$ originates from the long-period $\nu_L$, a residual from the imperfect mean-to-osculating transformation; the mapping adopts a simple mean-equals-osculating assumption that is not exact in general, so the frozen condition ($A_L = 0$) is not strictly preserved across the models. The second peak frequency in $e$ is located at $2\nu_M$, an inherent dynamical component arising from tidal deformation by the Earth: the gravitational forcing is aligned primarily with the Earth-Moon direction, producing a perturbation at twice the nodal precession rate. These two main sources of deformations are visible from Fig.~\ref{fig:strobo_cr3bp}. The non-zero $A_L$ expands the circle to an annulus with a finite width. Then, the impact from $2\nu_M$ is visible as the elongation in the $\hat{\bm{x}}$ direction at twice the rate of precession within the \acrshort{mrf}. In Table~\ref{tab:cr3bp_freq}, the peak corresponding to $\nu_L$ is labeled \emph{removable}, as it can be suppressed by targeting $A_L = 0$ via refining the initial condition (Sec.~\ref{sec:fddc}); the remaining peaks are inherent to the unaveraged dynamics (\acrshort{cr3bp}) and cannot be entirely eliminated through design intervention. These oscillatory deformations in the \acrshort{elfo} geometry in turn impact the constellation coverage performance. Before including additional perturbations from the realistic Earth-Moon motion (ephemerides), the \acrshort{cr3bp} supplies an intermediate step to separately gauge the impact from the unaveraged forcing from an idealized Earth motion. 

\begin{table}[h!]
\centering
\caption{Dominant frequency components in the \acrshort{cr3bp} (20-year propagation). Frequencies identified as linear combinations of $\nu_S = 15.399$~rad/nd and $\nu_M = 1.082$~rad/nd.}\label{tab:cr3bp_freq}
\setlength{\tabcolsep}{6pt}
\renewcommand{\arraystretch}{1.2}
\footnotesize
\begin{tabular}{@{} l l c c l c @{}}
\toprule
Signal & & Freq. [rad/nd] & Amplitude & Identification & Removable \\
\midrule
\multirow{5}{*}{$e$}
 & & 0.157 & 0.061 & $\nu_L$ (non-frozen residual) & \checkmark \\
 & & 2.164 & 0.034 & $2\nu_M$ & $\times$ \\
 & & 13.235 & 0.007 & $\nu_S - 2\nu_M$ & $\times$ \\
 & & 2.007 & 0.003 & $2\nu_M - \nu_L$ & $\times$ \\
 & & 15.399 & 0.002 & $\nu_S$ & $\times$ \\
\midrule
\multirow{2}{*}{$a/a_0$}
 & & 13.235 & 0.006 & $\nu_S - 2\nu_M$ & $\times$ \\
 & & 15.399 & 0.002 & $\nu_S$ & $\times$ \\
\bottomrule
\end{tabular}
\end{table}

\subsubsection{\acrshort{hfem}}

\begin{figure}[h!]
    \centering
    \begin{subfigure}[b]{0.48\textwidth}
        \centering
        \includegraphics[width=\linewidth]{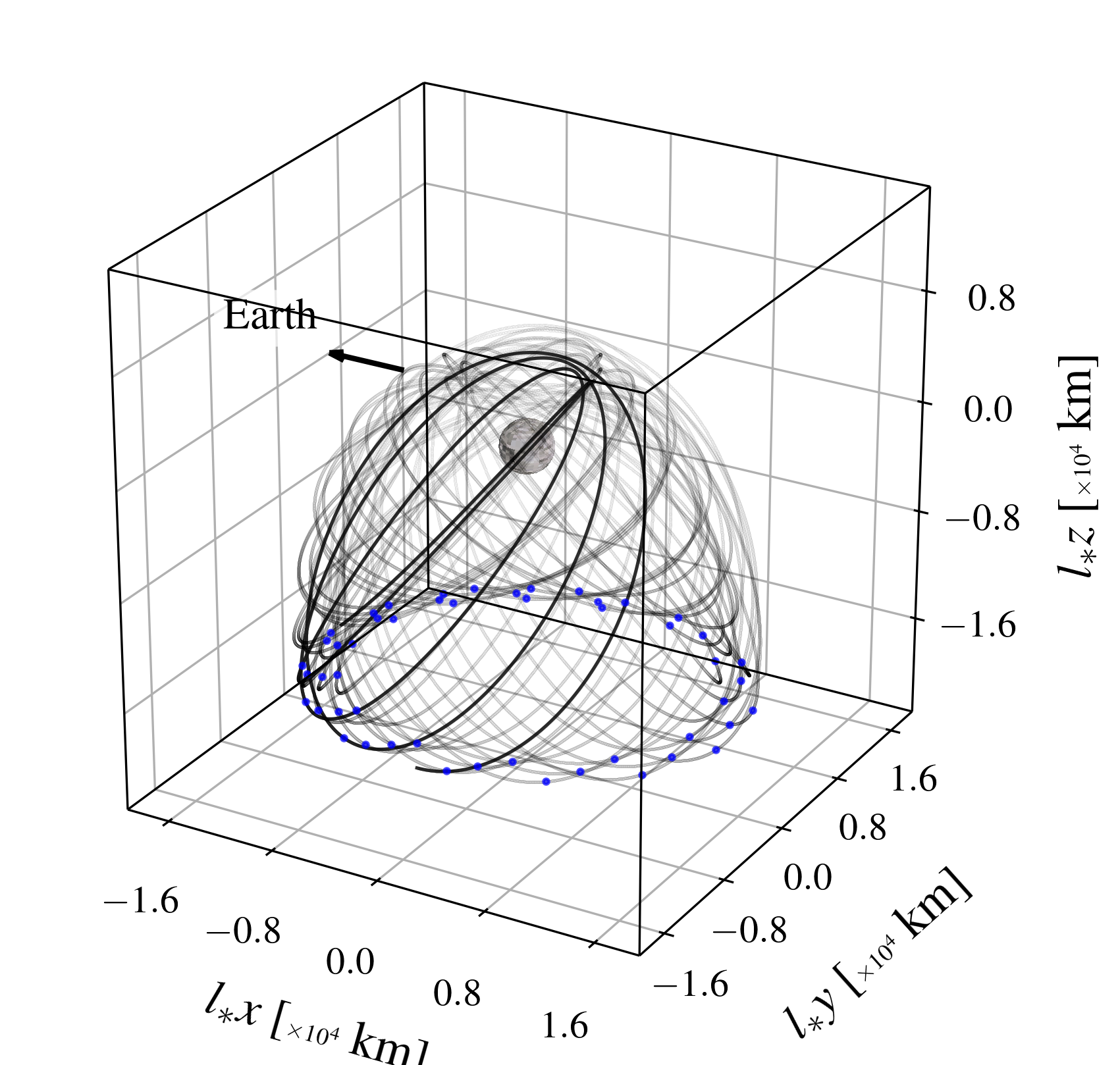}
        \caption{Three-month trajectory (\acrshort{mrf}).}\label{fig:hfem_geom}
    \end{subfigure}
    \hfill
    \begin{subfigure}[b]{0.48\textwidth}
        \centering
        \includegraphics[width=\linewidth]{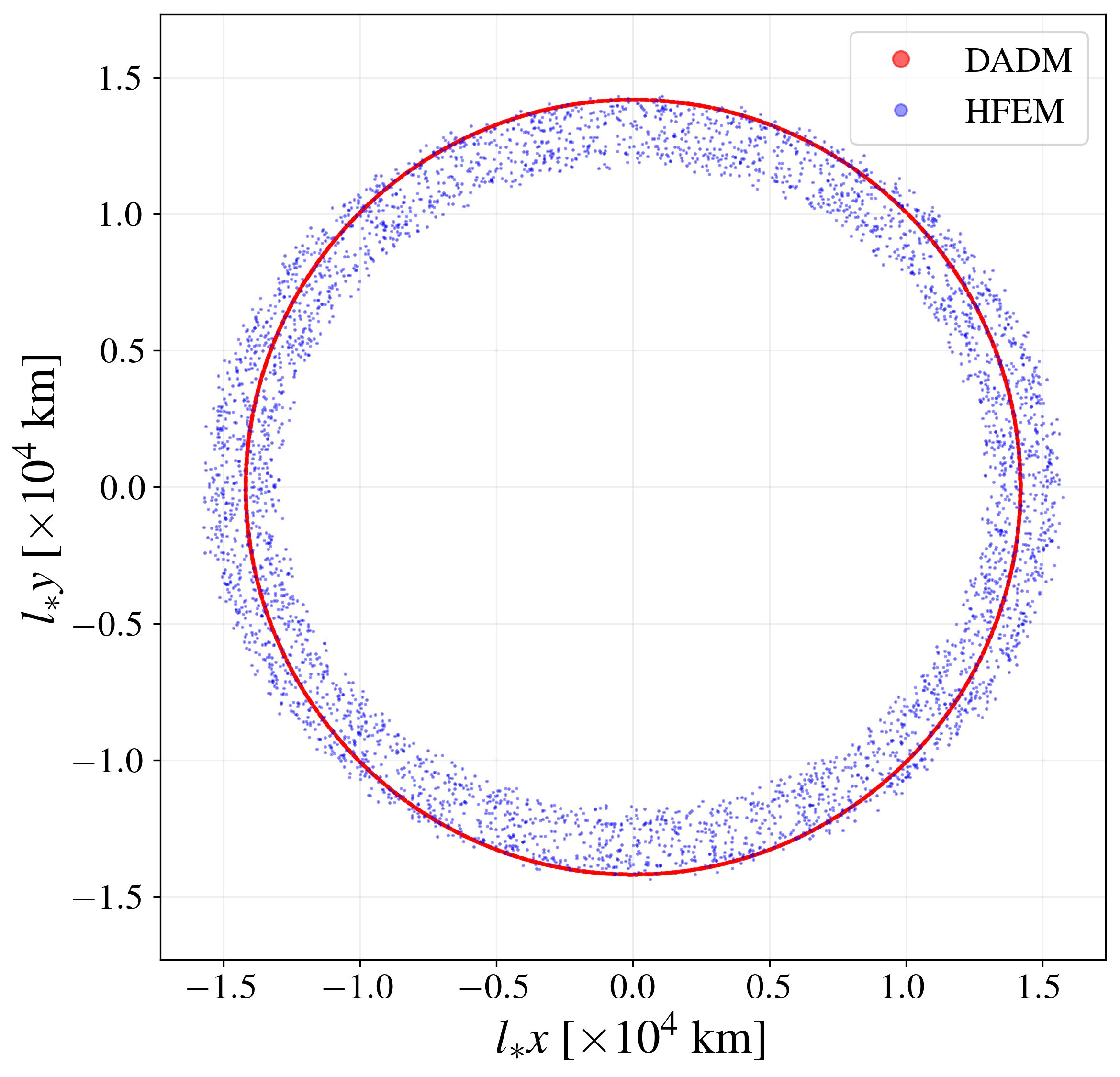}
        \caption{Apolune states (20 years), \acrshort{dadm} vs \acrshort{hfem}.}\label{fig:strobo_hfem}
    \end{subfigure}
    \\[6pt]
    \begin{subfigure}[b]{0.48\textwidth}
        \centering
        \includegraphics[width=\linewidth]{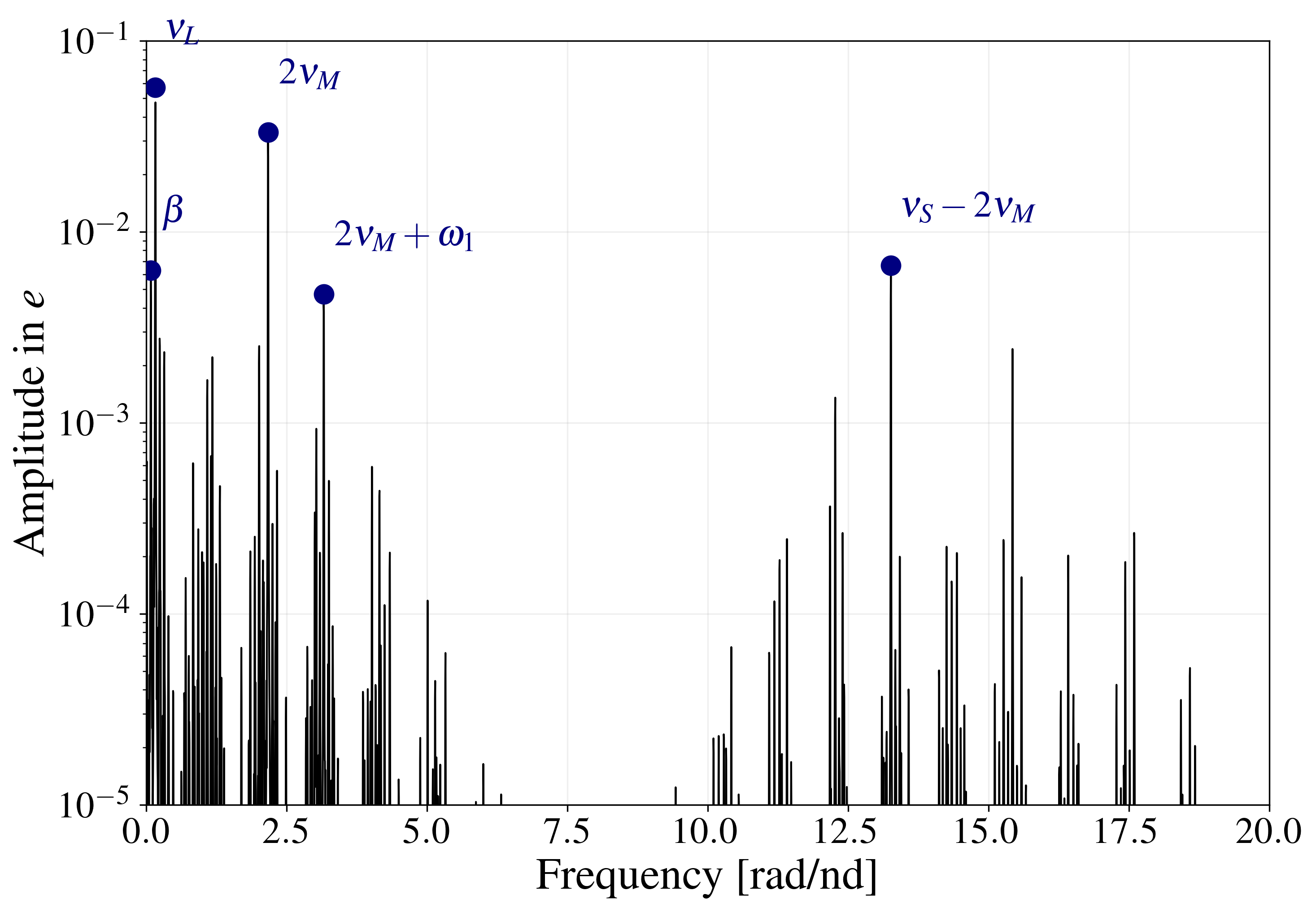}
        \caption{Spectral decomposition of $e$ with the five dominant frequencies annotated.}\label{fig:hfem_fft_e}
    \end{subfigure}
    \hfill
    \begin{subfigure}[b]{0.48\textwidth}
        \centering
        \includegraphics[width=\linewidth]{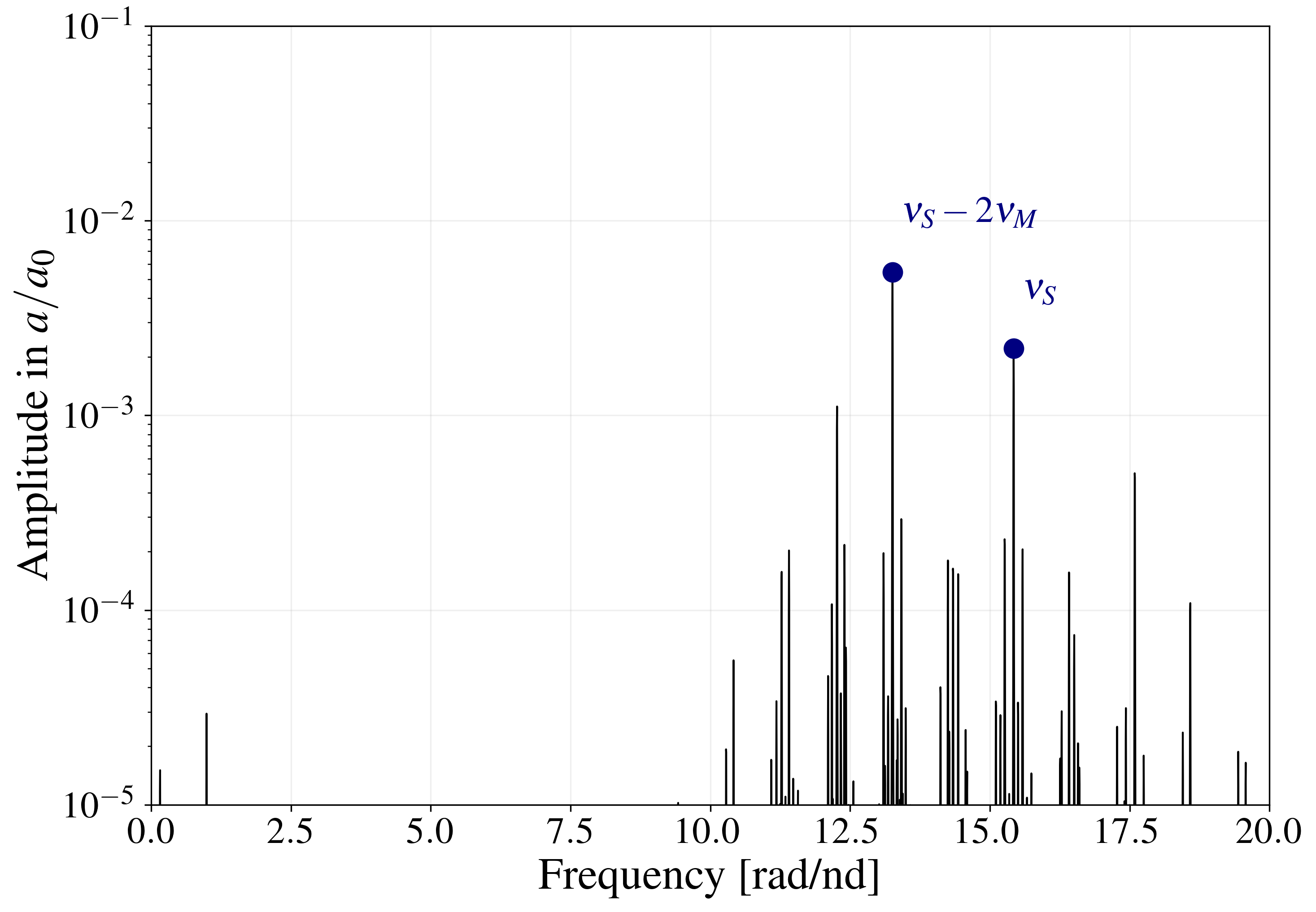}
        \caption{Spectral decomposition of $a/a_0$ with the two dominant frequencies annotated.}\label{fig:hfem_fft_a}
    \end{subfigure}
    \caption{Geometry and frequency content of a representative \acrshort{elfo} propagated in the \acrshort{hfem} ($A_L = 0$ configuration from Table~\ref{tab:baseline_orbit}; 20-year propagation for the apolune and spectral analyses).}\label{fig:hfem_summary}
\end{figure}

The quasi-periodic motion is also assumed to persist within the \acrshort{hfem} for the specific \acrshort{elfos} that are the focus of this investigation. Following the transformation scheme from Eq.~\eqref{eq:eof_mci_snapshot}, a Cartesian initial condition within the \acrshort{mci} is retrieved from the \acrshort{eof} elements at $A_L = 0$ in Table~\ref{tab:baseline_orbit} and then propagated within the \acrshort{hfem}. The initial reference epoch is at $JD_0 = 2460941.5$ (2025-09-23 00:00:00 UTC). Figures~\ref{fig:hfem_geom} and \ref{fig:strobo_hfem} demonstrate that the \acrshort{hfem} trajectory and apolune distribution remain qualitatively similar to the \acrshort{cr3bp} counterparts (Figs.~\ref{fig:cr3bp_geom} and \ref{fig:strobo_cr3bp}) within an ellipse-like annulus behavior emerging at the apolunes.

While the \acrshort{hfem} \acrshort{elfos} remain quasi-periodic (by assumption), the time-dependent \acrshort{hfem} introduces additional \emph{external} frequencies~\cite{jorba1997persistence}. In the \acrshort{hfem} dynamics, these frequencies originate from the gravitational (conservative) perturbations due to the realistic motion of the Sun-Earth-Moon system. The external frequencies draw contrast from the internal frequencies of the \acrshort{elfos} in that the frequencies are constant, predetermined from the Sun-Earth-Moon ephemerides, and the respective phase angles only depend on the epoch (Julian Date) \cite{park2025numerical, sanaga2026systematic},
   \begin{align}
      \label{eq:omega}  JD_0 \;\longrightarrow\; (\bm{\omega}, \bm{\alpha})
        = (\omega_1, \dots, \omega_{n_{ext}}, \alpha_1, \dots, \alpha_{n_{ext}}).
    \end{align}
with a linear relationship $d\bm{\alpha}/dt = \bm{\omega}$. While the internal frequency-phase structure $(\bm{\nu}, \bm{\theta}) = (\nu_S, \nu_M, \nu_L,\allowbreak \theta_S, \theta_M, \theta_L)$ introduced in Sec.~\ref{subsec:freq-modes} is determined by the initial condition, the external frequency-phase structure $(\bm{\omega}, \bm{\alpha})$ arises from the time dependence of the system and mirrors the epoch information. This distinction is fundamental: the internal frequencies describe the spacecraft motion and may shift in both value and phase via adjustment of initial condition, while the external frequencies are inherited from the quasi-periodic forcings of the gravitational environment and remain constant for a reference epoch. Three external frequencies are particularly relevant~\cite{park2025numerical,gomez2002solar} in governing the geometry change under the \acrshort{hfem}: $\omega_1 = 0.992$~rad/nd (anomalistic, apsidal precession of the lunar orbit; period $\approx 27.5$~day), $\omega_2 = 0.925$~rad/nd (synodic, Sun-Earth-Moon geometry; $\approx 29.5$~day), and $\omega_3 = 1.004$~rad/nd (draconic, lunar nodal regression; $\approx 27.2$~day).

\begin{figure}[h!]
    \centering
    \includegraphics[width=0.55\textwidth]{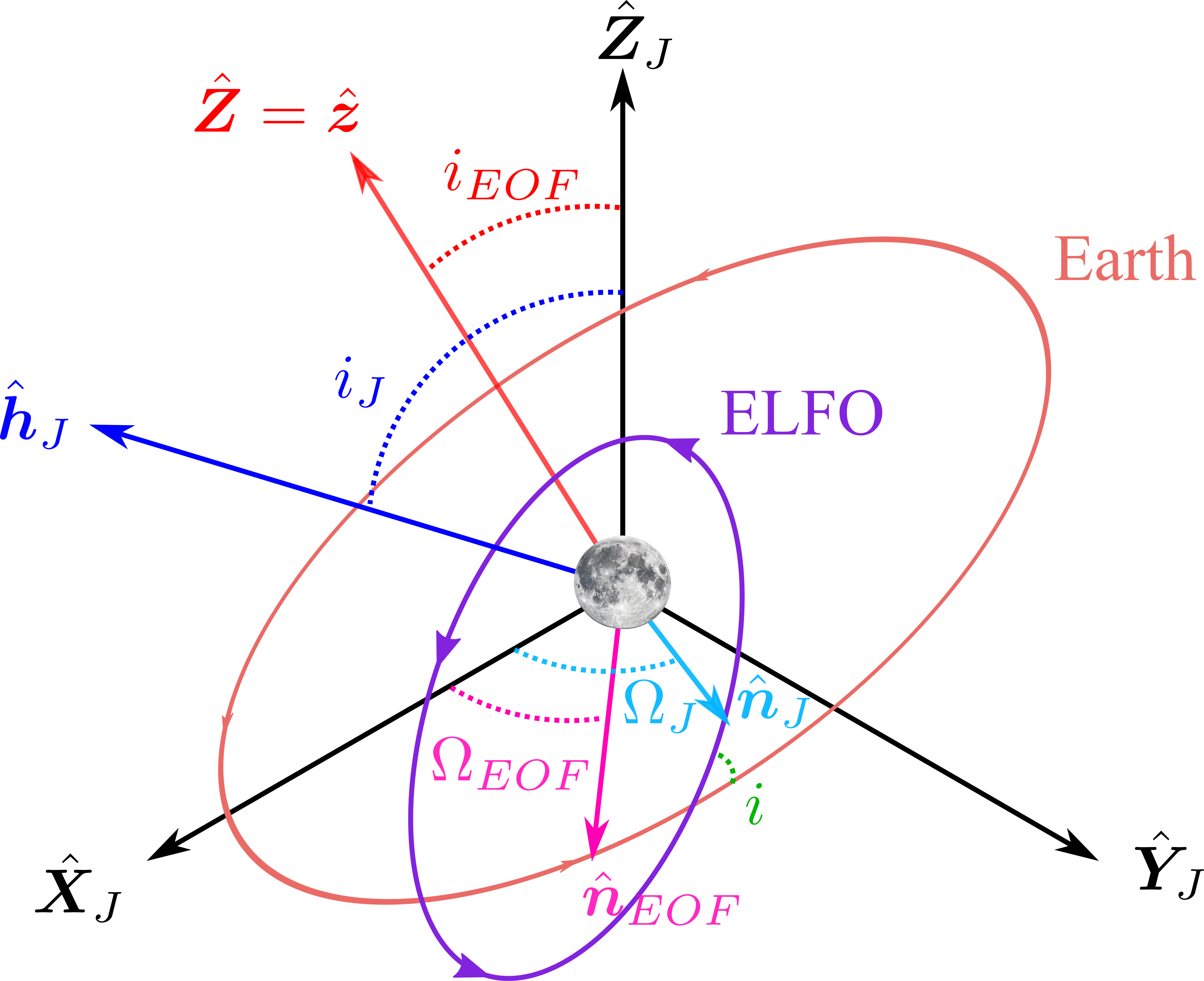}
    \caption{Geometry of the $\beta$ beat mechanism in the \acrshort{mci} (not to scale).}\label{fig:omega_lu_schematic}
\end{figure}

Noting the existence of the extra frequency structures within the \acrshort{hfem}, the spectra for osculating $e$ and $a$ are examined for the \acrshort{hfem} \acrshort{elfo} trajectory in Figs.~\ref{fig:hfem_fft_e}-\ref{fig:hfem_fft_a}. Compared to the \acrshort{cr3bp} counterparts (Figs.~\ref{fig:cr3bp_fft_e}-\ref{fig:cr3bp_fft_a}), two new frequency components are revealed in $e$. The first component is the oscillation at $\beta = \nu_M - \omega_3$, with an amplitude of $0.006$ in $e$ (fourth peak within Table~\ref{tab:hfem_freq}) and a period of approximately 349 days. Its geometric origin is traced by examining the evolution of the \acrshort{eof}, a non-inertial frame, within the inertial \acrshort{mci}. Figure~\ref{fig:omega_lu_schematic} illustrates the osculating orbital plane of an \acrshort{elfo} and the Earth orbital plane within the \acrshort{mci}. The spacecraft angular momentum direction $\hat{\bm{h}}_J$ and the ascending node $\hat{\bm{n}}_J$ define the \acrshort{elfo} orbital plane, while the Earth's osculating orbital plane around the Moon defines the instantaneous \acrshort{eof}, with angular momentum direction $\hat{\bm{Z}} = \hat{\bm{z}}$ and ascending node $\hat{\bm{n}}_{EOF}$. Let $i_J$ and $i_{EOF}$ denote the inclinations of the \acrshort{elfo} and the Earth-Moon plane within the \acrshort{mci}, respectively, and define $\Delta\Omega_J = \Omega_J - \Omega_{EOF}$ as their relative ascending node. The inclination of the \acrshort{elfo} with respect to the \acrshort{eof}, i.e., $i$, is not governed by the secular dynamics captured within the \acrshort{dadm} alone, since the \acrshort{eof} is itself in motion within the \acrshort{mci}. With respect to the ecliptic, the Earth-Moon plane maintains a nearly constant inclination $i_{EOF} \approx 5.15^\circ$ while its ascending node regresses at an approximate rate $\dot{\Omega}_{EOF} \approx 1 - \omega_3$ over the 18.6-year cycle \cite{gomez2002solar}. Since $i$ is measured from a plane that is itself rotating, its rate carries a contribution from that rotation, in addition to the gravitational torque on the \acrshort{elfo}. This contribution to the inclination is approximated as,
\begin{align}
    \label{eq:di_frame}
    \left(\frac{di}{dt}\right)_{\!\text{frame}} \approx -\dot{\Omega}_{EOF}\,\sin i_{EOF}\,\sin\Delta\Omega_J.
\end{align}
Noting that $\dot{\Omega}_J \approx 1 - \nu_M$ for the spacecraft, the relative node evolves at the beat frequency $\beta$ as,
\begin{align}
    \label{eq:beta_beat}
    d\Delta\Omega_J/dt = \dot{\Omega}_J - \dot{\Omega}_{EOF} \approx (1 - \nu_M) - (1 - \omega_3) = -(\nu_M - \omega_3) = -\beta.
\end{align}
The contribution in Eq.~\eqref{eq:di_frame} therefore oscillates at $|\beta|$. Such modulation couples into the eccentricity and the argument of perilune through Eqs.~\eqref{eq:dedt}-\eqref{eq:domegadt}, resulting in the $\beta$ peak within the spectrum for $e$ (Fig.~\ref{fig:hfem_fft_e}, Table~\ref{tab:hfem_freq}).\footnote{The gravitational torque largely cancels this frame contribution within the inclination itself, so the $\beta$ beat is considerably more visible in $e$ than in $i$.} The physical mechanism originates from the precession of the \acrshort{eof} within an inertial frame, driven primarily by the solar gravitational perturbation on the Earth-Moon system~\cite{gomez2002solar}. Accordingly, this beat frequency effect on the \acrshort{elfo} geometry does not exist within the \acrshort{dadm} or \acrshort{cr3bp}; rather, it acts within the \acrshort{hfem}, periodically perturbing the ``frozen'' geometry by driving a forced response of the $(e, \omega)$ oscillator at $\beta$, off-resonant from the free libration at $\nu_L$ (Sec.~\ref{subsec:freq-modes}). The $\beta$ beat is a structural mechanism that cannot be entirely removed through design intervention. To the authors' knowledge, the spectral identification of this beat and its coverage impact have not previously been emphasized in the \acrshort{elfo} constellation literature; the quantitative impact for the present reference configuration is examined in Sec.~\ref{subsec:hf-gravity-impact}.

\begin{table}[h!]
\centering
\caption{Dominant frequency components in the \acrshort{hfem} (20-year propagation).}\label{tab:hfem_freq}
\setlength{\tabcolsep}{6pt}
\renewcommand{\arraystretch}{1.2}
\footnotesize
\begin{tabular}{@{} l l c c l c @{}}
\toprule
Signal & & Freq. [rad/nd] & Amplitude & Identification & Removable \\
\midrule
\multirow{5}{*}{$e$}
 & & 0.159 & 0.057 & $\nu_L$ (non-frozen residual) & \checkmark \\
 & & 2.167 & 0.033 & $2\nu_M$ & $\times$ \\
 & & 13.258 & 0.007 & $\nu_S - 2\nu_M$ & $\times$ \\
 & & 0.078 & 0.006 & $\beta = \nu_M - \omega_3$ & $\times$ \\
 & & 3.158 & 0.005 & $2\nu_M + \omega_1$ & $\times$ \\
\midrule
\multirow{2}{*}{$a/a_0$}
 & & 13.258 & 0.006 & $\nu_S - 2\nu_M$ & $\times$ \\
 & & 15.425 & 0.002 & $\nu_S$ & $\times$ \\
\bottomrule
\end{tabular}
\end{table}

\newpage
The second frequency component absent in the \acrshort{cr3bp} is the anomalistic sideband $2\nu_M + \omega_1$, with an amplitude 0.005 in $e$ and a period of approximately 9 days. The apsidal precession of the lunar orbit at rate $\omega_1$ modulates the instantaneous Earth-Moon distance, thereby varying the amplitude of the unaveraged tidal forcing responsible for the $2\nu_M$ harmonic and producing sidebands at $2\nu_M \pm \omega_1$, of which only the upper sideband appears among the dominant peaks in Table~\ref{tab:hfem_freq}. Together with the $2\nu_M$ harmonic inherited from the \acrshort{cr3bp} and the $\beta$ beat, these sidebands emerge from perturbations inherent to the dynamics and constitute \emph{unremovable} components that cannot be eliminated entirely through design intervention. In contrast, $\nu_L$ remains \emph{removable} across all three models, as it can be suppressed by enforcing $A_L = 0$ via adjustment in the initial conditions (Sec.~\ref{sec:fddc}). Although the present analysis is based on a single representative configuration, the identified frequency components originate from the structure of the dynamics rather than from the particular choice of the orbital elements; the mechanisms themselves are therefore expected to persist, while their relative amplitudes and coverage impact remain configuration-dependent and do not generalize numerically from the present case.

\subsection{Lunar Obliquity and User Geometry}
\label{subsec:obliquity}
\begin{figure}[h!]
    \centering
    \begin{subfigure}[b]{0.48\textwidth}
        \centering
        \includegraphics[width=\linewidth]{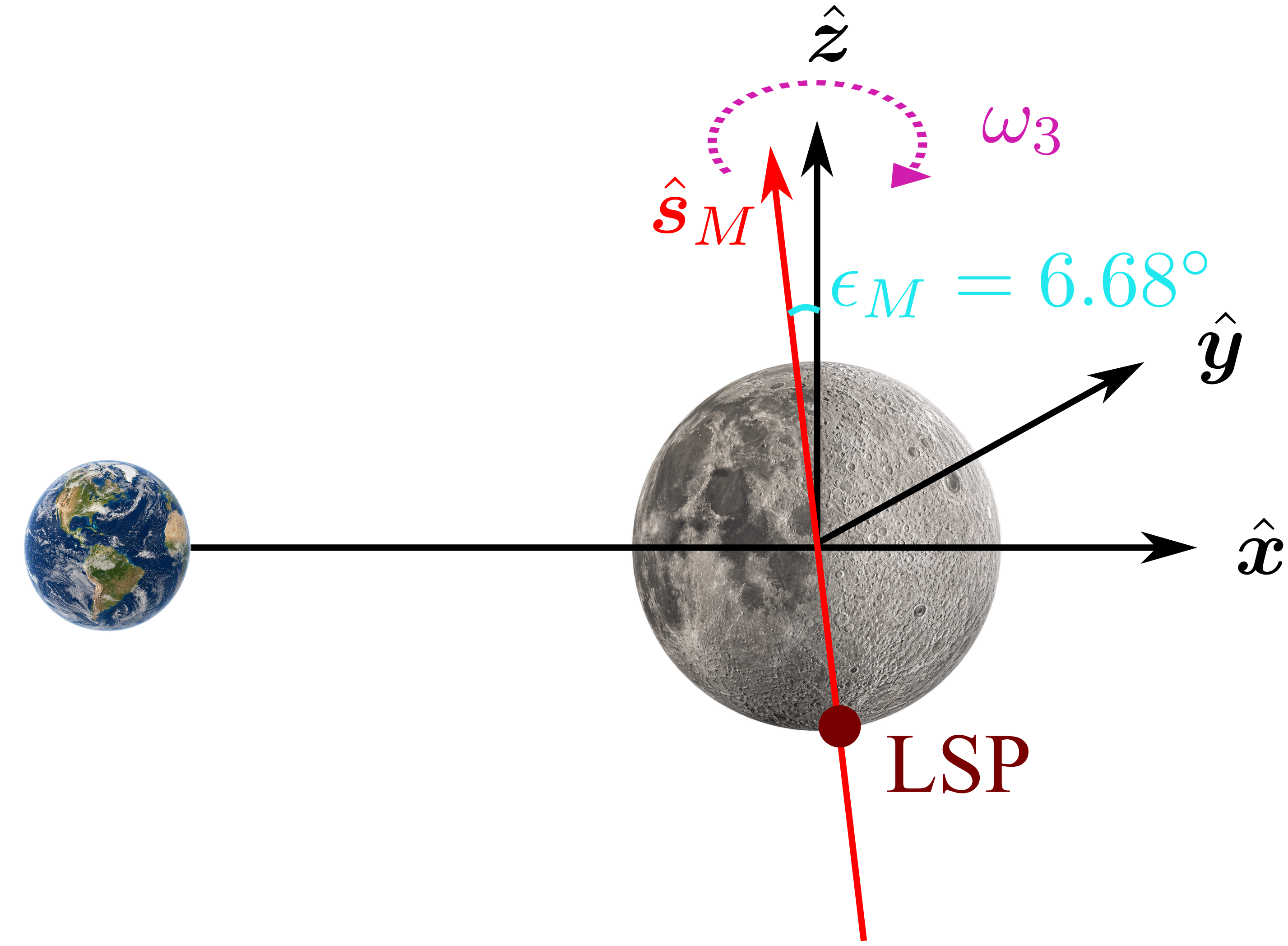}
        \caption{Schematic of the lunar spin axis $\hat{\bm{s}}_M$ tilted by $\approx 6.68^\circ$ from $\hat{\bm{z}}$ within the \acrshort{mrf} (not to scale).}\label{fig:user_geom_schematic}
    \end{subfigure}
    \hfill
    \begin{subfigure}[b]{0.48\textwidth}
        \centering
        \includegraphics[width=\linewidth]{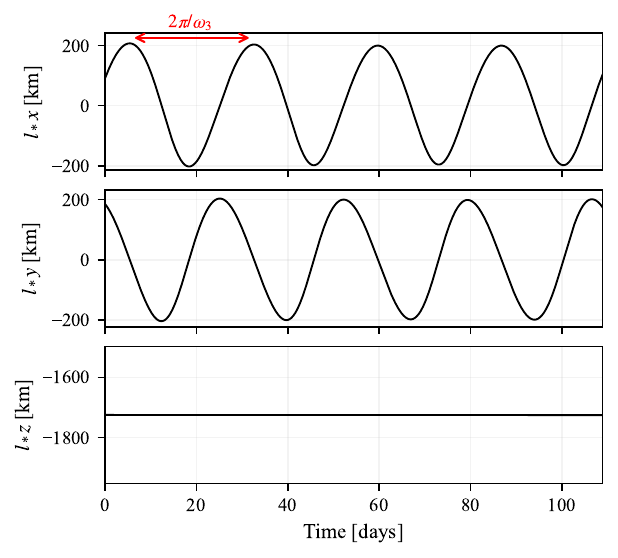}
        \caption{\acrshort{lsp} position components in the \acrshort{mrf} over four draconic months, evaluated with DE440 \cite{park2021jpl}.}\label{fig:user_geom_history}
    \end{subfigure}
    \caption{User-geometry frequency content induced by the lunar obliquity.}\label{fig:user_geom}
\end{figure}

Any fixed location on the lunar surface, e.g., the \acrshort{lsp}, evolves within a reference frame due to the precession of the lunar axis of rotation. This user-geometry evolution is an independent channel leading to time-dependent oscillations in the constellation coverage performance\footnote{Note the difference with Sec.~\ref{subsec:freq-across-models}, where a \emph{dynamical} evolution of \acrshort{elfo} geometry examined under perturbations; the current subsection details the \emph{kinematical} changes in the user location relative to the \acrshort{elfos}.}, analyzed within the \acrshort{mrf}. The lunar axis of rotation ($\hat{\bm{s}}_M$) is tilted with respect to the (instantaneous) \acrshort{eof} normal vector ($\hat{\bm{Z}} = \hat{\bm{z}}$) by approximately $\epsilon_M = 6.68^\circ$ (Fig.~\ref{fig:user_geom_schematic}), corresponding to the Cassini-state spin-orbit angle that remains nearly constant in the Moon's rotational dynamics~\cite{peale1969generalized}. The tilted direction aligns with the lunar ascending node that drifts at the rate $\dot{\Omega}_{EOF} = -(\omega_3 - 1)$ (Fig.~\ref{fig:omega_lu_schematic}) in the inertial frame \cite{gomez2002solar}, and $-\omega_3$ within the \acrshort{mrf}. As such, lunar spin axis $\hat{\bm{s}}_M$ precesses with $\omega_3$ with respect to $\hat{\bm{z}}$ in the clockwise direction, tracing an approximate cone of half-angle $\approx 6.68^\circ$ within the \acrshort{mrf} as observed in Fig.~\ref{fig:user_geom_schematic}. For illustration, the \acrshort{lsp}'s apparent location within the \acrshort{mrf} is traced in Fig.~\ref{fig:user_geom_history} for 4 draconic months ($2\pi/\omega_3$)\footnote{The \texttt{MOON\_PA} frame from the JPL kernel \texttt{moon\_pa\_de440\_200625.bpc} is leveraged \cite{park2021jpl}.}; the $z$-component of the \acrshort{lsp} position remains essentially fixed at $\approx -1725$~km (a constant effective latitude of $\varphi = -83.32^\circ$), while the $(x, y)$ components oscillate at the draconic period $2\pi/\omega_3$ with amplitudes of approximately $\pm 200$~km. Although slower modulations exist (e.g., physical libration and free wobble of the lunar pole~\cite{rambaux2011moon}), the constellation operation horizon is dominated by the nearly periodic user-geometry change at $\omega_3$. In the context of the relative geometry between a fixed user in the lunar body-fixed frame (\texttt{MOON\_PA}) and an \acrshort{elfo}, the drift of $\hat{\bm{s}}_M$ with respect to the orbital plane plays a significant role. This relative configuration is modulated by the same beat frequency $\beta = \nu_M - \omega_3$ introduced in Sec.~\ref{subsec:freq-across-models}, with the rotating-frame construction viewed in Fig.~\ref{fig:user_geom}\subref{fig:user_geom_schematic} mirroring the geometric evolution in Fig.~\ref{fig:omega_lu_schematic}. The coverage-level manifestation of this user-geometry contribution is examined further in Sec.~\ref{subsec:hf-gravity-impact} and Appendix~\ref{app:obliquity}.  

\section{Design Space Exploration in the DADM}
\label{sec:design-exploration}
The frequency structures identified in Sec.~\ref{sec:freq-structure} facilitate a systematic exploration of the constellation design space. The assumed constellation architecture is first described, followed by a symmetry-reduced design framework that leverages the quasi-periodic behavior within the \acrshort{dadm}, serving as a lower-fidelity model. The global design space is then surveyed under a specific performance metric. Lastly, the role of the long-period amplitude $A_L$ is investigated for a sample constellation. 

\subsection{Assumed Constellation Architecture and Performance Metric}
\label{sec:assumed_constellation}
Diverse constellation configurations may suit a variety of objectives. The current analysis is motivated by a specific architecture envisioned for NASA's \acrfull{lcrns}, planned to be launched and operated by Intuitive Machines~\cite{brack2025}. Rather than adopting the full mission-specific details, the current analysis retains only the generic aspects of the proposed constellation as a representative architecture that leverages \acrshort{elfos}. The constellation consists of five satellites with a nominal operational period of ten years. Coverage is assessed at the \acrfull{lsp} using the \acrfull{gdop}, a standard metric that quantifies the geometric quality of a set of ranging sources. An elevation mask of $5^\circ$ is applied to exclude satellites near the local horizon, and at each evaluation epoch the number of visible satellites above the mask is denoted $n_{\mathrm{vis}}$. The \acrshort{gdop} is constructed from the satellite-user geometry as
\begin{align}
    \label{eq:gdop}
    \text{GDOP} = \sqrt{\text{tr}\!\left((\bm{H}^T \bm{H})^{-1}\right)}.
\end{align}
The observation matrix $\bm{H} \in \mathbb{R}^{n_{\mathrm{vis}} \times 4}$ has a $k$-th row ($1\leq k \leq n_{\mathrm{vis}}$) consisting of the unit direction vector from the user to the $k$-th visible satellite (first three components) concatenated with unity (fourth component, accommodating the receiver clock bias). When $n_{\mathrm{vis}} < 4$, $\bm{H}^T \bm{H}$ is rank-deficient and \acrshort{gdop} is undefined; such epochs are treated as non-covered. A lower \acrshort{gdop} value indicates a geometrically well-distributed set of visible satellites, implying that pseudorange errors to the satellites are less amplified when inverted to yield the user's position and clock-bias estimate. The current analysis adopts a threshold of \acrshort{gdop} $\leq 6$ following Brack et al.~\cite{brack2025}. The \emph{coverage fraction} is defined as the proportion of evaluation epochs satisfying the \acrshort{gdop} criterion, with undefined-\acrshort{gdop} epochs counted as non-covered, and serves as the metric for evaluating a given constellation.

\subsection{Symmetry-Leveraging Design Framework}
\label{subsec:design-framework}
The \acrshort{gdop} is a discontinuous function of the satellite-user geometry because of the line-of-sight visibility requirement. This property typically renders gradient-based optimization ineffective under \acrshort{gdop}-based performance metrics, so meta-heuristic methods often serve as a practical alternative in supplying the \acrshort{elfo} constellations. Such strategies, while effective for a specific point design, scale poorly with the numerical propagation horizons required in higher-fidelity environments and are therefore ill-suited for a preliminary survey of the design space.

In the initial design stage, dynamical fidelity may be traded for a tractable constellation evaluation strategy compatible with meta-heuristic optimization. Under the assumptions introduced below, the \acrshort{gdop} itself exhibits quasi-periodic behavior, admitting a parameterization over the angular domain $(\theta_S, \theta_M) \in [0, 2\pi)^2$ rather than the time domain and thereby compressing a potentially long propagation window into a finite two-dimensional torus; the terms \emph{angular domain} and \emph{torus} denote this same domain throughout the current analysis. Exact symmetries further collapse the effective design space, enabling a rapid global survey of the constellation landscape. Many of these symmetries are intuitive and partially leveraged in literature~\cite{ceresoli2025design, brack2025}; the present work supplies a comprehensive and explicit treatment of the possible symmetry structures. The framework follows in three steps: Sec.~\ref{subsec:symmetry-assumptions} states the assumptions that render the constellation symmetric, Sec.~\ref{subsec:design-vars} identifies the resulting design variables together with the three symmetries that reduce them, and Sec.~\ref{subsec:optimization-framework} poses the optimization over the reduced space.

\subsubsection{Symmetry-Inducing Assumptions}
\label{subsec:symmetry-assumptions}
The following set of assumptions renders the constellation symmetric in several respects:
\begin{itemize}
    \item \ul{Time-independent dynamics} (\acrshort{dadm} or \acrshort{cr3bp}): Both models admit the same frequency-angle parametrization $(\nu_S, \nu_M, A_L, \theta_S, \theta_M, \theta_L)$ for a satellite state on an \acrshort{elfo} (Sec.~\ref{subsec:freq-modes}) independent from an epoch.
    \item \ul{Same short- and medium-period frequencies}: To supply prolonged coverage, consistent values of the short- and medium-period frequencies are enforced across the constellation. Otherwise, the satellites eventually drift in the respective phases, deteriorating the overall constellation performance. Under this assumption, the phase offsets in $\theta_S, \theta_M$ are maintained.
    \item \ul{Zero long-period amplitude} ($A_L = 0$): A non-zero long-period amplitude $A_L$ introduces oscillations in the coupled triad $(e, i, \omega)$ (Sec.~\ref{subsec:freq-modes}), lifting each spacecraft state on a 2D torus to a 3D torus within time-independent dynamics. For such motions, since $\nu_L$ may not be consistent within the constellation, the overall \acrshort{gdop} behaves as a multi-dimensional quasi-periodic function. In either time-domain or angular-domain evaluations, capturing this additional degree of freedom demands a substantially larger number of samples. Adopting the frozen condition ($A_L = 0$) is therefore a practical selection that maintains a tractable angular parametrization. This symmetry-inducing assumption is revisited in Sec.~\ref{subsec:al-dadm}, where suppressing the long-period amplitude is not only tractable but also favorable for the adopted \acrshort{gdop}-threshold metric as examined in a case study.
    \item \ul{Constant lunar obliquity} ($\epsilon_M = 6.68^\circ$): The rotational axis of the Moon ($\hat{\bm{s}}_M$, Fig.~\ref{fig:user_geom_schematic}) is assumed to be tilted by a constant obliquity $\epsilon_M = 6.68^\circ$ with respect to the orbit normal ($\hat{\bm{z}}$). An additional assumption may be introduced that eliminates this obliquity, i.e., $\hat{\bm{s}}_M = \hat{\bm{z}}$. The impact of such an assumption in \acrshort{elfo} constellation design is further examined in Appendix~\ref{app:onaxis}. 
    \item \ul{Observer location at the LSP}: With the fixed user location at the \acrshort{lsp}, the latitude is also fixed within the \acrshort{mrf} at $\varphi = -83.32^\circ$. This assumption may be relaxed to examine coverage at other latitudes or specific surface locations at an increased evaluation cost; the main body of the current analysis retains the \acrshort{lsp} focus throughout, whereas Appendix~\ref{app:non-lsp} examines the non-\acrshort{lsp} user locations.
\end{itemize}
Although these assumptions do not capture the full realistic dynamical and mission environment, they are collectively justified for a preliminary analysis aimed at a global survey of the design space rather than a point-solution optimization within a higher-fidelity model. Explicitly recording these assumptions is also beneficial in interpreting the constellation behaviors within realistic models as deviations from this simplified, symmetric reference modeling environment.

\subsubsection{Design Variables and Symmetries}
\label{subsec:design-vars}
With the previous assumptions, the design variables for the constellation reduce to (1) the orbital geometry, specified by the frequency pair $(\nu_S, \nu_M)$, constant across all satellites, and (2) the phase variables $\theta_{S,k}, \theta_{M, k}$ for $1 \leq k \leq 5$ (satellite index). Under the frozen condition ($A_L = 0$) with $\omega = \pi/2$ and $i < 90^\circ$\footnote{\acrshort{elfos} are assumed to be in a prograde motion, although it is not a general constellation requirement.}, the frequencies $(\nu_S, \nu_M)$ are determined from $a$ and $e$ in closed forms within the \acrshort{dadm} \cite{park2026bridging},
\begin{align}
    \label{eq:nu_s} \nu_S &= \frac{n}{n_E} = \sqrt{\frac{\tilde{\mu}_M}{a^3}}\,t_*, \\
    \label{eq:nu_m} \nu_M &= 1 + \frac{1}{4}(1-\mu)\frac{n_E}{n}\,\sqrt{\tfrac{3}{5}}\,(20\sin^2 i - 5).
\end{align}
The leading value of unity in the $\nu_M$ expression reflects the rotation rate of the \acrshort{mrf}. When $\nu_S$ and $\nu_M$ are not in resonance, the constellation state does not repeat itself. Exploiting the quasi-periodicity, however, the reference phase angles $(\theta_{S,1}, \theta_{M,1})$ of satellite~1 serve as surrogates for time: as $t$ evolves, these phases sweep the torus $[0, 2\pi)^2$, so sampling this torus is equivalent to examining the constellation over an infinite horizon within the assumed dynamics and symmetric configurations (Sec.~\ref{subsec:symmetry-assumptions}). Along with this formulation, three symmetries (invariances) further reduce the effective design space:
\begin{itemize}
    \item \ul{Translational symmetry}: Over an infinite horizon, a uniform shift in $\theta_S$ or $\theta_M$ does not change the constellation; therefore, the phase angles for satellite~1 may be set to zero, i.e., $\theta_{S,1} = \theta_{M, 1} = 0$ at $t = 0$, and removed from the design variables. This invariance is also leveraged in \citet{brack2025}.
    \item \ul{Permutational symmetry}: Permuting any two satellites in their respective phase angles does not change the constellation, i.e., $(\theta_{S, i}, \theta_{M, i}) \leftrightarrow (\theta_{S, j}, \theta_{M, j})$ for $i\neq j$ represents the same configuration. Therefore, satellites 2-5 may be organized in ascending order with respect to the short-period phase $\theta_S$ without loss of generality, reducing the effective search space by a factor of $4! = 24$.
    \item \ul{Longitudinal symmetry} (\acrshort{dadm} only): Denote the longitude for the user (at the \acrshort{lsp}) within the \acrshort{mrf} as $\lambda$. Evaluation of the constellation over an infinite time horizon yields identical performance at all longitudes for a given latitude within the \acrshort{mrf} as orbital precession rate is constant within the \acrshort{dadm}. As such, $\lambda = 0$ evaluates all possible configurations. 
\end{itemize}
Given these symmetries combined, the quasi-periodic \acrshort{gdop} behavior is compactly represented over the two-dimensional angular grid $(\theta_{S, 1}, \theta_{M, 1})$. 

\subsubsection{Optimization Formulation}
\label{subsec:optimization-framework}
\begin{figure}[htb]
    \centering
    \begin{subfigure}[b]{0.48\textwidth}
        \centering
        \includegraphics[width=\linewidth]{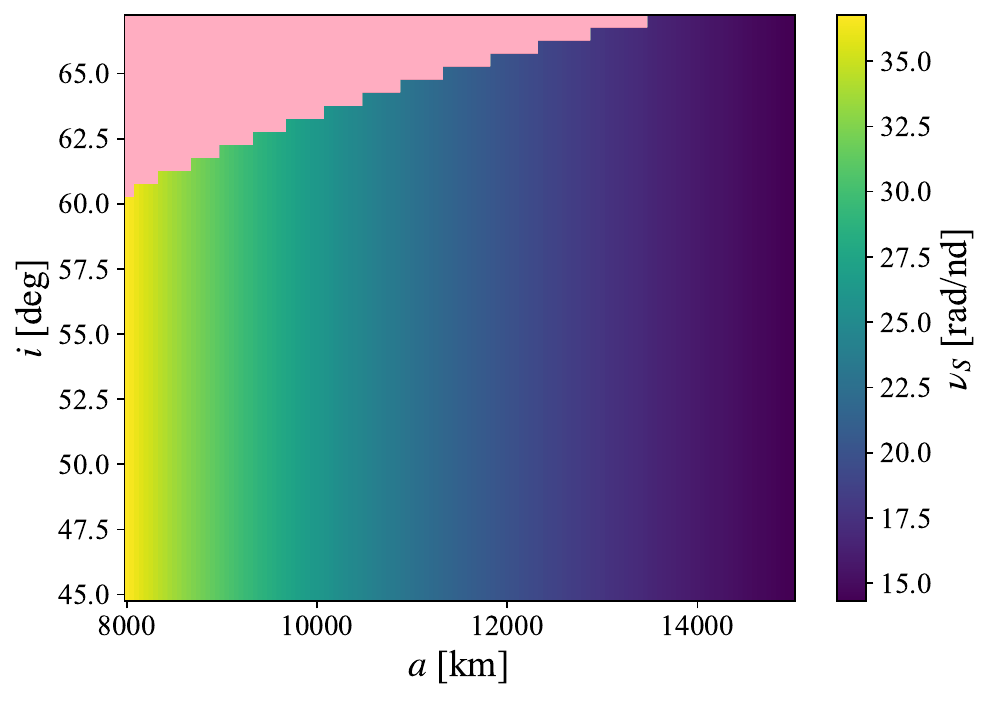}
        \caption{$\nu_S = \nu_S(a)$ from Eq.~\eqref{eq:nu_s}.}
        \label{fig:freq_grid_a}
    \end{subfigure}
    \hfill
    \begin{subfigure}[b]{0.48\textwidth}
        \centering
        \includegraphics[width=\linewidth]{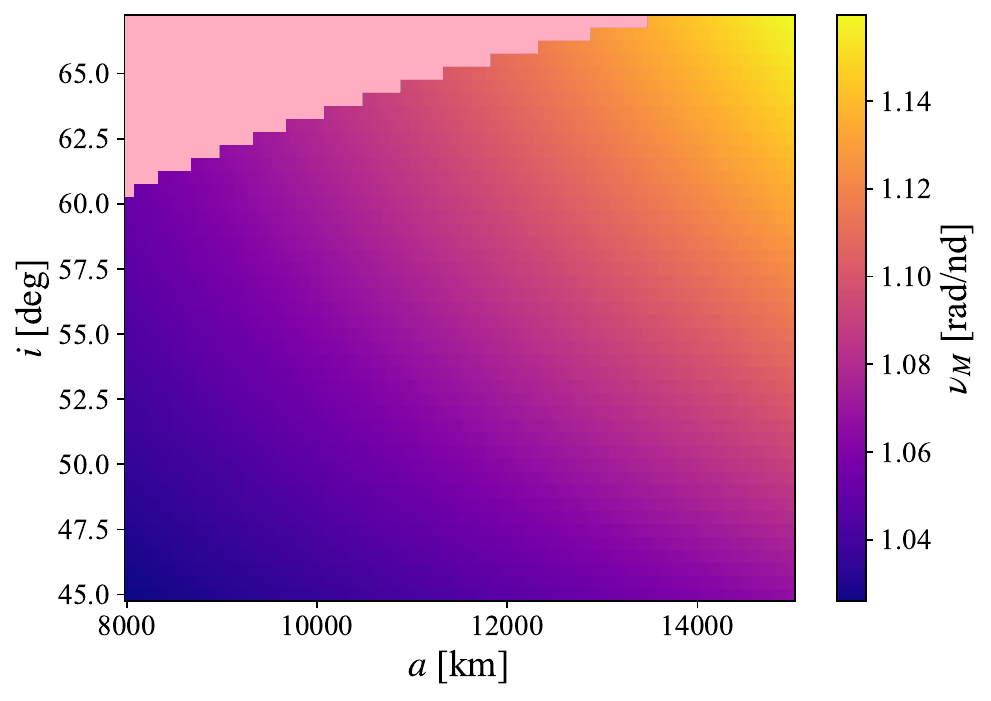}
        \caption{$\nu_M = \nu_M(a, i)$ from Eq.~\eqref{eq:nu_m}.}
        \label{fig:freq_grid_b}
    \end{subfigure}
    \caption{Frequency mapping over the $(a, i)$ design grid under the frozen-orbit condition (Eqs.~\eqref{eq:nu_s}-\eqref{eq:nu_m}). Pink regions indicate altitude violations (perilune below 100~km).}
    \label{fig:freq_grid}
\end{figure}

With the design variables and symmetries established in Sec.~\ref{subsec:design-vars}, the preliminary optimization procedure is posed as a global survey over a symmetry-reduced design space. At each node of an $(a, i)$ grid (with $50$~km and $0.5^\circ$ spacing in $a$ and $i$, respectively, over ranges encompassing the values specified by Brack et al.~\cite{brack2025}), the frequency pair $(\nu_S, \nu_M)$ is fixed through Eqs.~\eqref{eq:nu_s}-\eqref{eq:nu_m}, together with $e = \sqrt{1 - \tfrac{5}{3}\cos^2 i}$ and $\omega = \pi/2$; only the inter-satellite phase offsets $(\theta_{S,k}, \theta_{M,k})$ for $k = 2, \ldots, 5$ remain free. Figure~\ref{fig:freq_grid} illustrates the $(a, i) \to (\nu_S, \nu_M)$ mapping utilized for the survey. 

The design vector, $\bm{\xi}$, explicitly embeds the symmetries. With $\theta_{S,1} = \theta_{M,1} = 0$ (from the translational symmetry), the short-period phases for satellites 2 through 5 are parametrized as incremental offsets $\Delta\theta_{S,k} = \theta_{S,k} - \theta_{S,k-1} \geq 0$ for $k = 2, \ldots, 5$ (from the permutational symmetry). The design vector is then
\begin{align}
    \label{eq:design_vector_ref} \bm{\xi} = \{\Delta\theta_{S,2-5}, \theta_{M,2-5}\},
\end{align}
subject to $\sum_{k=2}^{5} \Delta\theta_{S,k} < 2\pi$; the medium-period phases $\theta_{M,k}$ are bounded within $[0, 2\pi)$ without ordering constraints. The objective is the fraction of the $(\theta_S, \theta_M) \in [0, 2\pi)^2$ torus satisfying the \acrshort{gdop} criterion (GDOP $\leq 6$) defined in Sec.~\ref{sec:assumed_constellation}, evaluated on a uniform $50 \times 50$ angular discretization. Leveraging the longitudinal symmetry, the user location is fixed within the \acrshort{mrf} at $\varphi = -83.32^\circ, \lambda = 0^\circ$. 

\subsection{Design Space Survey: User at the \acrshort{lsp}}
\label{sec:design-survey}
The design framework is applied to the \acrshort{dadm}, proceeding with the previous assumptions. The compact angular parametrization $(\theta_S, \theta_M) \in [0, 2\pi)^2$ admits exhaustive torus evaluation without long-horizon numerical propagation. The symmetry reduction developed in Sec.~\ref{subsec:design-framework} as well as the (nearly) closed-form evaluation of the \acrshort{elfo} states within the \acrshort{dadm} altogether render a global survey computationally tractable across the entire $(a, i)$ grid. A differential evolution algorithm (Python's \texttt{scipy.\allowbreak optimize.\allowbreak differential\_evolution}) serves as a proof-of-concept optimizer; a population of $100$ individuals is evolved over $200$ generations, and three independent trials are conducted per grid point to mitigate sensitivity to the initial population.\footnote{Differential evolution is used here as a proof-of-concept optimizer; the specific algorithm and its hyperparameters are not central to the framework.} The survey here is presented within the \acrshort{dadm}, with extensions to another time-independent model, i.e., the \acrshort{cr3bp}, supplied in Appendix~\ref{app:onaxis}.

\begin{figure}[h!]
    \centering
    \begin{subfigure}[b]{0.50\textwidth}
        \centering
        \includegraphics[width=\textwidth]{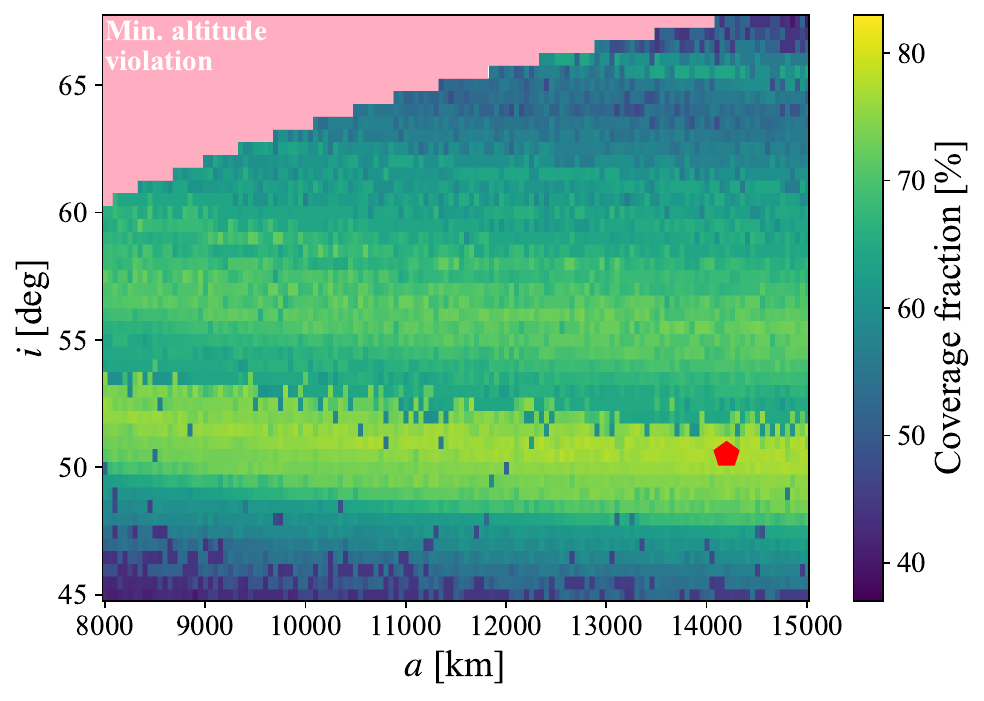}
        \caption{Coverage heatmap.}
        \label{fig:sweep_obl_heatmap}
    \end{subfigure}
    \hfill
    \begin{subfigure}[b]{0.40\textwidth}
        \centering
        \includegraphics[width=\textwidth]{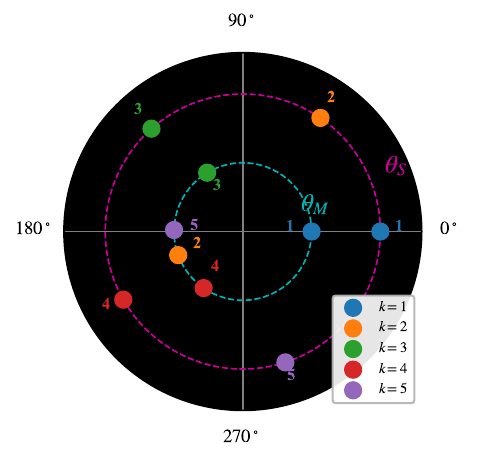}
        \caption{Phase distribution at \redpentagon.}
        \label{fig:sweep_obl_phases}
    \end{subfigure}
    \\[6pt]
    \begin{subfigure}[b]{0.47\textwidth}
        \centering
        \includegraphics[width=\textwidth]{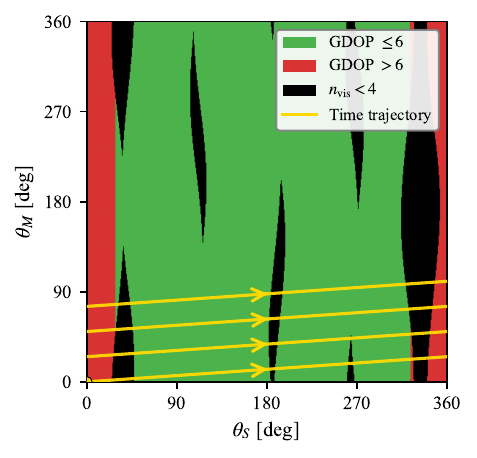}
        \caption{\acrshort{gdop} torus map at \redpentagon, ($\lambda = 0^\circ$, coverage $73.9\%$).}
        \label{fig:sweep_obl_torus}
    \end{subfigure}
    \hfill
    \begin{subfigure}[b]{0.47\textwidth}
        \centering
        \includegraphics[width=\textwidth]{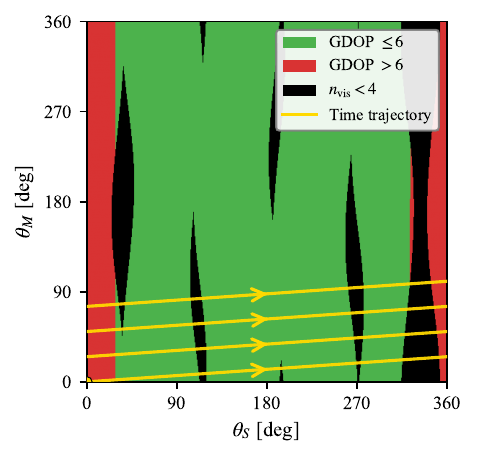}
        \caption{\acrshort{gdop} torus map at \redpentagon, ($\lambda = 180^\circ$, coverage $73.9\%$).}
        \label{fig:sweep_obl_torus_lon180}
    \end{subfigure}
    \caption{\acrshort{dadm} constellation design survey with the user at the \acrshort{lsp} ($\varphi = -83.32^\circ$ in the \acrshort{mrf}). The red pentagon (\redpentagon) in Fig.~\ref{fig:sweep_obl_heatmap} marks the reference configuration at $a = 14{,}200$~km, $i = 50.5^\circ$.}
    \label{fig:sweep_obl}
\end{figure}

The global design space is illustrated in Fig.~\ref{fig:sweep_obl}. The heatmap in Fig.~\ref{fig:sweep_obl_heatmap} colors each combination in $(a, i)$ with the best coverage fraction located via the optimization process. Note that two inclination bands near $i = 51^\circ$ and $i = 56^\circ$ supply favorable performance over a wide range of semi-major axis ($a$) values. The reference configuration at $(a, i) = (14{,}200~\text{km}, 50.5^\circ)$, employed throughout the current investigation, is highlighted with the red pentagon. Its (optimal) phase distribution is reported in Fig.~\ref{fig:sweep_obl_phases}; the inner and outer rings correspond to the medium- and short-period phases, respectively, with satellite~1 ({\textcolor[rgb]{0,0.447,0.741}{\raisebox{-0.25ex}{\Large$\bullet$}}}) fixed at the origin by the translational symmetry. The corresponding \acrshort{gdop} torus map (Fig.~\ref{fig:sweep_obl_torus}) is re-evaluated on a finer $500 \times 500$ grid; green regions indicate \acrshort{gdop} $\leq 6$, red regions indicate \acrshort{gdop} $> 6$, and black regions correspond to $n_{\mathrm{vis}} < 4$ where the \acrshort{gdop} is ill-defined. The underlying time evolution is recovered as a straight line on the torus starting from the origin with the slope $d\theta_M / d\theta_S = \nu_M / \nu_S$. This time revolution wraps around the torus and eventually fills it densely in the non-resonant case; four such wraps (yellow line) in $\nu_S$ are overlaid in Fig.~\ref{fig:sweep_obl_torus} to indicate the temporal direction. The longitudinal invariance is illustrated through symmetric behaviors in Figs.~\ref{fig:sweep_obl_torus}-\ref{fig:sweep_obl_torus_lon180}, evaluated with the same optimal phasing in Fig.~\ref{fig:sweep_obl_phases}. While Fig.~\ref{fig:sweep_obl_torus} is generated for the longitude $\lambda = 0^\circ$ (fixed within the \acrshort{mrf}), an identical coverage fraction ($73.9\%$) is observed at $\lambda = 180^\circ$ (also fixed within the \acrshort{mrf}) with the pattern shifted by $180^\circ$ in the $\theta_M$ direction. In general, $\lambda$ evolves with a non-resonant frequency $\omega_3$ and, thus, time-average of the coverage is accessed via a single reference value at $\lambda= 0^\circ$.  

\begin{figure}[h!]
    \centering
    \begin{subfigure}[b]{0.48\textwidth}
        \centering
        \includegraphics[width=\textwidth]{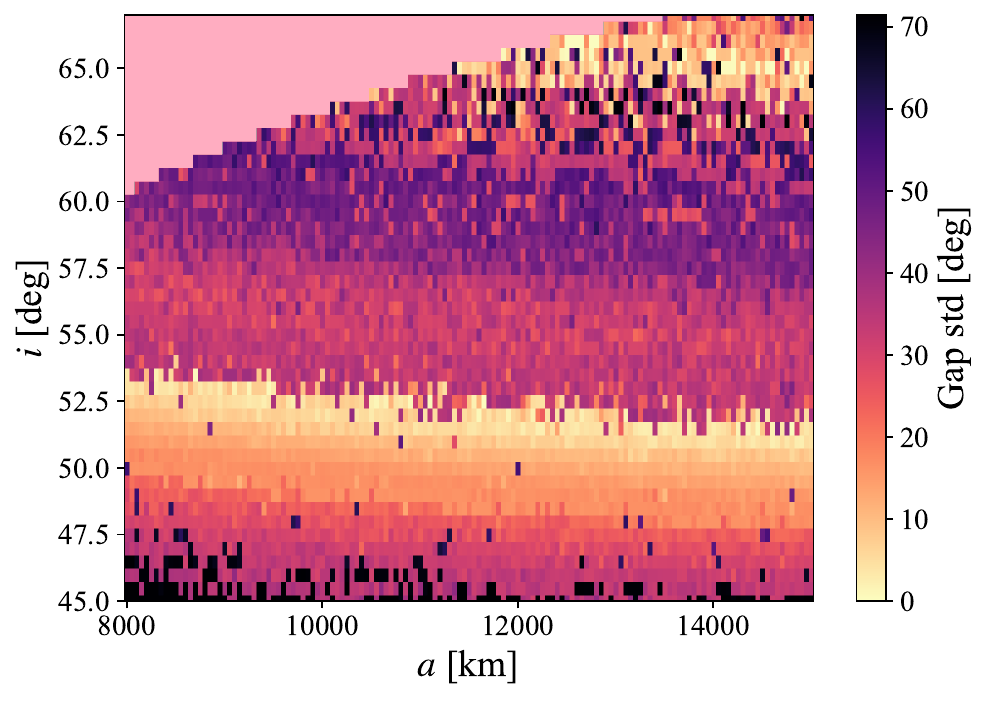}
        \caption{$\theta_S$.}
        \label{fig:spread_dadm_obl_thetaS}
    \end{subfigure}
    \hfill
    \begin{subfigure}[b]{0.48\textwidth}
        \centering
        \includegraphics[width=\textwidth]{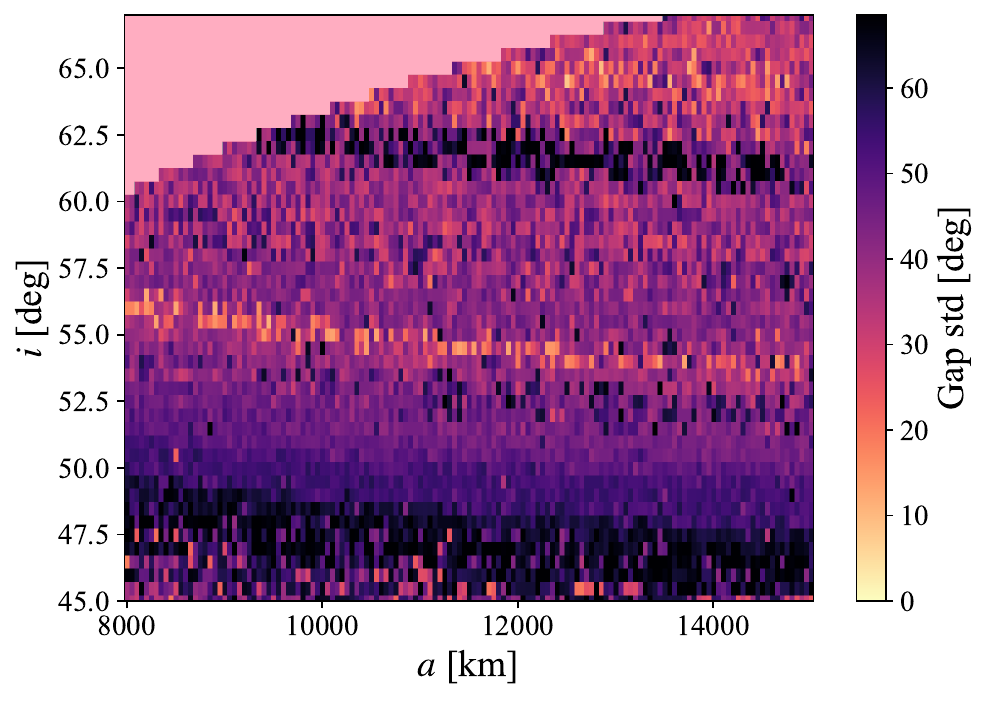}
        \caption{$\theta_M$.}
        \label{fig:spread_dadm_obl_thetaM}
    \end{subfigure}
    \caption{Standard deviation of consecutive angular gaps in the optimized \acrshort{dadm} phase distributions at the obliquity-corrected user latitude $\varphi = -83.32^\circ$ within the \acrshort{mrf}. (a) Short-period phase $\theta_S$; (b) medium-period phase $\theta_M$. Lower values (lighter) indicate more uniformly distributed satellites.}
    \label{fig:phase_spread_obl}
\end{figure}

The uniformity of the optimized phase distributions varies across the design space. Figure~\ref{fig:phase_spread_obl} displays the standard deviation for the five consecutive angular gaps among the satellite phases (sorted in ascending order on $[0, 2\pi)$ with the wrap-around gap included) across the $(a, i)$ grid for the respective optimal phases emerging in Fig.~\ref{fig:sweep_obl_heatmap}. Lower values correspond to more uniformly spaced satellites, with $\mathrm{std} = 0^\circ$ denoting a perfectly uniform five-satellite constellation with $\approx 72^\circ$ offsets in the respective phases. In Fig.~\ref{fig:spread_dadm_obl_thetaS}, the high-uniformity $\theta_S$ band near $i \approx 51^\circ$ is depicted. At higher inclinations near $i \approx 56^\circ$, the optimizer rather converges to markedly non-uniform $\theta_S$ spacing, consistent with an optimized configuration reported by \citet{brack2025}. This trend appears characteristic of five-satellite constellations at higher inclinations; see Appendix~\ref{app:uniform-phase} for an additional examination of the spacing in $\theta_S$. The distribution in $\theta_M$ is illustrated in Fig.~\ref{fig:spread_dadm_obl_thetaM}. It remains comparatively non-uniform across the design space. Therefore, unless operational constraints require a specific number of orbital planes, a diverse $\theta_M$ distribution may enhance overall constellation coverage. 

The global design space survey reveals useful patterns and insights regarding the favorable inclination bands (Fig.~\ref{fig:sweep_obl_heatmap}) as well as the phase spacing patterns (Fig.~\ref{fig:phase_spread_obl}). These insights facilitate preliminary design and adaptations in higher-level mission design serving as warm starting points for subsequent optimization schemes. Recognition of the frequency structure of the \acrshort{elfo} and the constellation is therefore central to a tractable global survey of the design space. The patterns reported here are tied to the adopted threshold-based \acrshort{gdop} metric and the five-satellite, \acrshort{lsp}-targeting configuration; alternative architectures or performance objectives may reshape the favored bands and phase distributions. Indeed, the possibility of such adaptations reinforces the value of the present tractable global survey as a preliminary tool.

\subsection{Long-Period Amplitude ($A_L$) on the Coverage Metric: a Case Study}
\label{subsec:al-dadm}
This subsection revisits the assumption to enforce $A_L = 0$ among the constellation satellites. The long-period oscillation introduces an additional degree of freedom per satellite: the amplitude $A_{L,k}$ and phase $\theta_{L,k}$, altogether governing the coupled triad $(e, i, \omega)$ (Fig.~\ref{fig:coupling}). Each satellite $k$ may possess an independent amplitude $A_{L,k}$ and phase $\theta_{L,k}$, so the constellation state (or \acrshort{gdop}) no longer lies on the two-dimensional $(\theta_S,\theta_M)$ torus exploited in the preceding subsections. For distinct $A_{L,k} \neq 0$ ($\nu_{L,k} \neq 0$ are also distinct), the \acrshort{gdop} becomes a quasi-periodic function with at least $7 = 2+5$ independent frequencies. Angle-domain sampling is therefore less effective in general, motivating a finite-horizon, time-domain evaluation as an alternative.

The reference constellation at $a = 14{,}200$~km, $i = 50.5^\circ$, $e \approx 0.571$ from Fig.~\ref{fig:sweep_obl_phases} is re-assessed with potential non-zero $A_{L,k}$ components as a case study. First, the angular domain solution from Fig.~\ref{fig:sweep_obl_phases} is re-optimized within the time domain. The state is sampled at a cadence of $T_S/(24\sqrt{2}) \approx 1.24$~hr, selected to be incommensurate with the short-period $T_S$ and thereby avoid aliasing artifacts. Each candidate configuration is optimized over an $\approx 2\,T_L$ horizon ($332$~d)\footnote{The libration period $T_L$ is constructed for each sampled $A_L$ via the singularity-free quadrature derived in Appendix~\ref{app:lp-quadrature}.} using differential evolution, and the optimum is then re-evaluated over a $20\,T_L$ horizon ($\approx 9$~yr) to ensure that the reported coverage is not an artifact of the shorter optimization window. The baseline yields a coverage of $74.3\%$ over the re-evaluation horizon and serves as the reference point for all subsequent analyses.

A fair comparison against the frozen equilibrium ($A_L = 0$) baseline requires that the inter-satellite phasing $(\theta_{S,k}, \theta_{M,k})$ not be held fixed at its $A_L = 0$ optimum when $A_L \neq 0$ is introduced. As such, the corresponding design vector is $\bm{\xi} \in \mathbb{R}^{18}$, consisting of $(\Delta \theta_{S,2-5}, \theta_{M,2-5}, A_{L,1-5}, \theta_{L,1-5})$. The osculating initial conditions for \acrshort{dadm} propagation are constructed from $\bm{\xi}$ via the three-step procedure detailed in Appendix~\ref{app:lp-ics}. Since the new design vector is higher-dimensional as compared to Eq.~\eqref{eq:design_vector_ref}, three numerical experiments are examined in sequence, each isolating a specific aspect of the long-period descriptors. In the first (P1), a uniform amplitude is imposed across all five satellites, $A_{L,k} = A_L$ for $k = 1,\ldots,5$, and the thirteen-dimensional inter-satellite design vector $(\Delta \theta_{S,2-5}, \theta_{M,2-5}, \theta_{L,1-5})$ is optimized jointly via differential evolution. Five different levels of amplitude are examined at $A_L \in \{0.01, 0.02, 0.03, 0.04, 0.05\}$. The second experiment (P2) examines a more realistic scenario where the satellites carry different long-period amplitudes rather than a common value. The mixed configuration $A_{L,k} = [0, 0.01, 0.02, 0.03, 0.04]$ is adopted as an arbitrary case, and the same thirteen-dimensional inter-satellite design vector is jointly optimized. The third experiment (P3) represents the most generic scenario, where the amplitudes are released as free variables alongside the phases, yielding an eighteen-dimensional joint optimization with $A_{L,k} \in [0, 0.05]$.

The results are consolidated in Table~\ref{tab:lp_experiments}. For P1, the jointly optimized coverage fails to recover the equilibrium baseline at any nonzero amplitude, and the deficit grows monotonically with $A_L$, from $-1.6$~pp at $A_L = 0.01$ to $-6.4$~pp at $A_L = 0.05$. For P2, the jointly optimized coverage of $70.3\%$ remains $4.1$~pp below the baseline, consistent with the monotone loss observed in P1. In both experiments, re-distributing the inter-satellite phasing (in all $\theta_S, \theta_M, \theta_L$) modulates the deficit but cannot eliminate it. In P3, if any combination of nonzero amplitude and inter-satellite phasing improves coverage relative to the equilibrium configuration, the differential evolution process is expected to locate such a solution. The optimizer rather drives all amplitudes toward zero (mean $A_{L,k} = 0.003$, maximum $0.006$), and the resulting coverage of $74.1\%$ nearly matches, but still falls short of, the equilibrium baseline within $0.3$~pp. Within this tested five-satellite, \acrshort{lsp}-targeting configuration under a threshold-based \acrshort{gdop} metric ($\leq 6$), the long-period amplitude is therefore detrimental: neither inter-satellite phase tuning nor a jointly optimized amplitude configuration is preferred to $A_L = 0$.

The degradation originates in the interaction between the \acrshort{gdop} distribution and the threshold-based metric. Across the reference constellation, the well-conditioned regions cluster in a narrow band just below the threshold, whereas the ill-conditioned regions extend into a hyperbolic tail far above the threshold as the visibility geometry approaches singular ($n_{\mathrm{vis}}<4$). A small perturbation induced by the long-period oscillation therefore flips well-conditioned epochs across the threshold far more readily than it rescues ill-conditioned ones, and the good-to-bad transitions dominate the reverse direction by a large margin. This structural asymmetry leads to net coverage loss for $A_L > 0$. The scope of this conclusion is worth stating explicitly. The $A_L$ dependency may differ under alternative constellation architectures, mission objectives, or performance metrics; however, the lower-fidelity analysis developed here supplies a systematic strategy for isolating and examining the impact of $A_L$ in any such setting. The present finding motivates the suppression of $A_L$ for the reference constellation within the \acrshort{hfem} as well, examined in Sec.~\ref{sec:fddc}.

\begin{table}[htb]
    \centering
    \caption{Coverage at the \acrshort{lsp} for the three long-period ($A_L\neq 0$) experiments (P1, P2, and P3), evaluated over $20\,T_L$ in the \acrshort{dadm}.}
    \label{tab:lp_experiments}
    \begin{tabular}{l l c c}
        \toprule
        Experiment & Configuration & Jointly optimized [\%] & vs.\ baseline [pp] \\
        \midrule
        Baseline & $A_L = 0$ & $74.3$ & - \\
        \midrule
        P1 & $A_L = 0.01$  & $72.8$ & $-1.6$ \\
        P1 & $A_L = 0.02$           & $72.3$ & $-2.0$ \\
        P1 & $A_L = 0.03$           & $69.8$ & $-4.5$ \\
        P1 & $A_L = 0.04$           & $69.2$ & $-5.2$ \\
        P1 & $A_L = 0.05$           & $67.9$ & $-6.4$ \\
        \midrule
        P2 & $A_{L,k} = \{0,\,0.01,\,0.02,\,0.03,\,0.04\}$ & $70.3$ & $-4.1$ \\
        \midrule
        P3 & $A_{L,k} \in [0,\,0.05]$ & $74.1$ & $-0.3$ \\
        \bottomrule
    \end{tabular}
\end{table}

\section{Refinement and Analysis in a Higher-Fidelity Model}
\label{sec:fddc}

The lower-fidelity analysis in Sec.~\ref{sec:design-exploration} supplies a time-independent examination of constellation performance through a coverage metric embedded within the $(\theta_S, \theta_M)$ torus. Deployment in a mission-specific environment, however, requires re-examination of this performance in a higher-fidelity dynamical model, where reducing \acrshort{gdop} to a lower-dimensional quasi-periodic function is challenging. A typical refinement approach expands the design vector to include the full state of each satellite at a common initial epoch, $\bm{\xi} = (\bm{X}_1, \ldots, \bm{X}_5) \in \mathbb{R}^{6\cdot 5}$, and optimizes over a finite time horizon with numerical integrations in the loop. Beyond the computational cost of repeated higher-fidelity propagations, such a search operates on a high-dimensional design space without explicit handles on the frequency components, e.g., $A_L$, that impact coverage.

The analysis in this section alternatively pursues three frequency-domain routes, each exploiting the spectral structure established in Sec.~\ref{sec:freq-structure}. First, the \acrfull{fddc}~\cite{park2025frequency} refines the initial conditions via targeting a vector of differentiable frequency components, decoupling the per-satellite correction from any joint optimization of all five satellites (Sec.~\ref{subsec:fddc-demo}). Second, a Fourier surrogate constructed from an \acrshort{hfem} trajectory enables phase re-optimization without costly numerical propagation in the loop (Sec.~\ref{subsec:fourier-reopt}). Third, the residual coverage gap is examined in the frequency domain and attributed to forcing components that no initial-condition refinement can remove (Sec.~\ref{subsec:hf-gravity-impact}). Common to all three is a single structural premise: under the symmetry-inducing assumptions of Sec.~\ref{subsec:symmetry-assumptions}, the five satellites lie on a common invariant torus within the \acrshort{hfem}. The \acrshort{fddc} drives them onto that torus with prescribed coordinates, the Fourier surrogate evaluates coverage anywhere on it without propagation, and the residual analysis isolates the forced content that no choice of torus coordinates removes. These three distinct directions confirm the usefulness of the proposed torus-based, frequency-domain strategies in both refining and interpreting the \acrshort{hfem} \acrshort{elfo} constellations. Coverage is evaluated as the daily fraction of samples with GDOP $\leq 6$ at the \acrshort{lsp}, a finite-horizon metric consistent with~\citet{brack2025}; comparisons against the torus-averaged coverage of Sec.~\ref{sec:design-exploration} are therefore interpreted as indicators rather than identical statistics.

\subsection{Targeting Specific Frequency Structures via FDDC}
\label{subsec:fddc-demo}
Prior work~\cite{park2025frequency} introduces the \acrshort{fddc} as a model-agnostic mechanism that drives a quasi-periodic trajectory toward a prescribed set of spectral components; the \acrshort{fddc} is now applied to the \acrshort{elfo} constellation, with the target spectral components encoding the inter-satellite design constraints associated with the lower-fidelity solution from Sec.~\ref{sec:design-exploration}. A satellite's initial state $\freeVE_k \in \mathbb{R}^6$\footnote{$\freeVE_k$ consists of $\bm{R}_J, \bm{V}_J$ in the current analysis, although a different parameterization, e.g., osculating orbital elements, may be leveraged.} is the free variable, and a vector constraint $\constraintV{k}(\freeVE_k) = \bm{0}$ is assembled from frequency components extracted from the propagated trajectory in the target dynamical model (\acrshort{hfem}),
\begin{align}
    \label{eq:fddc_constraint}
    \constraintV{k}(\freeVE_k) =
    \mb \nu_{S, k}(\freeVE_k) - \nu_{S}^* \\
        \theta_{S,k}(\freeVE_k) - \theta_{S,k}^* \\
        \nu_{M, k}(\freeVE_k) - \nu_{M}^* \\
        \theta_{M,k}(\freeVE_k) - \theta_{M,k}^* \\
        A_{L,k}(\freeVE_k) - A_{L}^*
    \me = \bm{0}, \quad k = 1, \ldots, 5,
\end{align}
where the target frequencies $(\nu_S^*, \nu_M^*)$ and phase angles $(\theta_{S,k}^*, \theta_{M,k}^*)$ are derived from the lower-fidelity solution. These first four components preserve the inter-satellite phasing across the fidelity transition. The fifth component is engaged when suppressing the long-period amplitude to a small value, i.e., $A_L^*\ll 1$. Although the long-period phase angle $\theta_L$ may also be targeted, it is unnecessary and ill-defined for $A_L \ll 1$. As such, $\theta_L$ is omitted from the target frequency structure.

Following the representative scalars in the last column of Table~\ref{tab:freq_summary}, each frequency component is extracted from a dedicated scalar signal sampled along the propagated trajectory: the short-period pair $(\nu_{S,k}, \theta_{S,k})$ from the out-of-plane component $z$ within the \acrshort{mrf}, the medium-period pair $(\nu_{M,k}, \theta_{M,k})$ from the in-plane component $x$ in the \acrshort{mrf}, and the long-period amplitude $A_{L,k}$ from the osculating eccentricity $e$. The targeted peak is identified, respectively, as the first (dominant) peak of the $z$ spectrum, the peak of the $x$ spectrum near $\nu_M$, and the peak of the $e$ spectrum nearest the long-period frequency $\nu_L$ prescribed from the lower-fidelity analysis (Sec.~\ref{sec:freq-structure}, Appendix~\ref{app:lp-quadrature}).

Peak locations from the \acrfull{dft} of the signal are refined to sub-bin accuracy using one of two complementary strategies available within the \acrshort{fddc} framework: iterative refinement in the spirit of \citet{laskar1999introduction} that maximizes the magnitude of the continuous spectrum at the peak, and a collocation-based approach following \citet{gomez2010collocation} that enforces exact DFT values at the peak bin and an adjacent bin. Both strategies yield \emph{differentiable} dependencies of $\nu_j$, $A_j$, and $\theta_j$ on the underlying signal, and therefore on the initial state through the propagated dynamics, i.e., the \acrshort{hfem}. The sensitivity matrix $\partial \constraintV{k}/\partial \freeVE_k$ is assembled by chain rule through the state transition matrix: for a sample at time $t_i$, the contribution to each spectral component is mediated by the state transition matrix $\stm{}{}{t_i}{t_0}$, and the full Jacobian of a frequency-phase-amplitude triplet relative to $\freeVE_k$ reduces to a $3 \times 6$ block per targeted peak irrespective of the overall phase-space dimension of the quasi-periodic trajectory~\cite{park2025frequency}. The resulting Newton update is
\begin{equation}
    \label{eq:fddc_newton}
    \Delta \freeVE_k = -\eta\,\left(\tfrac{\partial \constraintV{k}}{\partial \freeVE_k}\right)^{\dagger}\,\constraintV{k}(\freeVE_k),
\end{equation}
with attenuation factor $\eta \in (0, 1]$ and $(\cdot)^{\dagger}$ denoting the Moore-Penrose pseudoinverse. The update is applied independently per satellite. For the full derivation and broader algorithmic details, see~\citet{park2025frequency}.

Iteration terminates when $\|\constraintV{k}\|_\infty$ drops below a prescribed tolerance ($10^{-6}$ in the present work, expressed in normalized frequency and radian phase units). For the reference five-satellite configuration in the \acrshort{hfem}, a two-stage procedure is adopted for robust convergence behavior. The first stage refines $(\nu_S, \theta_S, \nu_M, \theta_M)$ over a $1.5$-year arc with $\eta = 1$ (typically less than $10$ iterations), and the second stage additionally engages $A_L$ over a $5$-year arc with $\eta = 0.5$ (typically less than $20$ iterations). Accurate detection of $A_L$ requires a horizon on the order of years; the spectral resolution inversely scales with the propagation horizon and $\nu_L < 0.2$~rad/nd is a relatively small value.

\begin{figure}[h!]
    \centering
    \begin{subfigure}[b]{0.48\textwidth}
        \centering
        \includegraphics[width=\textwidth]{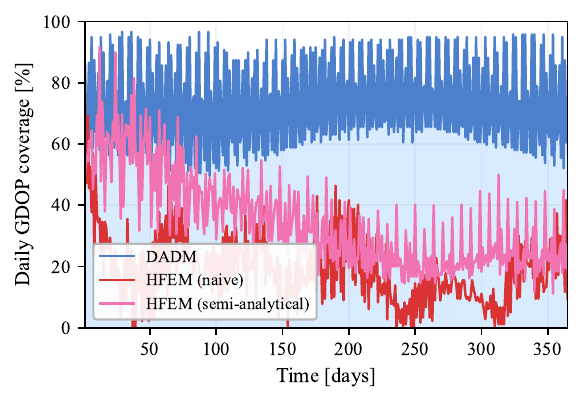}
        \caption{\acrshort{gdop}, na\"{\i}ve and semi-analytical transformation.}
        \label{fig:fddc_gdop_before}
    \end{subfigure}
    \hfill
    \begin{subfigure}[b]{0.48\textwidth}
        \centering
        \includegraphics[width=\textwidth]{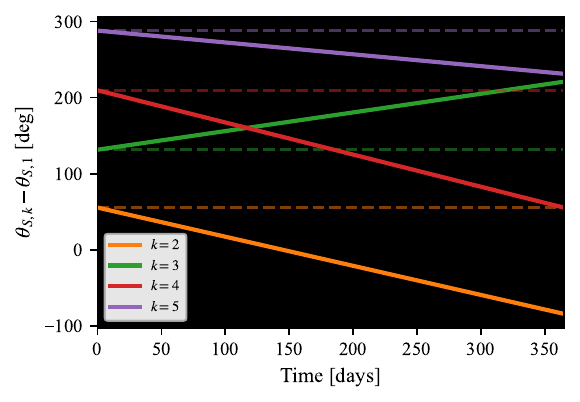}
        \caption{Short-period phase drift before \acrshort{fddc} (na\"{\i}ve).}
        \label{fig:fddc_phase_before}
    \end{subfigure}
    \\[6pt]
    \begin{subfigure}[b]{0.48\textwidth}
        \centering
        \includegraphics[width=\textwidth]{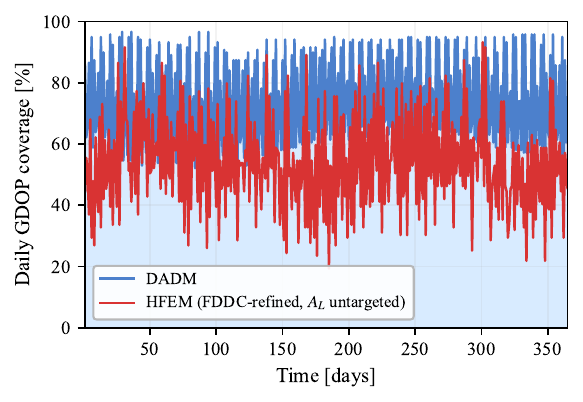}
        \caption{\acrshort{gdop}, \acrshort{fddc}-refined: $A_L$ untargeted.}
        \label{fig:fddc_gdop_after}
    \end{subfigure}
    \hfill
    \begin{subfigure}[b]{0.48\textwidth}
        \centering
        \includegraphics[width=\textwidth]{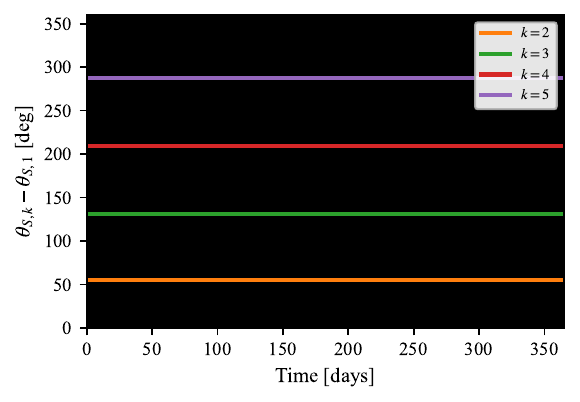}
        \caption{Short-period phase drift after \acrshort{fddc}.}
        \label{fig:fddc_phase_after}
    \end{subfigure}
    \\[6pt]
    \begin{subfigure}[b]{0.48\textwidth}
        \centering
        \includegraphics[width=\textwidth]{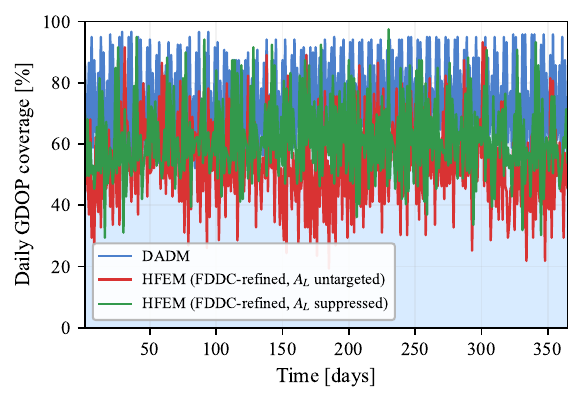}
        \caption{\acrshort{gdop}, \acrshort{fddc}-refined: $A_L$ untargeted vs.\ suppressed.}
        \label{fig:fddc_AL_gdop}
    \end{subfigure}
    \hfill
    \begin{subfigure}[b]{0.48\textwidth}
        \centering
        \includegraphics[width=\textwidth]{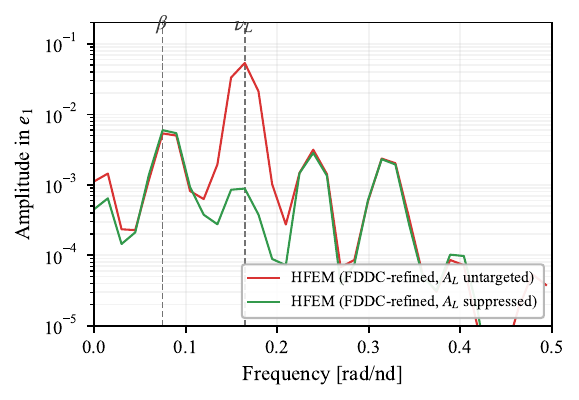}
        \caption{Satellite 1 eccentricity spectrum.}
        \label{fig:fddc_AL_fft}
    \end{subfigure}
    \caption{Refinement within the \acrshort{hfem} through the \acrshort{fddc}.}
    \label{fig:fddc}
\end{figure}

The value of a \acrshort{fddc} refinement step is illustrated through two transformation strategies: (1) na\"{\i}ve transformation, where the mean elements from the \acrshort{dadm} are assumed equivalent to the osculating elements at the starting epoch, $JD_0 = 2460941.5$ (2025-09-23 00:00:00 UTC), and (2) a semi-analytical transformation based on the Von Zeipel method~\cite{nie2018lunar}, correcting for first-order error. The analysis continues with the same constellation as originally designed within the \acrshort{dadm}, illustrated in Figs.~\ref{fig:sweep_obl_phases}-\ref{fig:sweep_obl_torus}. These two strategies are compared in Table~\ref{tab:vz_comparison} in terms of their short-period frequency $\nu_S$; the mismatch in $\nu_S$ is the greatest source of \acrshort{gdop} degradation over a short time span. Over a one-year evaluation horizon ($\approx 84$~nd), the na\"{\i}ve transformation, with a worst-case inter-satellite spread of $\max|\Delta\nu_S| = 0.05$~rad/nd, accumulates a relative phase drift of up to $4.2$~rad ($\approx 0.67$ full cycles) across the constellation, misaligning the inter-satellite phasing. The semi-analytical transformation~\cite{nie2018lunar} reduces this spread to $\max|\Delta\nu_S| = 0.03$~rad/nd and raises the coverage to $34.9\%$, yet the residual mismatch continues to misalign the inter-satellite phasing over the same horizon. The \acrshort{gdop} performance is illustrated in Fig.~\ref{fig:fddc_gdop_before}, where the blue line is the reference constellation from the \acrshort{dadm} analysis and the red line is the constellation initialized from the na\"{\i}ve transformation; the vertical axis reports the daily fraction of samples satisfying \acrshort{gdop} $\leq 6$. The short-period phase drift ($\theta_{S,k} - \theta_{S,1}$, $2 \leq k \leq 5$) is examined in Fig.~\ref{fig:fddc_phase_before}, where dashed and solid lines denote the lower-fidelity design values and \acrshort{hfem} results, respectively. The underlying cause for the degrading performance in Fig.~\ref{fig:fddc_gdop_before} is apparent, where the mismatched $\nu_S$ in Table~\ref{tab:vz_comparison} results in a short-period phase drift that misaligns the inter-satellite phasing over the one-year horizon. 

The \acrshort{fddc} eliminates the phase drift in two stages, each examined below against the baseline from the na\"{\i}ve transformation: the first restores the inter-satellite phasing, and the second additionally suppresses the long-period amplitude. The first stage targets only the short- and medium-period frequencies and phases (Fig.~\ref{fig:fddc}\subref{fig:fddc_gdop_after}-\subref{fig:fddc_phase_after}). The short-period phase drift is eliminated over the horizon (Fig.~\ref{fig:fddc_phase_after}), leading to more sustained coverage over the time span in Fig.~\ref{fig:fddc_gdop_after}. Table~\ref{tab:vz_comparison} includes the refined frequency structure at $(\nu_S^*, \nu_M^*)$. However, the coverage ($54.1\%$) still falls short of the \acrshort{dadm} baseline ($73.5\%$). This gap is partly overcome via the second \acrshort{fddc} stage by suppressing the long-period motion with the target $A_L = A_L^*$. In the current analysis, $A_L^* = 10^{-3}$. Below this value, numerical detection of the frequency component within the \acrshort{fddc} process becomes increasingly difficult\footnote{As seen in Fig.~\ref{fig:fddc_AL_fft}, $\nu_L$ is generally close to $2\beta$, challenging highly accurate frequency-component extraction in general.}. The results are depicted in Fig.~\ref{fig:fddc_AL_gdop} and Fig.~\ref{fig:fddc_AL_fft}: the $A_L$-suppressed constellation outperforms the $A_L$-untargeted constellation (red: $A_L$ untargeted; green: $A_L$ suppressed). In the eccentricity spectrum, the peak near $\nu_L \approx 0.159$ is effectively suppressed to the designated value ($A_L^* = 10^{-3}$) for all five satellites (Fig.~\ref{fig:fddc_AL_fft}). This improvement is consistent with the result from Sec.~\ref{subsec:al-dadm}, where a non-zero $A_L$ degrades constellation coverage for this reference solution. After this second stage, the satellite frequency components are targeted as desired, illustrated in the last column of Table~\ref{tab:vz_comparison}. On top of the common values for $(\nu_S^*, \nu_M^*)$, $A_L = A_L^*$ is enforced across all satellites.

\begin{table}[htb]
    \centering
    \caption{Per-satellite frequency-amplitude diagnostics evaluated in the \acrshort{hfem} (5-year propagation, frequency-collocation refinement). Design values: $\nu_S^* = 15.55$~rad/nd, $\nu_M^* = 1.085$~rad/nd; the frozen-orbit reference has $A_L = 0$. Full initial conditions are supplied in Appendix~\ref{app:initial-conditions}.}
    \label{tab:vz_comparison}
    \begin{tabular}{cl cccc}
        \toprule
        Frequency & Satellite & Na\"{\i}ve & Semi-analytical~\cite{nie2018lunar} & \acrshort{fddc} & \acrshort{fddc} ($A_L^* = 10^{-3}$) \\
        \midrule
        \multirow{6}{*}{$\nu_S$~[rad/nd]}
            & 1 & 15.57 & 15.53 & 15.55 & 15.55 \\
            & 2 & 15.54 & 15.52 & 15.55 & 15.55 \\
            & 3 & 15.59 & 15.54 & 15.55 & 15.55 \\
            & 4 & 15.54 & 15.55 & 15.55 & 15.55 \\
            & 5 & 15.56 & 15.53 & 15.55 & 15.55 \\
            & $\max|\Delta\nu_S|$ & 0.05 & 0.03 & $\approx 0$ & $\approx 0$ \\
        \midrule
        \multirow{6}{*}{$\nu_M$~[rad/nd]}
            & 1 & 1.084 & 1.083 & 1.085 & 1.085 \\
            & 2 & 1.082 & 1.082 & 1.085 & 1.085 \\
            & 3 & 1.082 & 1.082 & 1.085 & 1.085 \\
            & 4 & 1.082 & 1.081 & 1.085 & 1.085 \\
            & 5 & 1.082 & 1.082 & 1.085 & 1.085 \\
            & $\max|\Delta\nu_M|$ & 0.002 & 0.002 & $\approx 0$ & $\approx 0$ \\
        \midrule
        \multirow{6}{*}{$A_L$}
            & 1 & 0.071 & 0.039 & 0.054 & 0.001 \\
            & 2 & 0.053 & 0.026 & 0.049 & 0.001 \\
            & 3 & 0.022 & 0.033 & 0.013 & 0.001 \\
            & 4 & 0.028 & 0.031 & 0.026 & 0.001 \\
            & 5 & 0.053 & 0.023 & 0.050 & 0.001 \\
            & $\mathrm{mean}\,A_L$ & 0.045 & 0.030 & 0.038 & 0.001 \\
        \midrule
        \multicolumn{2}{l}{Coverage (365~d)} & 19.2\% & 34.9\% & 54.1\% & 61.7\% \\
        \multicolumn{2}{l}{} & (Fig.~\ref{fig:fddc_gdop_before}) & & (Fig.~\ref{fig:fddc_gdop_after}) & (Fig.~\ref{fig:fddc_AL_gdop}) \\
        \bottomrule
    \end{tabular}
\end{table}

The correction of inter-satellite phasing (i.e., targeting $\theta_{S,k}$ and $\theta_{M,k}$) is, in principle, achievable through a differential correction process within the orbital element domain. For example, linear fits to the osculating mean anomaly $M$ or ascending node $\Omega$ supply the slope (frequency) and intercept (phase). An iterative process can then be devised to adjust the initial orbital elements until the desired phase offsets are achieved. For example, the semi-major axis correction pursued by \citet{ely2006constellations} as well as Ceresoli et al.~\cite{ceresoli2025design} represents an approach for targeting $\nu_S$. \citet{ceresoli2025design} further introduce a station-keeping strategy based on $M$ control boxes. The \acrshort{fddc} offers practical advantages beyond this conceptual similarity. First, suppression of the long-period amplitude $A_L$ has no one-to-one direct counterpart in the orbital elements. It requires control over the magnitude of the coupled oscillation in $(e, i, \omega)$, a non-trivial process within the \acrshort{hfem}. Moreover, the \acrshort{fddc} targets multiple components $(\bm{\nu}, \bm{\theta})$ in a single coupled correction with well-defined gradient information, all aiding rapid, automated refinement within realistic models. While illustrated for a point-mass based \acrshort{hfem}, the strategy extends to any other conservative forcing functions, e.g., the lunar harmonics.

The results suggest three practical roles for the \acrshort{fddc} within an \acrshort{elfo} constellation design workflow: (1) generating warm-start seeds for subsequent higher-fidelity optimization, (2) transitioning lower-fidelity solutions into baseline trajectories in the target dynamics, and (3) refining solutions optimized over a finite horizon such that the targeted frequency structure is preserved beyond the optimization window. The third role is relevant even for configurations already optimized within a higher-fidelity model: without an explicit frequency-component handle, integration-based searches can leave residual spectral mismatches that accumulate into phase drift outside the optimization horizon.

\subsection{Phase Re-Optimization via a Fourier Surrogate}
\label{subsec:fourier-reopt}

The persisting gap between the \acrshort{hfem} and the \acrshort{dadm} coverage after targeting the frequency structure, as illustrated in Fig.~\ref{fig:fddc_AL_gdop}, partly originates from the different optimality domains; \acrshort{dadm}-optimal phases may be displaced from the locally ``best'' configuration under the \acrshort{hfem}. Although a direct optimization with $\bm{\xi} \in \mathbb{R}^{30}$ is possible to bridge this gap, costly numerical propagations without explicit constraints to bind the solution to a common $(\nu_S, \nu_M)$ with $A_L = 0$ remains a key disadvantage.

A surrogate model offers an alternative route, constructed from the Fourier decomposition of a single reference satellite trajectory within the \acrshort{hfem}. Under the constellation assumptions introduced in Sec.~\ref{subsec:symmetry-assumptions} (shared internal frequencies $\nu_S$ and $\nu_M$ across all five satellites, with $A_L = 0$), the satellites lie on a common quasi-periodic torus within the \acrshort{hfem} \cite{jorba1997persistence}, differing only in their torus coordinates $(\theta_{S,k}, \theta_{M,k})$. Recall that the external frequencies and phases $(\bm{\omega}, \bm{\alpha})$ (Eq.~\eqref{eq:omega}) are fixed for a given epoch and shared among the constellation. A Fourier decomposition of a reference satellite orbit therefore captures this common torus, and the remaining four satellites are recovered through per-peak phase rotations keyed to their relative offsets $(\Delta\theta_{S,k}, \Delta\theta_{M,k})$, leveraging the multi-frequency spectral characterization of \acrshort{hfem} trajectories on quasi-periodic tori~\cite{lian2013note,park2025numerical}. The construction proceeds in three steps: (1) a reference trajectory is decomposed into a finite spectral sum; (2) a per-peak phase rotation generates the remaining satellites on the same quasi-periodic torus; (3) the resulting analytical positions are evaluated at arbitrary times to evaluate the \acrshort{gdop}, forgoing an explicit numerical integration within the \acrshort{hfem}. The full coefficient tables for the reference decomposition are reported in Appendix~\ref{app:fourier}; the method summary required for the re-optimization and accuracy assessment is supplied here.

The reference trajectory is the \acrshort{hfem} \acrshort{mrf} position history of satellite~1 over a 10-year horizon, with the initial condition extracted from the \acrshort{fddc}-refined constellation of Sec.~\ref{subsec:fddc-demo} (with $A_L$ suppressed; last column in Table~\ref{tab:vz_comparison}). Each Cartesian component is decomposed via the Laskar refinement algorithm~\cite{laskar1999introduction} into $N_\zeta$ sub-bin-refined peaks, ordered by descending amplitude, each carrying an amplitude $A_{j,\zeta}$, a phase $\phi_{j,\zeta}$, and an integer label $(m_j, n_j, p_j, q_j, s_j)$ identifying it as a linear combination of the fundamental frequencies $(\nu_S, \nu_M, \omega_1, \omega_2, \omega_3)$\footnote{While other frequencies exist \cite{gomez2002solar}, they are not needed in decomposing the relatively high-amplitude peaks.} introduced in Sec.~\ref{subsec:freq-across-models}. The decomposition captures the positional variance with the first $N_x = 40$, $N_y = 40$, and $N_z = 18$ peaks. 

\begin{figure}[h!]
    \centering
    \includegraphics[width=0.65\textwidth]{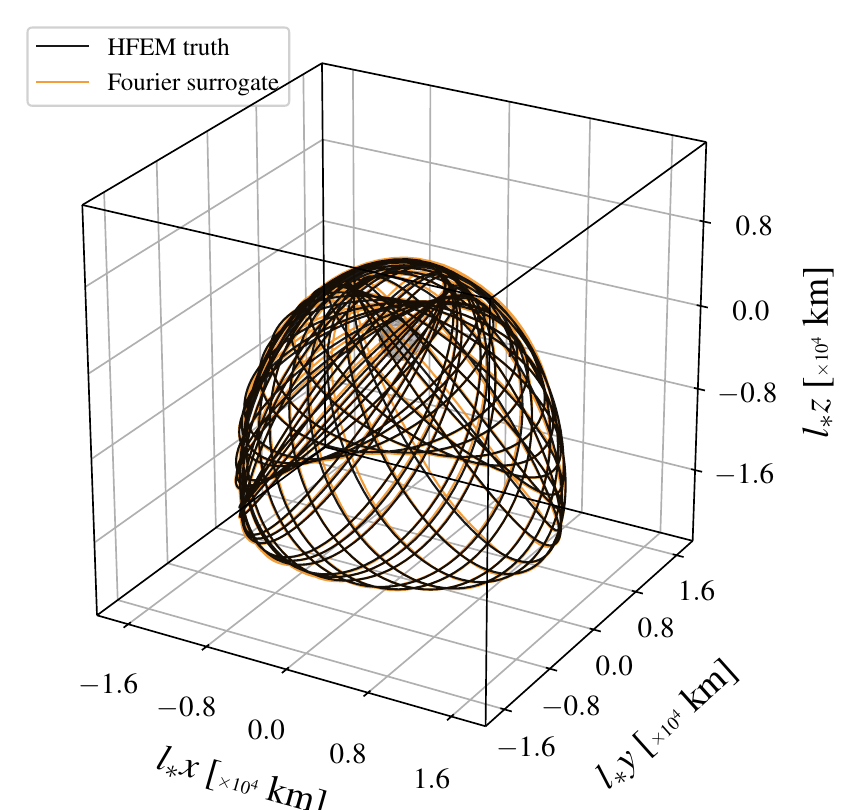}
    \caption{Comparison of the Fourier surrogate trajectory (orange) against the \acrshort{hfem} truth (black) in the \acrshort{mrf} over two months, satellite 1 (with $A_L$ suppressed; last column in Table~\ref{tab:vz_comparison}).}
    \label{fig:fourier-accuracy}
\end{figure}

The remaining four satellites follow from the reference decomposition by a phase rotation alone, and the resulting surrogate reproduces the \acrshort{hfem} coverage closely enough to replace direct propagation. Given the decomposition of satellite~1, the dimensional position history of any satellite $k$ on the same quasi-periodic torus is supplied via a per-peak phase rotation,
\begin{align}\label{eq:fourier-shift}
    l_*r_{k,\zeta}(t) = l_*\bar{r}_\zeta + \sum_{j=1}^{N_\zeta} A_{j,\zeta}\,\cos\,\!\bigl(\mathfrak{f}_{j,\zeta}\, t + \phi_{j,\zeta} + m_j\,\Delta\theta_{S,k} + n_j\,\Delta\theta_{M,k}\bigr), \quad \zeta \in \{x, y, z\},
\end{align}
Here, $l_*$ is the characteristic length scale (Sec.~\ref{sec:prelim}), $(\Delta\theta_{S,k}, \Delta\theta_{M,k})$ are the torus phase offsets of satellite~$k$ relative to satellite~1, and $\bar{r}_\zeta$ is the time-averaged nondimensional position (constant offset). With this surrogate, a full five-satellite \acrshort{gdop} evaluation results in a speedup of several orders of magnitude relative to \acrshort{hfem} numerical propagation with a trade-off in accuracy. The surrogate fidelity is illustrated directly against the \acrshort{hfem} truth trajectory used to extract the decomposition. An overlay in Fig.~\ref{fig:fourier-accuracy} illustrates the two trajectories for satellite~1 in the \acrshort{mrf} over a two-month window, with a mean position difference of $\approx 180$~km. This offset is small relative to the orbit scale ($a \approx 14{,}200$~km), and the surrogate retains the overall geometric structure of the truth trajectory. As a complementary check at the \acrshort{gdop} level, the surrogate is evaluated at the reference \acrshort{fddc}-refined, $A_L$-suppressed phasing against the \acrshort{hfem} truth (Fig.~\ref{fig:fddc_AL_gdop}): the surrogate yields $61.9\%$ against the truth value of $61.7\%$, a discrepancy of only $0.2$~pp. This gap, attributable to the finite spectral resolution of the truncated decomposition, is small relative to the optimization gains sought and confirms that the surrogate rather faithfully represents the \acrshort{hfem} coverage landscape across the torus.

\begin{figure}[h!]
    \centering
    \includegraphics[width=0.48\textwidth]{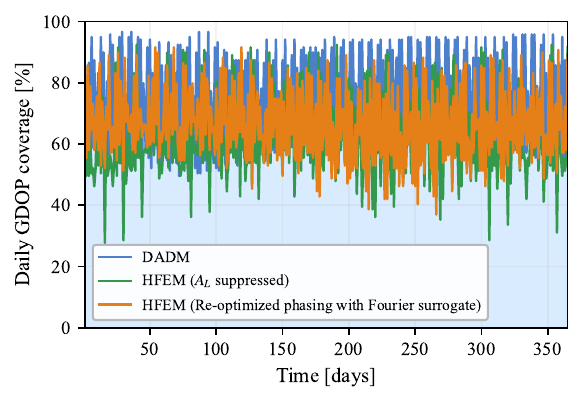}
    \caption{Daily \acrshort{gdop} coverage for the reference constellation before and after Fourier-surrogate phase re-optimization in the \acrshort{hfem}. Blue and green plots are identical to those in Fig.~\ref{fig:fddc_AL_gdop}.}
    \label{fig:fourier-reopt-gdop}
\end{figure}

With the surrogate in place, the phasing is re-optimized via differential evolution over the eight free parameters $(\Delta\theta_{S,2-5},\allowbreak \Delta\theta_{M,2-5})$, with satellite~1 fixed at the origin. The population is seeded with the \acrshort{dadm}-optimal phasing to ensure that the optimizer originates within the same basin of attraction. Across three independent seeds, with a population of 60 and 100 generations, the surrogate model identifies an optimum coverage of $68.4\%$, an improvement of $+6.5$~pp over the same surrogate evaluated at the baseline phasing ($61.9\%$, reported in the preceding paragraph). The optimized phasing is transferred to the \acrshort{hfem} initial conditions through the same \acrshort{fddc} infrastructure utilized in Sec.~\ref{sec:fddc}. The phase targets $(\theta_{S,k}^*, \theta_{M,k}^*)$ identified by the surrogate are imposed as constraints, alongside the $A_L$ suppression target at $A_L^* = 10^{-3}$. After convergence, the resulting five-satellite constellation is propagated in the \acrshort{hfem} and the coverage is evaluated. The outcome is summarized in Table~\ref{tab:fourier-reopt}, and, in Fig.~\ref{fig:fourier-reopt-gdop}, the daily coverage is illustrated as a time history. The re-optimized, \acrshort{fddc}-refined constellation achieves a coverage of $67.2\%$, an improvement of $+5.5$~pp over the baseline phasing. The surrogate prediction of $68.4\%$ is not fully realized as the Fourier surrogate uses a finite series. Nevertheless, the improvement is notable and is achieved solely through phasing adjustment, with no change to the overall orbital geometry $(\nu_S, \nu_M)$.

The two optimal phasings located independently in the \acrshort{dadm} and via the \acrshort{hfem}-Fourier surrogate are cross-evaluated in Table~\ref{tab:fourier-reopt}. The phase distributions differ by at most $11^\circ$ per satellite (Table~\ref{tab:fourier-reopt-phases}), yet the coverage outcomes vary by up to $19.3$~pp depending on the dynamical model used for evaluation (Table~\ref{tab:fourier-reopt-cov}). In each row of Table~\ref{tab:fourier-reopt-cov}, the boldface entry marks the model where that phasing was optimized; it is also the higher value in the row. The coverage landscape is thus model-specific: optimizing within one dynamical model does not guarantee improvement in another, even when the phase distributions are close. Nonetheless, the existence of an \acrshort{hfem} optimum in the neighborhood of the \acrshort{dadm} solution suggests the usefulness of warm-start strategies, including the lower-fidelity analysis (Sec.~\ref{sec:design-exploration}) and the Fourier surrogate in the current section.

\begin{table}[htb]
    \centering
    \caption{Comparison of two phasing strategies for the reference constellation ($a = 14{,}200$~km, $i = 50.5^\circ$, $\nu_S \approx 15.55$~rad/nd, $\nu_M \approx 1.08$~rad/nd).}
    \label{tab:fourier-reopt}
    \begin{subtable}[t]{\textwidth}
        \centering
        \caption{Optimized phase distributions [deg] (satellite~1 fixed at the origin).}
        \label{tab:fourier-reopt-phases}
        \begin{tabular}{c r@{\hskip 1em}r r r@{\hskip 1em}r}
            \toprule
            & \multicolumn{2}{c}{\acrshort{dadm}} & & \multicolumn{2}{c}{HFEM-Fourier} \\
            \cmidrule(lr){2-3} \cmidrule(lr){5-6}
            Sat & $\theta_S$ & $\theta_M$ & & $\theta_S$ & $\theta_M$ \\
            \midrule
            1 &   0.0 &   0.0 & &   0.0 &   0.0 \\
            2 &  55.7 & 200.1 & &  51.8 & 189.5 \\
            3 & 131.6 & 121.2 & & 132.4 & 118.3 \\
            4 & 209.7 & 235.3 & & 214.3 & 240.1 \\
            5 & 288.0 & 178.7 & & 279.1 & 172.0 \\
            \bottomrule
        \end{tabular}
    \end{subtable}

    \vspace{8pt}
    \begin{subtable}[t]{\textwidth}
        \centering
        \caption{Coverage [\%] within dynamical models evaluated on the optimal phasing from Table~\ref{tab:fourier-reopt-phases} ($365$-day propagation).}
        \label{tab:fourier-reopt-cov}
        \begin{tabular}{l cc}
            \toprule
            Evaluation model & \acrshort{dadm} phasing & HFEM-Fourier surrogate phasing \\
            \midrule
            \acrshort{dadm}  & \textbf{73.5}\textsuperscript{$\ddagger$} & 54.2 \\
            \acrshort{hfem}  & 61.7 & \textbf{67.2} \\
            \bottomrule
        \end{tabular}
    \end{subtable}

    \vspace{6pt}
    \footnotesize
    \textsuperscript{$\ddagger$}\,The 365-day \acrshort{dadm} coverage of $73.5\%$ is not identical to the torus-averaged $73.9\%$ in Fig.~\ref{fig:sweep_obl_torus} (Sec.~\ref{sec:design-survey}) or to the re-optimized time-domain solution $74.3\%$ in Sec.~\ref{subsec:al-dadm}, evaluated over $20T_L$. The resulting differences merely reflect the inherent sensitivity of threshold-based coverage to evaluation domain and horizon, and comparisons within a fixed procedure (e.g., the four cells of Table~\ref{tab:fourier-reopt-cov}) remain meaningful.
\end{table}

\subsection{Unremovable Coverage Oscillations in the \acrshort{hfem}}
\label{subsec:hf-gravity-impact}

Frequency structures impact the constellation performance, producing oscillations that cannot be entirely removed in the \acrshort{hfem}. Specifically, \acrshort{elfos} deform under additional perturbations that are absent from the \acrshort{dadm} (Sec.~\ref{subsec:freq-across-models}). In contrast to the long-period amplitude ($A_L$), these components arise from physical effects that the \acrshort{dadm} averaging procedure removes, and they cannot be eliminated by initial-condition adjustment within the \acrshort{hfem}. The residual gap between the \acrshort{hfem}-Fourier-optimized constellation and the reference \acrshort{dadm} constellation (Fig.~\ref{fig:fourier-reopt-gdop}) substantiates this behavior\footnote{The degradation of \acrshort{cr3bp} constellation coverage in comparison with the \acrshort{dadm} (Figs.~\ref{fig:sweep_heatmap} and \ref{fig:cr3bp_sweep_heatmap}) also illustrates the impact on coverage from dynamics, examined in Appendix~\ref{app:onaxis}.}. This subsection identifies the dominant contributors and characterizes their role in the residual gap. 

\begin{figure}[h!]
    \centering
    \begin{subfigure}[b]{0.48\textwidth}
        \centering
        \includegraphics[width=\textwidth]{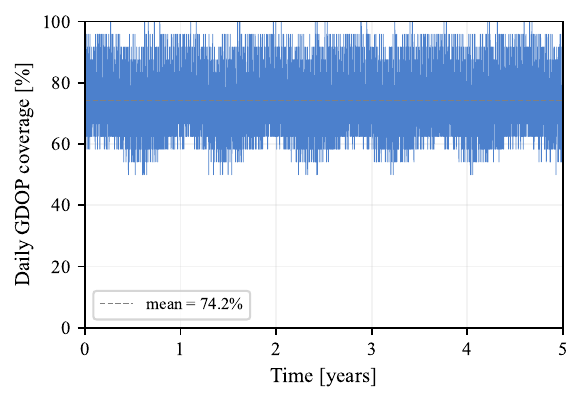}
        \caption{Daily \acrshort{gdop} coverage, \acrshort{dadm}.}
        \label{fig:gdop-fft-5yr-dadm-ts}
    \end{subfigure}
    \hfill
    \begin{subfigure}[b]{0.48\textwidth}
        \centering
        \includegraphics[width=\textwidth]{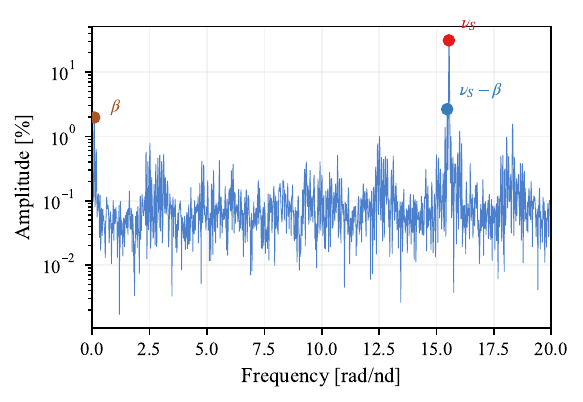}
        \caption{Spectrum of quarter-day coverage, \acrshort{dadm}.}
        \label{fig:gdop-fft-5yr-dadm-fft}
    \end{subfigure}
    \\[6pt]
    \begin{subfigure}[b]{0.48\textwidth}
        \centering
        \includegraphics[width=\textwidth]{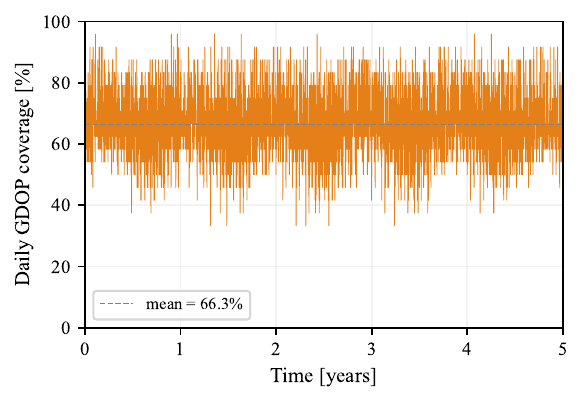}
        \caption{Daily \acrshort{gdop} coverage, \acrshort{hfem}.}
        \label{fig:gdop-fft-5yr-ts}
    \end{subfigure}
    \hfill
    \begin{subfigure}[b]{0.48\textwidth}
        \centering
        \includegraphics[width=\textwidth]{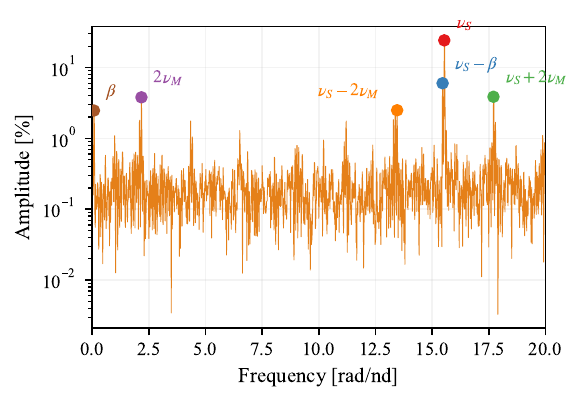}
        \caption{Spectrum of quarter-day coverage, \acrshort{hfem}.}
        \label{fig:gdop-fft-5yr-fft}
    \end{subfigure}
    \caption{Five-year \acrshort{gdop} coverage and spectrum for the reference constellation: \acrshort{dadm} vs \acrshort{hfem}.}
    \label{fig:gdop-fft-5yr}
\end{figure}

The coverage examined in Fig.~\ref{fig:gdop-fft-5yr-dadm-ts} is originally identified from Figs.~\ref{fig:sweep_obl_torus} and \ref{fig:fourier-reopt-gdop}. The same blue plot from Fig.~\ref{fig:fourier-reopt-gdop} is extended to a 5-year evaluation. The oscillations in the daily coverage are further examined in the frequency domain in Fig.~\ref{fig:gdop-fft-5yr-dadm-fft}\footnote{Note, the raw \acrshort{gdop} signal is not smooth and does not constitute a proper signal for frequency analysis. A quarter-day (averaged) coverage is used as the signal to reduce aliasing and allow detection of $\nu_S$.}, with the first three peaks marked with their linear frequency decomposition. They consist of $\nu_S, \nu_S-\beta$, and $\beta$. The first frequency $\nu_S$ is intuitively understood from the pattern in Fig.~\ref{fig:sweep_obl_torus} that mainly evolves with the short-period. Then, the $\beta$-modulation follows from the evolution of user geometry within the \acrshort{mrf} as discussed in Sec.~\ref{subsec:obliquity}.

In contrast to the \acrshort{dadm} constellation, the \acrshort{hfem}-Fourier (and \acrshort{fddc}-refined) solution is examined for a 5-year horizon in Figure~\ref{fig:gdop-fft-5yr-ts}. The oscillations in the daily coverage are further examined in the frequency domain in Fig.~\ref{fig:gdop-fft-5yr-fft}, with the six largest peaks and their linear-combination labels summarized in Table~\ref{tab:gdop-peaks}. The last column indicates the existence of the corresponding peak within the \acrshort{dadm}; the peaks $\nu_S, \nu_S - \beta$, and $\beta$ are matched, whereas the other peaks that involve $2\nu_M$ occur only specifically to the \acrshort{hfem} since they originate from the unaveraged Earth forcing functions absent from the \acrshort{dadm} (Sec.~\ref{sec:freq-structure}). Some of the $\beta$-modulation is common across the models, reflecting the $\omega_3$ wobble of the user in the \acrshort{mrf} (Sec.~\ref{subsec:obliquity}). However, the amplitudes are larger for the \acrshort{hfem} for both peaks at $\nu_S - \beta$ and $\beta$. The larger oscillations are expected to originate from the geometry change of the \acrshort{elfos} under the \acrshort{hfem} forcing effects, also modulated at the $\beta$ frequency (Sec.~\ref{subsec:freq-across-models}). As a whole, while both the \acrshort{dadm} and \acrshort{hfem} demonstrate long-term ($\beta$, $\sim$ 1 year) oscillations in the coverage, these oscillations are much larger for the \acrshort{hfem}. The \acrshort{hfem} also introduces $2\nu_M$-related oscillations that modulate the performance over shorter timescales. 

\begin{table}[h!]
    \centering
    \caption{Six largest peaks in the \acrshort{hfem} \acrshort{gdop} spectrum (Fig.~\ref{fig:gdop-fft-5yr-fft}) and the corresponding \acrshort{dadm} amplitudes (Fig.~\ref{fig:gdop-fft-5yr-dadm-fft}). $\checkmark$ marks peaks present in both spectra; $\times$ marks peaks absent from the \acrshort{dadm}.}
    \label{tab:gdop-peaks}
    \begin{tabular}{c c r r l r r}
        \toprule
        \multirow{2}{*}{No.} & Marker & Freq. & Period & \multirow{2}{*}{Mode} & \acrshort{hfem} amp & \acrshort{dadm} amp \\
                             & (Figs.~\ref{fig:gdop-fft-5yr-dadm-fft},\,\ref{fig:gdop-fft-5yr-fft}) & [rad/nd] & [d] &  & [\%] & [\%] \\
        \midrule
        1 & \pkmark{peakA} & $\approx 15.54$ & $\approx 1.76$  & $\nu_S$          & $\approx 24.1$ & $\checkmark\ \approx 31.1$ \\
        2 & \pkmark{peakB} & $\approx 15.47$ & $\approx 1.77$  & $\nu_S - \beta$  & $\approx 6.0$  & $\checkmark\ \approx 2.6$  \\
        3 & \pkmark{peakC} & $\approx 17.71$ & $\approx 1.54$  & $\nu_S + 2\nu_M$ & $\approx 3.9$  & $\times$                   \\
        4 & \pkmark{peakD} & $\approx\phantom{0}2.17$  & $\approx 12.6$  & $2\nu_M$         & $\approx 3.8$  & $\times$                   \\
        5 & \pkmark{peakE} & $\approx 13.46$ & $\approx 2.03$  & $\nu_S - 2\nu_M$ & $\approx 2.5$  & $\times$                   \\
        6 & \pkmark{peakF} & $\approx\phantom{0}0.08$  & $\approx 342$   & $\beta$          & $\approx 2.5$  & $\checkmark\ \approx 2.0$  \\
        \bottomrule
    \end{tabular}
\end{table}

These unremovable oscillatory behaviors in coverage have multiple implications. First, controllable and uncontrollable (unremovable) components coexist in a given \acrshort{elfo} constellation. Although the \acrshort{elfo} geometry, expressed through $\nu_S, \nu_M, A_L$ and the phases $\theta_S, \theta_M, \theta_L$, is controlled by the initial satellite states, the geometry deformations from realistic perturbations are inherent to the dynamics and inherited from the epoch through Eq.~\eqref{eq:omega}. The resulting classification is predictive for design within the investigated reference configuration: the controllable components flag where initial condition and phasing refinement may reduce coverage loss, whereas the unremovable components identify the residual floor that lies beyond initial-condition design. The largest such variation originates from the beat frequency $\beta$, together with its $\nu_S - \beta$ sideband, resulting in a non-negligible modulation for approximately a one-year period. This long-period modulation is detrimental to short-window optimization: if the optimization horizon is shorter than $\beta^{-1}$, the optimizer may settle on a configuration whose phase of the $\beta$ envelope happens to be favorable within that window, selecting a locally high-coverage interval at the cost of degraded performance as the envelope evolves over the remainder of the mission horizon. Removing these modulations requires active maneuver strategies that remain out of scope; without such interventions, the periodic degradations in coverage are natural and must be expected (Fig.~\ref{fig:gdop-fft-5yr-ts}). Although the current analysis focuses on a single configuration, the $\beta$-deformation arises from the \acrshort{hfem} dynamics and may affect other constellation configurations to different degrees depending on the specific performance metric and user geometries.

\section{Concluding Remarks}
\label{sec:conclusions}

The current investigation examines \acrshort{elfo} constellation design and analysis through a torus-based, frequency-domain lens, treating each \acrshort{elfo} as an invariant torus in the underlying dynamics. Within the \acrshort{dadm}, each \acrshort{elfo} is characterized by three intrinsic frequency-angle pairs, and their evolution across models of increasing fidelity is predictable; while the intrinsic structure persists, additional spectral components arise from the higher-fidelity perturbations. The same perspective extends to the satellite-user geometry, where lunar obliquity injects modulations at the lunar draconic period (Sec.~\ref{subsec:obliquity}). Two complementary roles emerge from this fundamental spectral perspective within the context of \acrshort{elfo} constellation design and analysis.

The first role is a rapid, fidelity-bridging framework for constellation design. The angular parametrization $(\theta_S, \theta_M)$ collapses an otherwise unbounded time horizon into a compact two-dimensional torus, and three symmetries (translational, permutational, and longitudinal) further reduce the effective design space. Together, these features enable an epoch-free global survey within time-independent, lower-fidelity models. The \acrshort{fddc} transitions designs across model fidelity by simultaneously targeting the short- and medium-period frequencies, their phases, and the long-period amplitude $A_L$; a Fourier surrogate constructed from a single higher-fidelity reference trajectory then re-optimizes the inter-satellite phasing within the target dynamics without propagation in the loop. These components carry a constellation design from the \acrshort{dadm} to the \acrshort{hfem} without entering a high-dimensional, integration-in-the-loop optimization. Several of these ingredients individually draw from prior work; the present contribution is their integration into a coherent fidelity-bridging workflow.

The second role is an interpretive analysis of the constellation coverage behavior. The spectral analysis separates the coverage of a higher-fidelity constellation into design-refinable and unremovable components. Phase drift in $(\theta_S, \theta_M)$ is reduced by the \acrshort{fddc} in Sec.~\ref{subsec:fddc-demo}; non-zero $A_L$ is demonstrated to degrade the adopted \acrshort{gdop} coverage metric in Sec.~\ref{subsec:al-dadm} and is suppressed in the higher-fidelity transfer through the \acrshort{fddc}; residual phasing sub-optimality is reduced by the \acrshort{hfem} Fourier surrogate in Sec.~\ref{subsec:fourier-reopt}. Once these controllable effects are addressed, Sec.~\ref{subsec:hf-gravity-impact} attributes the residual gap primarily to the unremovable oscillations originating from the unaveraged Earth forcing terms and the beat frequency mechanism ($\beta$). Among the unremovable components, the beat frequency $\beta = \nu_M - \omega_3$, arising from the differential nodal regression between the spacecraft and the Earth-Moon plane, is identified as the dominant ${\sim}1$-year modulation of coverage in the \acrshort{hfem}, a structural feature not previously emphasized in \acrshort{elfo} constellation literature. This breakdown supplies a structural account of coverage behavior and a diagnostic lens for related \acrshort{elfo} constellation design studies.

Together, the rapid, fidelity-bridging framework and the diagnostic decomposition complement existing approaches to \acrshort{elfo} constellation design. Although established here on a reference five-satellite, \acrshort{lsp}-targeting configuration under a threshold-based \acrshort{gdop} metric and a point-mass \acrshort{hfem}, the perspective points toward an adaptive, rapid design framework for \acrshort{elfo} constellations, with broader extensions across architectures, coverage metrics, and additional dynamical perturbations as a natural next step.

\section*{Acknowledgments}
Dr. Beom Park acknowledges support from the Apollo 11 Postdoctoral Fellowship from the School of Aeronautics and Astronautics at Purdue University. Portions of this work were completed at Purdue University under Intuitive Machines LLC Award 21123283.

\bibliographystyle{unsrtnat}
\bibliography{references}

\begin{thebibliography}{32}
\providecommand{\natexlab}[1]{#1}
\providecommand{\url}[1]{\texttt{#1}}
\expandafter\ifx\csname urlstyle\endcsname\relax
  \providecommand{\doi}[1]{doi: #1}\else
  \providecommand{\doi}{doi: \begingroup \urlstyle{rm}\Url}\fi

\bibitem[von Zeipel(1910)]{von1910application}
Hugo von Zeipel.
\newblock Sur l'application des s{\'e}ries de {M}. {L}indstedt {\`a}
  l'{\'e}tude du mouvement des com{\`e}tes p{\'e}riodiques.
\newblock \emph{Astronomische Nachrichten}, 183:\penalty0 345--418, 1910.
\newblock \doi{10.1002/asna.19091832202}.

\bibitem[Lidov(1962)]{lidov1962evolution}
Michael~L. Lidov.
\newblock The evolution of orbits of artificial satellites of planets under the
  action of gravitational perturbations of external bodies.
\newblock \emph{Planetary and Space Science}, 9\penalty0 (10):\penalty0
  719--759, 1962.
\newblock \doi{10.1016/0032-0633(62)90129-0}.

\bibitem[Kozai(1959)]{kozai1959effects}
Yoshihide Kozai.
\newblock On the effects of the sun and the moon upon the motion of a close
  earth satellite.
\newblock \emph{SAO Special Report}, 22\penalty0 (part 2), 1959.

\bibitem[Ely(2005)]{ely2005stable}
Todd~A. Ely.
\newblock Stable constellations of frozen elliptical inclined lunar orbits.
\newblock \emph{The Journal of the Astronautical Sciences}, 53:\penalty0
  301--316, 2005.
\newblock \doi{10.1007/BF03546355}.

\bibitem[Ely and Lieb(2006)]{ely2006constellations}
Todd~A. Ely and Erica Lieb.
\newblock Constellations of elliptical inclined lunar orbits providing polar
  and global coverage.
\newblock \emph{The Journal of the Astronautical Sciences}, 54\penalty0
  (1):\penalty0 53--67, 2006.
\newblock \doi{10.1007/BF03256476}.

\bibitem[Folta and Quinn(2006)]{folta2006lunar}
David Folta and David Quinn.
\newblock Lunar frozen orbits.
\newblock In \emph{AIAA/AAS Astrodynamics Specialist Conference and Exhibit},
  page 6749, 2006.
\newblock \doi{10.2514/6.2006-6749}.

\bibitem[Israel et~al.(2020)Israel, Mauldin, Roberts, Mitchell, Pulkkinen,
  Cooper, Johnson, Christe, and Gramling]{israel2020lunanet}
David~J Israel, Kendall~D Mauldin, Christopher~J Roberts, Jason~W Mitchell,
  Antti~A Pulkkinen, La~Vida~D Cooper, Michael~A Johnson, Steven~D Christe, and
  Cheryl~J Gramling.
\newblock Lunanet: a flexible and extensible lunar exploration communications
  and navigation infrastructure.
\newblock In \emph{2020 IEEE Aerospace Conference}, pages 1--14. IEEE, 2020.
\newblock \doi{10.1109/AERO47225.2020.9172509}.

\bibitem[Jones(2024)]{jones2024queqiao}
Andrew Jones.
\newblock China's {Queqiao-2} relay satellite enters lunar orbit.
\newblock SpaceNews, 2024.
\newblock URL
  \url{https://spacenews.com/chinas-queqiao-2-relay-satellite-enters-lunar-orbit/}.
\newblock Accessed: 2026-03-04.

\bibitem[Ceresoli et~al.(2025)Ceresoli, Maccari, Moresi, and
  Lavagna]{ceresoli2025design}
Michele Ceresoli, Fabrizio Maccari, Edoardo Moresi, and Mich{\`e}le Lavagna.
\newblock Design and station-keeping strategies for robust lunar navigation
  constellations.
\newblock \emph{Acta Astronautica}, 236:\penalty0 20--31, 2025.
\newblock \doi{10.1016/j.actaastro.2025.06.025}.

\bibitem[Zanotti et~al.(2024)Zanotti, Ceresoli, Pasquale, Prinetto, and
  Lavagna]{zanotti2024high}
Giovanni Zanotti, Michele Ceresoli, Andrea Pasquale, Jacopo Prinetto, and
  Mich{\`e}le Lavagna.
\newblock High performance lunar constellation for navigation services to moon
  orbiting users.
\newblock \emph{Advances in Space Research}, 73\penalty0 (11):\penalty0
  5665--5679, 2024.
\newblock \doi{10.1016/j.asr.2023.03.032}.

\bibitem[Brack et~al.(2025)Brack, Roorda, Mathur, Gaylor, Hartigan, Ryden,
  Crenshaw, Volle, and Stewart]{brack2025}
Daniel Brack, Tim Roorda, Ravishankar Mathur, David Gaylor, Mark Hartigan,
  Grant Ryden, Juan Crenshaw, Michael Volle, and Shaun Stewart.
\newblock Preliminary design of the lunar data network constellation under
  operational considerations.
\newblock In \emph{2025 AAS/AIAA Astrodynamics Specialist Conference, Boston,
  MA, August 10-14}. 2025.

\bibitem[Park et~al.(2025{\natexlab{a}})Park, Sanaga, and
  Howell]{park2025numerical}
Beom Park, Rohith~Reddy Sanaga, and Kathleen~C. Howell.
\newblock Numerical assessment of a frequency-based hierarchy for the cislunar
  domain.
\newblock \emph{Journal of Guidance, Control, and Dynamics}, 48\penalty0
  (11):\penalty0 2462--2479, 2025{\natexlab{a}}.
\newblock \doi{10.2514/1.G009176}.

\bibitem[Park and Howell(2026)]{park2026bridging}
Beom Park and Kathleen~C. Howell.
\newblock Linking averaged and unaveraged three-body dynamics near smaller
  primaries: Symmetric periodic orbits.
\newblock \emph{arXiv preprint arXiv:2606.08485}, 2026.
\newblock \doi{10.48550/arXiv.2606.08485}.

\bibitem[Park et~al.(2025{\natexlab{b}})Park, Howell, and
  Stewart]{park2025frequency}
Beom Park, Kathleen~C. Howell, and Shaun Stewart.
\newblock A frequency-domain differential corrector for quasi-periodic
  trajectory design and analysis.
\newblock \emph{Physica D: Nonlinear Phenomena}, page 134882,
  2025{\natexlab{b}}.
\newblock \doi{10.1016/j.physd.2025.134882}.

\bibitem[Park and Howell(2024)]{park2024assessment}
Beom Park and Kathleen~C. Howell.
\newblock {Assessment of dynamical models for transitioning from the Circular
  Restricted Three-Body Problem to an ephemeris model with applications}.
\newblock \emph{Celestial Mechanics and Dynamical Astronomy}, 136\penalty0 (6),
  2024.
\newblock \doi{10.1007/s10569-023-10178-9}.

\bibitem[Park et~al.(2021)Park, Folkner, Williams, and Boggs]{park2021jpl}
Ryan~S. Park, William~M. Folkner, James~G. Williams, and Dale~H. Boggs.
\newblock The {JPL} planetary and lunar ephemerides {DE440} and {DE441}.
\newblock \emph{The Astronomical Journal}, 161\penalty0 (3):\penalty0 105,
  2021.
\newblock \doi{10.3847/1538-3881/abd414}.

\bibitem[Vallado(2001)]{vallado2001fundamentals}
David~A. Vallado.
\newblock \emph{Fundamentals of Astrodynamics and Applications}.
\newblock Microcosm Press, 2 edition, 2001.

\bibitem[Longuski et~al.(2022)Longuski, Hoots, and
  Pollock]{longuski2022introduction}
James~M. Longuski, Felix~R. Hoots, and George~E. Pollock.
\newblock \emph{Introduction to Orbital Perturbations}.
\newblock Springer, 2022.
\newblock \doi{10.1007/978-3-030-89758-1}.

\bibitem[Broucke(2003)]{broucke2003long}
Roger~A. Broucke.
\newblock Long-term third-body effects via double averaging.
\newblock \emph{Journal of Guidance, Control, and Dynamics}, 26\penalty0
  (1):\penalty0 27--32, 2003.
\newblock \doi{10.2514/2.5041}.

\bibitem[Nie and Gurfil(2018)]{nie2018lunar}
Tao Nie and Pini Gurfil.
\newblock Lunar frozen orbits revisited.
\newblock \emph{Celestial Mechanics and Dynamical Astronomy}, 130\penalty0
  (10):\penalty0 61, 2018.
\newblock \doi{10.1007/s10569-018-9858-0}.

\bibitem[Russell and Brinckerhoff(2009)]{russell2009circulating}
Ryan~P. Russell and Adam~T. Brinckerhoff.
\newblock Circulating eccentric orbits around planetary moons.
\newblock \emph{Journal of Guidance, Control, and Dynamics}, 32\penalty0
  (2):\penalty0 423--435, 2009.
\newblock \doi{10.2514/1.38593}.

\bibitem[G{\'o}mez et~al.(2010)G{\'o}mez, Mondelo, and
  Sim{\'o}]{gomez2010collocation}
Gerard G{\'o}mez, Josep-Maria Mondelo, and Carles Sim{\'o}.
\newblock A collocation method for the numerical fourier analysis of
  quasi-periodic functions. i. numerical tests and examples.
\newblock \emph{Discrete and Continuous Dynamical Systems - B}, 14\penalty0
  (1):\penalty0 41--74, 2010.
\newblock \doi{10.3934/dcdsb.2010.14.41}.

\bibitem[Jorba and Villanueva(1997)]{jorba1997persistence}
Angel Jorba and Jordi Villanueva.
\newblock On the persistence of lower dimensional invariant tori under
  quasi-periodic perturbations.
\newblock \emph{Journal of Nonlinear Science}, 7:\penalty0 427--473, 1997.
\newblock \doi{10.1007/s003329900036}.

\bibitem[Sanaga et~al.(2026)Sanaga, Park, and Howell]{sanaga2026systematic}
Rohith~Reddy Sanaga, Beom Park, and Kathleen~C Howell.
\newblock Systematic numerical transitions of cislunar trajectories across
  dynamical models.
\newblock \emph{The Journal of the Astronautical Sciences}, 73\penalty0
  (1):\penalty0 10, 2026.
\newblock \doi{10.1007/s40295-026-00571-5}.

\bibitem[G{\'o}mez et~al.(2002)G{\'o}mez, Masdemont, and
  Mondelo]{gomez2002solar}
Gerard G{\'o}mez, Josep~J Masdemont, and Josep-Maria Mondelo.
\newblock Solar system models with a selected set of frequencies.
\newblock \emph{Astronomy \& Astrophysics}, 390\penalty0 (2):\penalty0
  733--749, 2002.
\newblock \doi{10.1051/0004-6361:20020625}.

\bibitem[Peale(1969)]{peale1969generalized}
Stanton~J. Peale.
\newblock Generalized {C}assini's laws.
\newblock \emph{The Astronomical Journal}, 74:\penalty0 483--489, 1969.
\newblock \doi{10.1086/110825}.

\bibitem[Rambaux and Williams(2011)]{rambaux2011moon}
Nicolas Rambaux and James~G. Williams.
\newblock The {M}oon's physical librations and determination of their free
  modes.
\newblock \emph{Celestial Mechanics and Dynamical Astronomy}, 109\penalty0
  (2):\penalty0 85--100, 2011.
\newblock \doi{10.1007/s10569-010-9314-2}.

\bibitem[Laskar(1999)]{laskar1999introduction}
Jacques Laskar.
\newblock Introduction to frequency map analysis.
\newblock In \emph{Hamiltonian systems with three or more degrees of freedom},
  pages 134--150. Springer, 1999.
\newblock \doi{10.1007/978-94-011-4673-9_13}.

\bibitem[Lian et~al.(2013)Lian, G{\'o}mez, Masdemont, and Tang]{lian2013note}
Yijun Lian, Gerard G{\'o}mez, Josep~J Masdemont, and Guojian Tang.
\newblock A note on the dynamics around the lagrange collinear points of the
  earth--moon system in a complete solar system model.
\newblock \emph{Celestial Mechanics and Dynamical Astronomy}, 115\penalty0
  (2):\penalty0 185--211, 2013.
\newblock \doi{10.1007/s10569-012-9459-2}.

\bibitem[G{\'o}mez and Mondelo(2001)]{gomez2001dynamics2}
G.~G{\'o}mez and J.M. Mondelo.
\newblock The dynamics around the collinear equilibrium points of the rtbp.
\newblock \emph{Physica D: Nonlinear Phenomena}, 157\penalty0 (4):\penalty0
  283--321, 2001.
\newblock ISSN 0167-2789.
\newblock \doi{10.1016/S0167-2789(01)00312-8}.

\bibitem[Olikara and Scheeres(2012)]{olikara2012numerical}
Zubin~P Olikara and Daniel~J Scheeres.
\newblock Numerical method for computing quasi-periodic orbits and their
  stability in the restricted three-body problem.
\newblock \emph{Advances in the Astronautical Sciences}, 145:\penalty0
  911--930, 2012.

\bibitem[Park et~al.(2025{\natexlab{c}})Park, Howell, and
  Stewart]{park2025elfo}
Beom Park, Kathleen~C. Howell, and Shaun Stewart.
\newblock Elliptical lunar frozen orbit constellation design within a model of
  evolving fidelity.
\newblock In \emph{35th AAS/AIAA Space Flight Mechanics Meeting, Kaua'i, HI,
  January 19--23}, 2025{\natexlab{c}}.

\end{thebibliography}

\clearpage
\appendix

\section{Long-Period Motion}
\label{app:lp-quadrature}

\subsection{Libration Period (Long-Period, $T_L$) via Quadrature}

\citet{russell2009circulating} express the libration period $T_L$ as a quadrature in $e$. Their formulation is first reviewed before a singularity-free adaptation is presented. For a given pair of integrals of motion $(C_1, C_2)$ introduced in Eqs.~\eqref{eq:C1_def}-\eqref{eq:C2_def}, the eccentricity oscillates between the turning points $e_{\min}$ and $e_{\max}$ (Fig.~\ref{fig:e_w_plane}), and both turning points admit closed-form expressions. The discriminant is defined as,
\begin{align}
\label{eq:Q}
Q &= 25(C_1 + C_2)^2 + 30(C_2 - C_1) + 9.
\end{align}
Then, for libration centered at $\omega = \pi/2$,
\begin{align}
\label{eq:emin_emax}
e_{\min,\,\max} &= \tfrac{1}{6}\sqrt{\mp 6\sqrt{Q} - 30(C_1 + C_2) + 18}.
\end{align}
At the frozen-orbit equilibrium $(e_{\text{eq}}, i_{\text{eq}}, \omega_{\text{eq}} = \pi/2)$, the two turning points coincide and $A_L = 0$. Squaring Eq.~\eqref{eq:dedt}, substituting $\sin^2 i$, $\sin^2\omega$, and $\cos^2\omega$ in terms of $(e, C_1, C_2)$ via Eqs.~\eqref{eq:C1_def}-\eqref{eq:C2_def}, and factoring the resulting polynomial in $e^2$ via Eq.~\eqref{eq:emin_emax} yields the quadrature \cite{russell2009circulating},
\begin{align}
\label{eq:TL-e}
T_L &= 2\int_{e_{\min}}^{e_{\max}} \frac{de}{|\dot{e}|}
    = \frac{8}{3\sqrt{3}}\,\frac{n}{(1-\mu)\, n_E}\,
      \int_{e_{\min}}^{e_{\max}}
      \frac{e\,de}{\sqrt{(2e^2-5C_2)(e^2-e_{\min}^2)(e_{\max}^2-e^2)}},
\end{align}
This expression is the libration counterpart for the circulating-orbit quadrature in Eq.~(11) from Russell and Brinckerhoff~\cite{russell2009circulating}.\footnote{Their Eq.~(11) shares the same integrand structure but is stated for the circulating case, with eccentricity bounds from their Eqs.~(12)-(13); the libration form given here substitutes the turning points of Eq.~\eqref{eq:emin_emax} (equivalent to their Eqs.~(13)-(14)). The $(1-\mu)$ factor, absorbed into the $(1-\mu)\approx 1$ approximation in~\citet{russell2009circulating}, is retained in the expression per Broucke~\cite{broucke2003long}, while the reported $T_L$ values adopt the same approximation (Sec.~\ref{sec:prelim}).} The integrand is singular at both $e = e_{\min}$ and $e = e_{\max}$, rendering direct numerical evaluation ill-conditioned for a small $A_L$. A singularity-free form follows from the substitution
\begin{align}
\label{eq:esub}
e^2 &= e_{\min}^2 + (e_{\max}^2 - e_{\min}^2)\sin^2 u, \qquad u \in \left[0,\,\tfrac{\pi}{2}\right],
\end{align}
For this substitution, $e\,de = \tfrac{1}{2}(e_{\max}^2 - e_{\min}^2)\sin 2u\,du$. The radical reduces to $\sqrt{(e^2-e_{\min}^2)(e_{\max}^2-e^2)} = \tfrac{1}{2}(e_{\max}^2-e_{\min}^2)\sin 2u$. As such, the singular factors cancel exactly, yielding,
\begin{align}
\label{eq:TL-final}
T_L &= \frac{8}{3\sqrt{3}}\,\frac{n}{(1-\mu)\, n_E}\,
      \int_0^{\pi/2}
      \frac{du}{\sqrt{K_1 + K_2\sin^2 u}},
\end{align}
where
\begin{align}
\label{eq:K1-K2}
K_1 &= 2e_{\min}^2 - 5C_2, \qquad
K_2 = 2\!\left(e_{\max}^2 - e_{\min}^2\right).
\end{align}
The integrand in Eq.~\eqref{eq:TL-final} is now smooth and bounded even at the limiting case $A_L \to 0$, where $K_2 \to 0$ and the integral reduces to $\pi/(2\sqrt{K_1})$, recovering the linearized libration period. The non-singular form enables accurate evaluation of $T_L$, and thereby the associated fundamental frequency $\nu_L = 2\pi/T_L$, across a range of $A_L$ values including the near-frozen regime where the original quadrature is ill-conditioned.

\subsection{Initial Condition Construction from $(A_{L,k}, \theta_{L,k})$}
\label{app:lp-ics}

The design vector ($\bm{\xi}$) in Sec.~\ref{subsec:al-dadm} assigns each satellite $k$ a long-period amplitude $A_{L,k}$ and phase $\theta_{L,k}$ alongside the inter-satellite phasing $(\theta_{S,k}, \theta_{M,k})$. Mean elements consistent with these descriptors are required to initiate the \acrshort{dadm} propagation (Eqs.~\eqref{eq:dadt}-\eqref{eq:dthetadt}), and the construction proceeds in three steps.

Firstly, the constants of motion $(C_{1,k}, C_{2,k})$ are determined. The amplitude $A_{L,k}$ alone does not uniquely locate the libration curve, as a one-parameter family of $(C_1, C_2)$ yields the same eccentricity swing. Moreover, the long-period-averaged RAAN drift rate $\langle\dot\Omega\rangle$ depends on $(C_1, C_2)$, so permitting each satellite to adopt its own $\langle\dot\Omega\rangle_k$ renders the medium-period frequency $\nu_{M,k}$ non-uniform across the constellation and breaks the designed phasing over time. The pair $(C_{1,k}, C_{2,k})$ is therefore supplied by solving the following,
\begin{align}
\label{eq:ic-target}
A_L(C_{1,k}, C_{2,k}) &= A_{L,k}, &
\langle\dot\Omega\rangle(C_{1,k}, C_{2,k}) &= 1-\nu_M,
\end{align}
where $\nu_M$ is fixed from the prior analysis. The long-period-averaged drift rate is evaluated as $(\int_{0}^{T_L} \dot{\Omega} \; dt )/ T_L$ from Eqs.~\eqref{eq:dOmegadt} and~\eqref{eq:TL-final}.

Secondly, the long-period phase $\theta_{L,k}$ parameterizes position along the libration curve in the $(e, i, \omega)$ phase space (Fig.~\ref{fig:e_w_plane}), with the convention $\theta_{L,k} = 0$ at the $e = e_{\max}$ turning point ($\omega = \pi/2$, $\dot e = 0$). Given $(C_{1,k}, C_{2,k})$ from the first step, the mean elements $(e_k, i_k, \omega_k, \Omega_k)$ at $\theta_{L,k}$ are supplied by integrating the \acrshort{dadm} equations of motion from the $e_{\max}$ turning point over the interval $\Delta t_k = \theta_{L,k}\,T_L / (2\pi)$, with $T_L$ from Eq.~\eqref{eq:TL-final}. 

Lastly, $\Omega_k$ combines both the secular drift and the long-period libration. Isolating the libration contribution yields
\begin{align}
\label{eq:domega-LP}
\delta\Omega_{LP,k} = \Omega_k - \langle\dot\Omega\rangle_k\,\Delta t_k,
\end{align}
The libration contribution is then added to the designed $\theta_{M,k}$ such that the mean inter-satellite RAAN spacing matches the Sec.~\ref{subsec:al-dadm} design even when $\theta_{L,k} \neq 0$. The full mean element vector for satellite $k$ is
\begin{align}
\label{eq:lp-ics}
[a,\; e_k,\; i_k,\; -\theta_{M,k} + \delta\Omega_{LP,k},\; \omega_k,\; \theta_{S,k}],
\end{align}
where $a$ is the common semi-major axis and the short-period phase $\theta_{S,k}$ is assigned to the mean anomaly. Note, $\theta_{M,k} = -\Omega_{k}$ at the initial epoch ($t_0 = 0$ from Eq.~\eqref{eq:eof_mrf}) without losing generality. Equation~\eqref{eq:lp-ics} is employed as the initial condition for the \acrshort{dadm} dynamics (Eqs.~\eqref{eq:dadt}-\eqref{eq:dthetadt}); numerical propagation then supplies spacecraft state within the constellation. 

\section{User-Geometry Assumptions}
\label{app:obliquity}

The analysis in the main body assumes a constant obliquity at $\epsilon_M = 6.68^\circ$ and a user at the \acrshort{lsp}. These assumptions are further examined in the current appendix.

\subsection{Zero Obliquity ($\epsilon_M = 0, \hat{\bm{s}}_M = \hat{\bm{z}}$)}
\label{app:onaxis}

For further analytical tractability and lower design dimensions, it is possible to assume a zero obliquity of the lunar spin axis, i.e., $\epsilon_M = 0, \hat{\bm{s}}_M = \hat{\bm{z}}$ (Fig.~\ref{fig:user_geom_schematic}). Under such an assumption, the effective angular dimensionality for the search space decreases by one. The user is assumed to be located at the \acrshort{lsp}, consistent with the main body.

\subsubsection{\acrshort{dadm}}
\label{subsec:dadm-results}
With the longitudinal invariance of the \acrshort{dadm} (Sec.~\ref{subsec:design-vars}), the survey collapses to a one-dimensional angular grid in $\theta_S$ alone (50 points). The resulting heatmap (Fig.~\ref{fig:sweep_heatmap}) preserves the inclination band near $i \approx 51^\circ$ identified in the obliquity-corrected case (Fig.~\ref{fig:sweep_obl_heatmap}), although with higher coverage fraction (brighter heatmap). The reference solution at $(a, i) = (14{,}200~\text{km}, 50.5^\circ)$ is examined in Fig.~\ref{fig:sweep_torus} that demonstrates vertical stripes independent of $\theta_M$, signifying the extra symmetry from $\epsilon_M=0$.

\subsubsection{\acrshort{cr3bp}}
\label{sec:cr3bp-torus}
\begin{figure}[h!]
    \centering
    \begin{subfigure}[b]{0.48\textwidth}
        \centering
        \includegraphics[width=\textwidth]{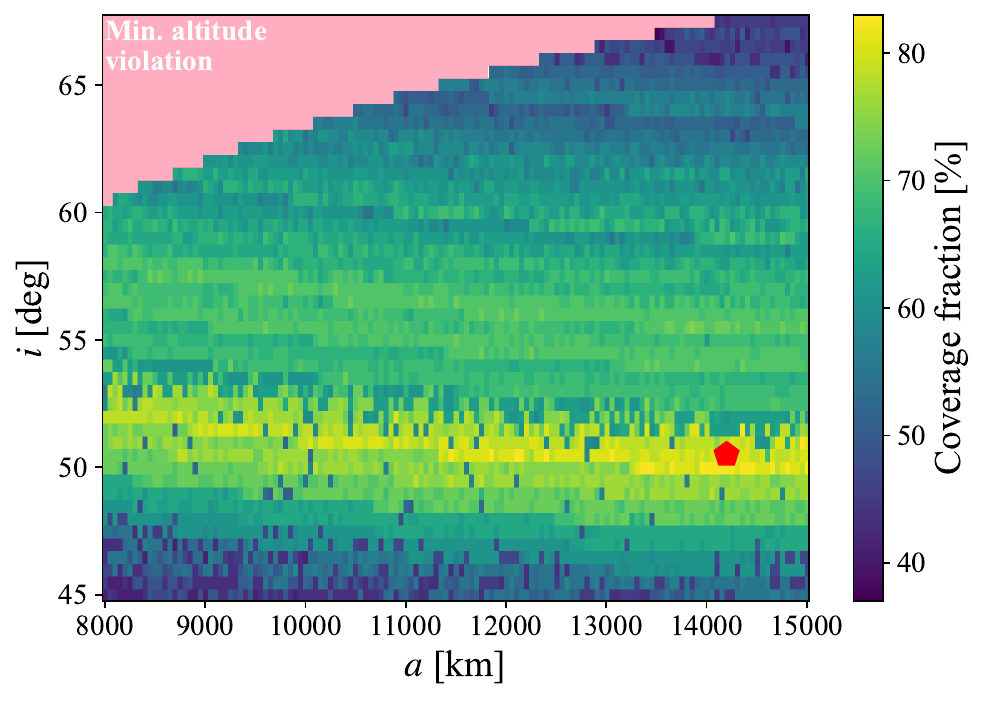}
        \caption{\acrshort{dadm} coverage heatmap.}
        \label{fig:sweep_heatmap}
    \end{subfigure}
    \hfill
    \begin{subfigure}[b]{0.48\textwidth}
        \centering
        \includegraphics[width=\textwidth]{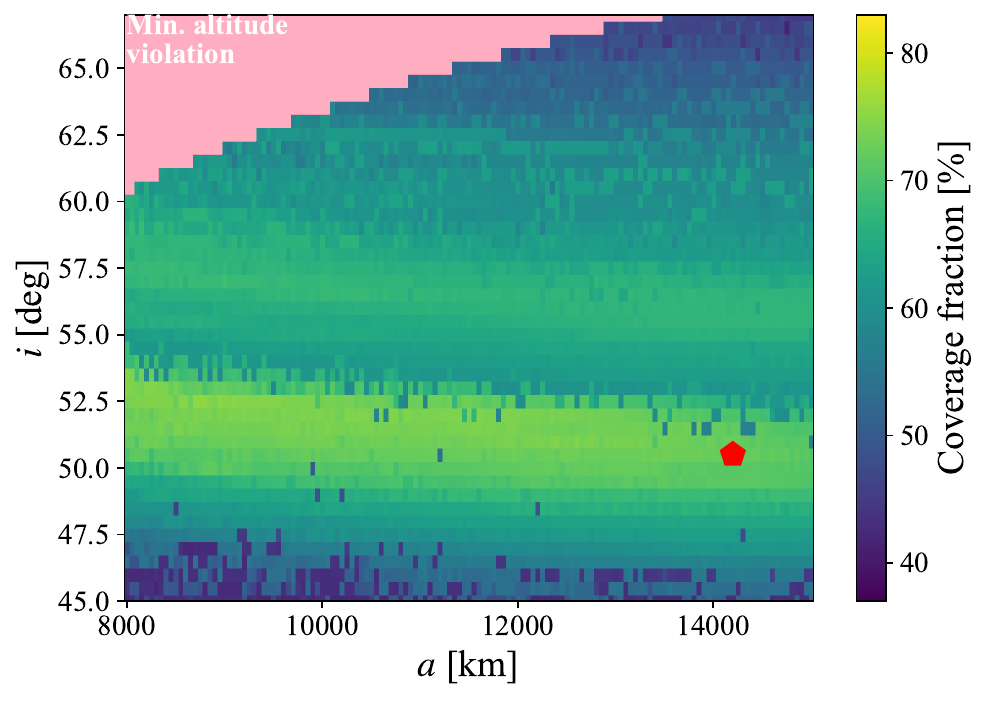}
        \caption{\acrshort{cr3bp} coverage heatmap.}
        \label{fig:cr3bp_sweep_heatmap}
    \end{subfigure}
    \\[6pt]
    \begin{subfigure}[b]{0.48\textwidth}
        \centering
        \includegraphics[width=\textwidth]{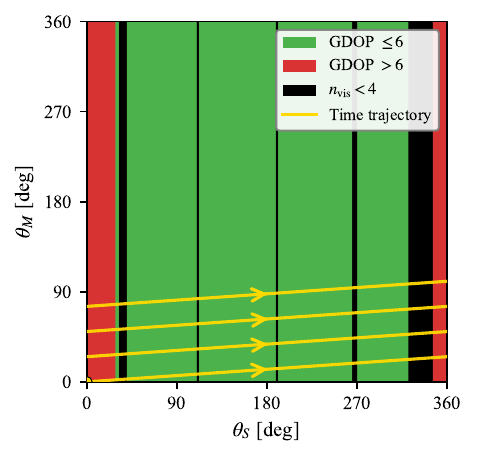}
        \caption{\acrshort{dadm} \acrshort{gdop} torus map (coverage $76.2\%$).}
        \label{fig:sweep_torus}
    \end{subfigure}
    \hfill
    \begin{subfigure}[b]{0.48\textwidth}
        \centering
        \includegraphics[width=\textwidth]{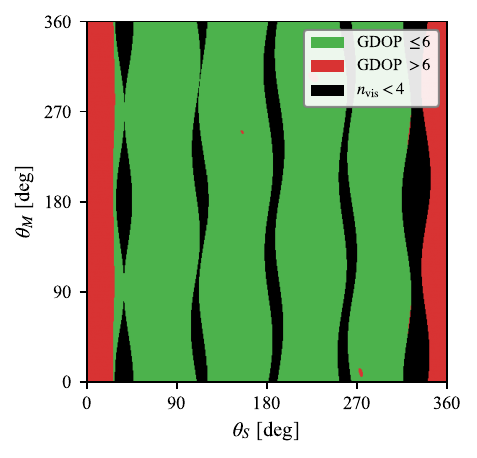}
        \caption{\acrshort{cr3bp}-optimal \acrshort{gdop} torus map (coverage $70.7\%$).}
        \label{fig:cr3bp_torus_cr3bp}
    \end{subfigure}
    \caption{Constellation design survey with $\epsilon_M = 0^\circ$ at the \acrshort{lsp}. \acrshort{dadm} (a),\,(c); \acrshort{cr3bp} (b),\,(d). The reference configuration $(a, i) = (14{,}200~\text{km}, 50.5^\circ)$ is marked by the red pentagon (\redpentagon) in (a),\,(b).}
    \label{fig:b1_sweep}
\end{figure}

\begin{figure}[h!]
    \centering
    \begin{subfigure}[b]{0.48\textwidth}
        \centering
        \includegraphics[width=\linewidth]{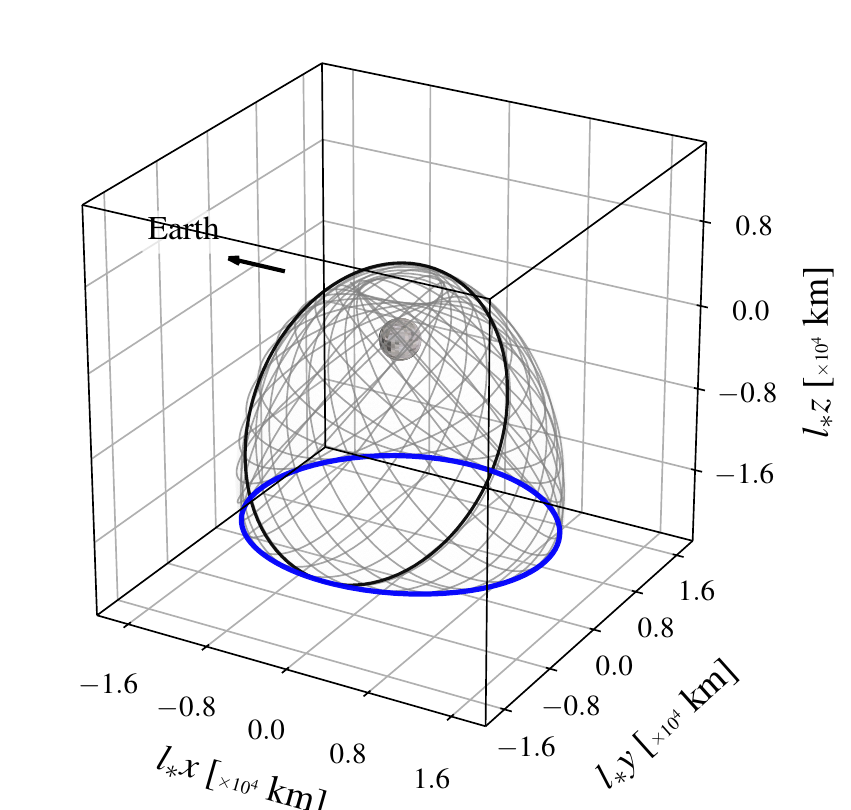}
        \caption{Geometry within the \acrshort{mrf} (25 invariant curves at evenly spaced $\theta_M$).}
        \label{fig:cr3bp_qpo_torus}
    \end{subfigure}
    \hfill
    \begin{subfigure}[b]{0.48\textwidth}
        \centering
        \includegraphics[width=\linewidth]{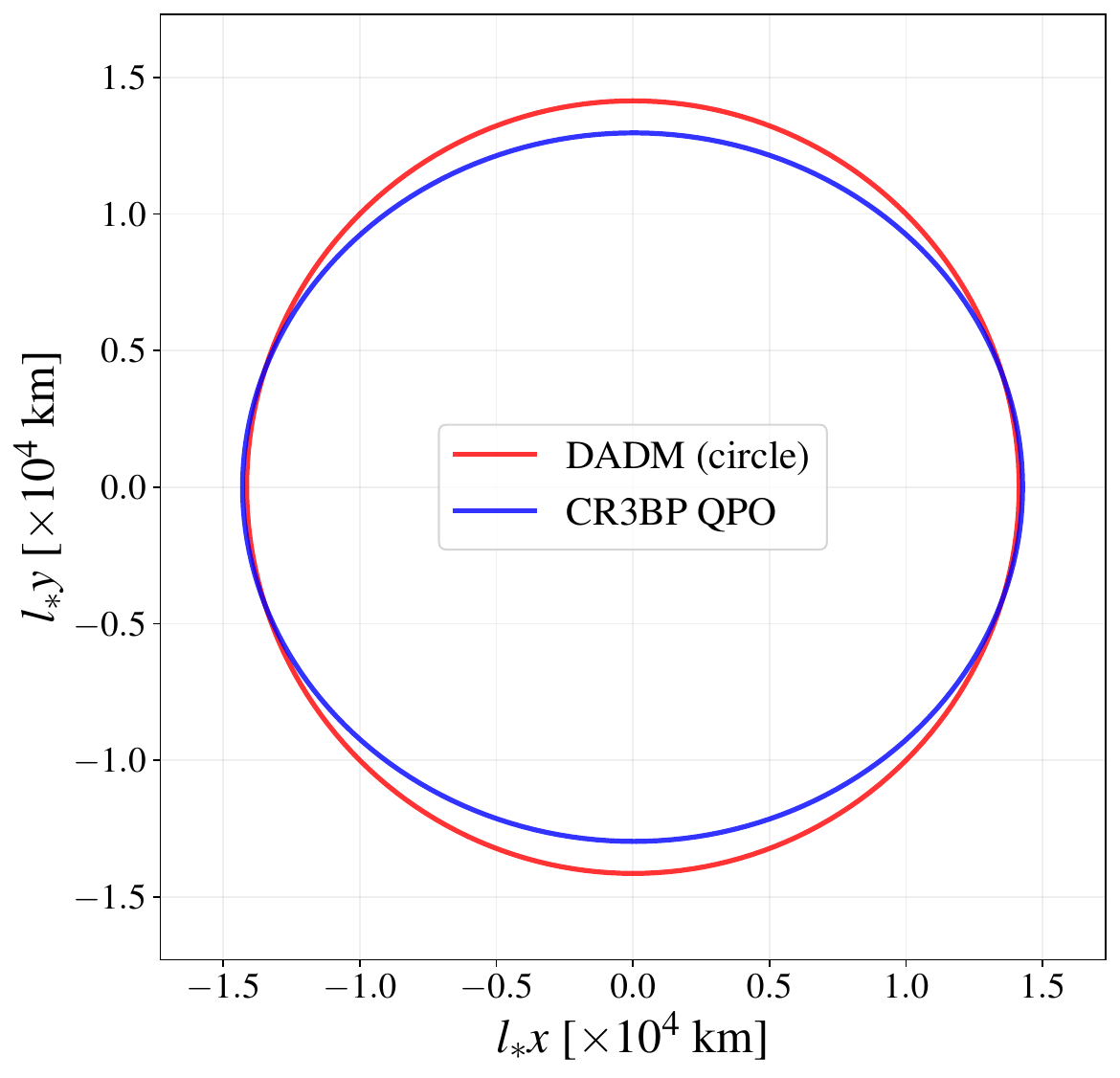}
        \caption{Apolune stroboscopic section.}
        \label{fig:cr3bp_qpo_strobo}
    \end{subfigure}
    \caption{\acrshort{cr3bp} 2D torus at the reference $(a, i) = (14{,}200~\text{km}, 50.5^\circ)$, $(\nu_S, \nu_M) \approx (15.55, 1.086)$~rad/nd.}
    \label{fig:cr3bp_qpo}
\end{figure}

With the reduced dimension from the zero obliquity, the search space within the \acrshort{cr3bp} becomes more tractable. Generally, the frequency pair $(\nu_S, \nu_M)$ specifies a structure along families of \acrshort{elfos} as 2D tori, i.e., \acrfull{qpos}, constructed via the stroboscopic mapping method supplied by G{\'o}mez and Mondelo~\cite{gomez2001dynamics2} as well as Olikara and Scheeres~\cite{olikara2012numerical}, with the application to \acrshort{elfos} following the authors' prior work~\cite{park2025elfo}. Figure~\ref{fig:cr3bp_qpo} illustrates a sample 2D torus at the reference $(\nu_S, \nu_M) \approx (15.55, 1.086)$~rad/nd. Geometry within the \acrshort{mrf} is illustrated in Fig.~\ref{fig:cr3bp_qpo_torus}. The gray closed curves are 25 invariant curves at evenly spaced $\theta_M$. The black and blue plots correspond to invariant curves constructed at $\theta_M = 0^\circ$ and $\theta_S = 180^\circ$, respectively. The stroboscopic section at $\theta_S = 180^\circ$ is examined further in Fig.~\ref{fig:cr3bp_qpo_strobo}. The 2D parametrization enforces $A_L = 0$ exactly, in contrast to the time-domain propagation of Fig.~\ref{fig:strobo_cr3bp} where the residual long-period libration broadens the apolune section into a finite-width annulus; the 2D torus therefore supplies the frozen configuration by construction where the phasing $(\theta_S, \theta_M)$ is optimized, decoupled from the long-period motion. For each grid point $(a,i)$, the corresponding 2D torus with the same $(\nu_S, \nu_M)$ (supplied from Eqs.~\eqref{eq:nu_s}-\eqref{eq:nu_m} in the \acrshort{dadm}) is constructed within the \acrshort{cr3bp}. In contrast to the \acrshort{dadm}, the \acrshort{cr3bp} breaks the longitudinal invariance, so the survey requires the full two-dimensional $(\theta_S, \theta_M)$ torus evaluation. The resulting heatmap (Fig.~\ref{fig:cr3bp_sweep_heatmap}) preserves the inclination band from the \acrshort{dadm} survey and additionally develops a secondary band near $i \approx 56^\circ$ (originally located within the \acrshort{dadm} with the obliquity, Fig.~\ref{fig:sweep_obl}). The coverage is examined at the reference configuration $(a, i) = (14{,}200~\text{km}, 50.5^\circ)$ (Fig.~\ref{fig:cr3bp_torus_cr3bp}). Within the \acrshort{cr3bp}, the longitudinal invariance is lost, and a repeating pattern dominated by the $2\nu_M$ harmonic emerges, consistent with the spectral analysis of Sec.~\ref{subsec:freq-across-models}. The optimized coverage of $70.7\%$ remains $5.5$~pp below the \acrshort{dadm} benchmark of $76.2\%$ (Fig.~\ref{fig:sweep_torus}), with the residual gap manifesting as $n_{\mathrm{vis}}<4$ (black) patterns with larger widths. The $5.5$~pp gap therefore characterizes the unremovable $2\nu_M$ content present in the \acrshort{cr3bp} but absent in the doubly-averaged dynamics for this reference configuration.

The zero-obliquity assumption serves two purposes in the preliminary analysis. Firstly, it maximizes the available symmetries; the \acrshort{dadm} survey collapses to a one-dimensional grid in $\theta_S$, enabling even more rapid global exploration of the $(a, i)$ space with some qualitative resemblance to the results from the obliquity-incorporated baseline (Sec.~\ref{sec:design-exploration}). Secondly, the side-by-side \acrshort{dadm}-\acrshort{cr3bp} comparison in Fig.~\ref{fig:b1_sweep} offers an intuitive understanding of the $2\nu_M$ modulation introduced by the unaveraged Earth forcing (Sec.~\ref{subsec:freq-across-models}).

\subsection{Non-\acrshort{lsp} User Locations with Obliquity ($\epsilon_M = 6.68^\circ$)}
\label{app:non-lsp}

Relaxing the fixed user location (at the \acrshort{lsp}) assumption introduces draconic-period ($\omega_3$) oscillations in both the user latitude and longitude as apparent within the \acrshort{mrf}, enlarging the \acrshort{gdop} dimensionality from $2$ $(\theta_S, \theta_M)$ to $3$ $(\theta_S, \theta_M, \alpha_3)$ within the angular domain. A complete optimization over this 3D quasi-periodic structure remains outside the scope of this paper; the present analysis rather examines the inherent complexity and identifies tractable approaches for navigating the enlarged design space, focusing on the \acrshort{dadm}.

First recall that in the case of \acrshort{lsp}, the apparent location within the \acrshort{mrf} traces a constant latitude $\varphi = -90 + \epsilon_M = -83.32^\circ$ and a longitude that drifts linearly at the draconic frequency, $\lambda(t) = -\omega_3 t$, where $\epsilon_M = 6.68^\circ$ is the lunar obliquity (Sec.~\ref{subsec:obliquity}, Fig.~\ref{fig:user_geom}). For a non-\acrshort{lsp} user, consider $(\varphi_\mathrm{PA}, \lambda_\mathrm{PA})$, where the subscript $\mathrm{PA}$ denotes the \texttt{MOON\_PA} frame, $\bm{r}_\mathrm{PA} = \rho_M(\cos\varphi_\mathrm{PA}\cos\lambda_\mathrm{PA},\, \cos\varphi_\mathrm{PA}\sin\lambda_\mathrm{PA},\, \sin\varphi_\mathrm{PA})^T$, where $\rho_M = 1{,}737.106$~km is the lunar radius. Then, the apparent user position within the \acrshort{mrf} evolves according to ${\bm{r}}(t) = \bm{C}_{\hat{\bm{z}}}(-\omega_3 t)\,\bm{C}_{\hat{\bm{y}}}(\epsilon_M)\,\bm{C}_{\hat{\bm{z}}}(\omega_3 t)\,{\bm{r}}_\mathrm{PA}$. Without loss of generality, $\lambda_\mathrm{PA} = 0$ is adopted in the expressions below, since the longitudinal invariance of the \acrshort{dadm} (Sec.~\ref{subsec:design-vars}) renders the coverage independent of the reference longitude; a nonzero $\lambda_\mathrm{PA}$ merely shifts the draconic phase. Algebra leads to the apparent \acrshort{mrf} latitude in the following relationship,
\begin{align}
\label{eq:lat-mrf}
\sin\varphi(t) = \sin\varphi_\mathrm{PA}\cos\epsilon_M - \cos\varphi_\mathrm{PA}\sin\epsilon_M\cos(\omega_3 t),
\end{align}
oscillating within the bound $[\varphi_\mathrm{PA} - \epsilon_M,\, \varphi_\mathrm{PA} + \epsilon_M]$ at the draconic period (extrema attained at $\omega_3 t = 0, \pi$). The longitude within the \acrshort{mrf} also follows as,
\begin{align}
\label{eq:lon-mrf}
\tan\big(\lambda(t) - \lambda_\mathrm{PA}\big) \;=\; \frac{\cos\varphi_\mathrm{PA}\,(1-\cos\epsilon_M)\,\cos(\omega_3 t)\,\sin(\omega_3 t) \;-\; \sin\varphi_\mathrm{PA}\,\sin\epsilon_M\,\sin(\omega_3 t)}{\cos\varphi_\mathrm{PA}\,\big[1 - (1-\cos\epsilon_M)\,\cos^2(\omega_3 t)\big] \;+\; \sin\varphi_\mathrm{PA}\,\sin\epsilon_M\,\cos(\omega_3 t)},
\end{align}
whose peak excursion grows with $|\tan\varphi_\mathrm{PA}|$. As $\varphi_\mathrm{PA} \to -90^\circ$, Eq.~\eqref{eq:lat-mrf} reduces to the constant $\sin\varphi = -\cos\epsilon_M$ and Eq.~\eqref{eq:lon-mrf} degenerates into $\tan(\lambda - \lambda_\mathrm{PA}) = -\tan(\omega_3 t)$, recovering the \acrshort{lsp} limit. Equations~\eqref{eq:lat-mrf}-\eqref{eq:lon-mrf} replace the trivial ``constant latitude, free longitude'' treatment of the \acrshort{lsp} user in Sec.~\ref{sec:design-survey} with a draconic-period oscillation in both angles as apparent within the \acrshort{mrf}.

\begin{figure}[h!]
    \centering
    \begin{subfigure}[b]{0.32\textwidth}
        \centering
        \includegraphics[width=\textwidth]{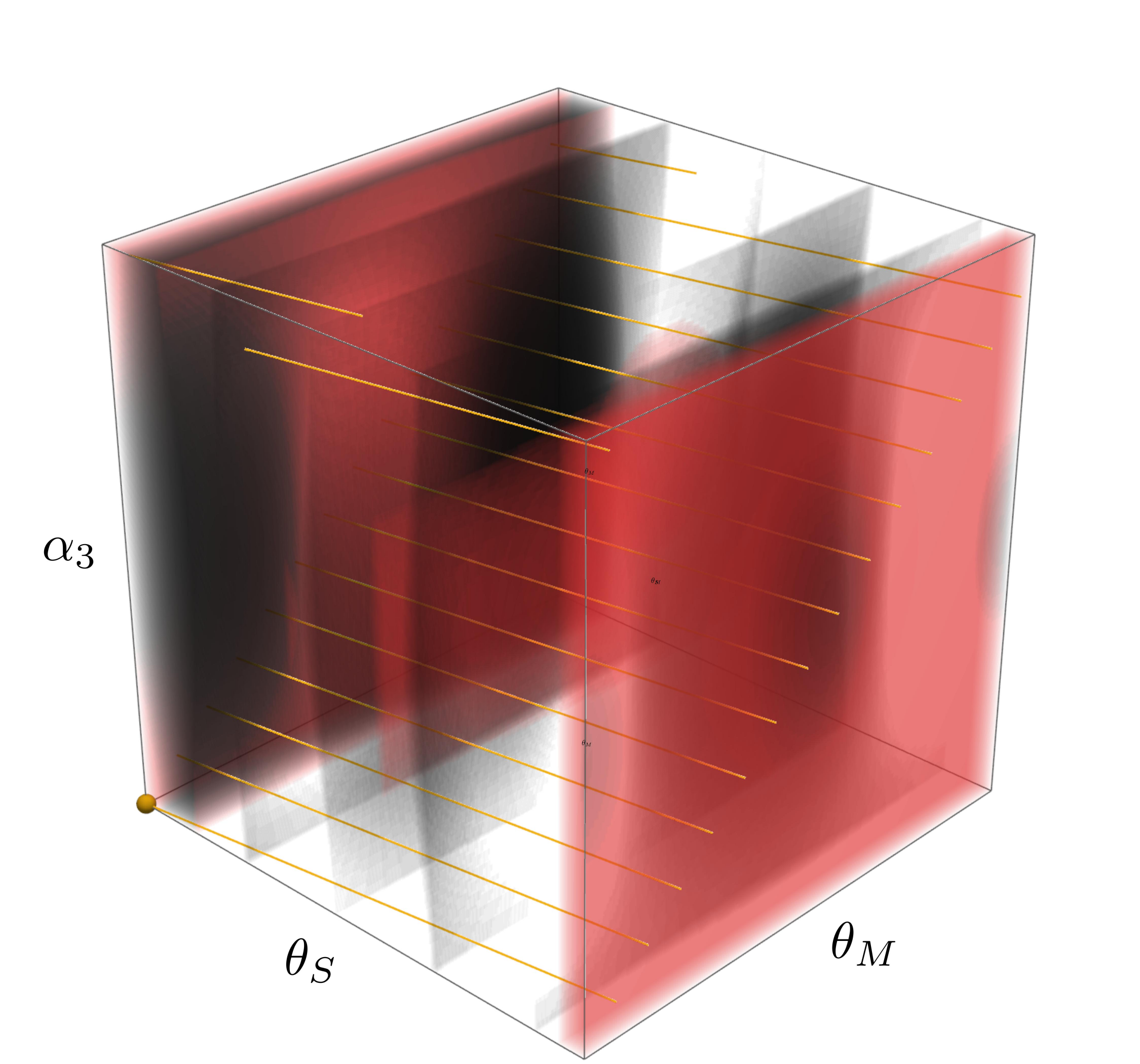}
        \caption{Combined volume: \acrshort{gdop} $> 6$ (red) and $n_{\mathrm{vis}} < 4$ (black).}
        \label{fig:non_lsp_3d_all}
    \end{subfigure}
    \hfill
    \begin{subfigure}[b]{0.32\textwidth}
        \centering
        \includegraphics[width=\textwidth]{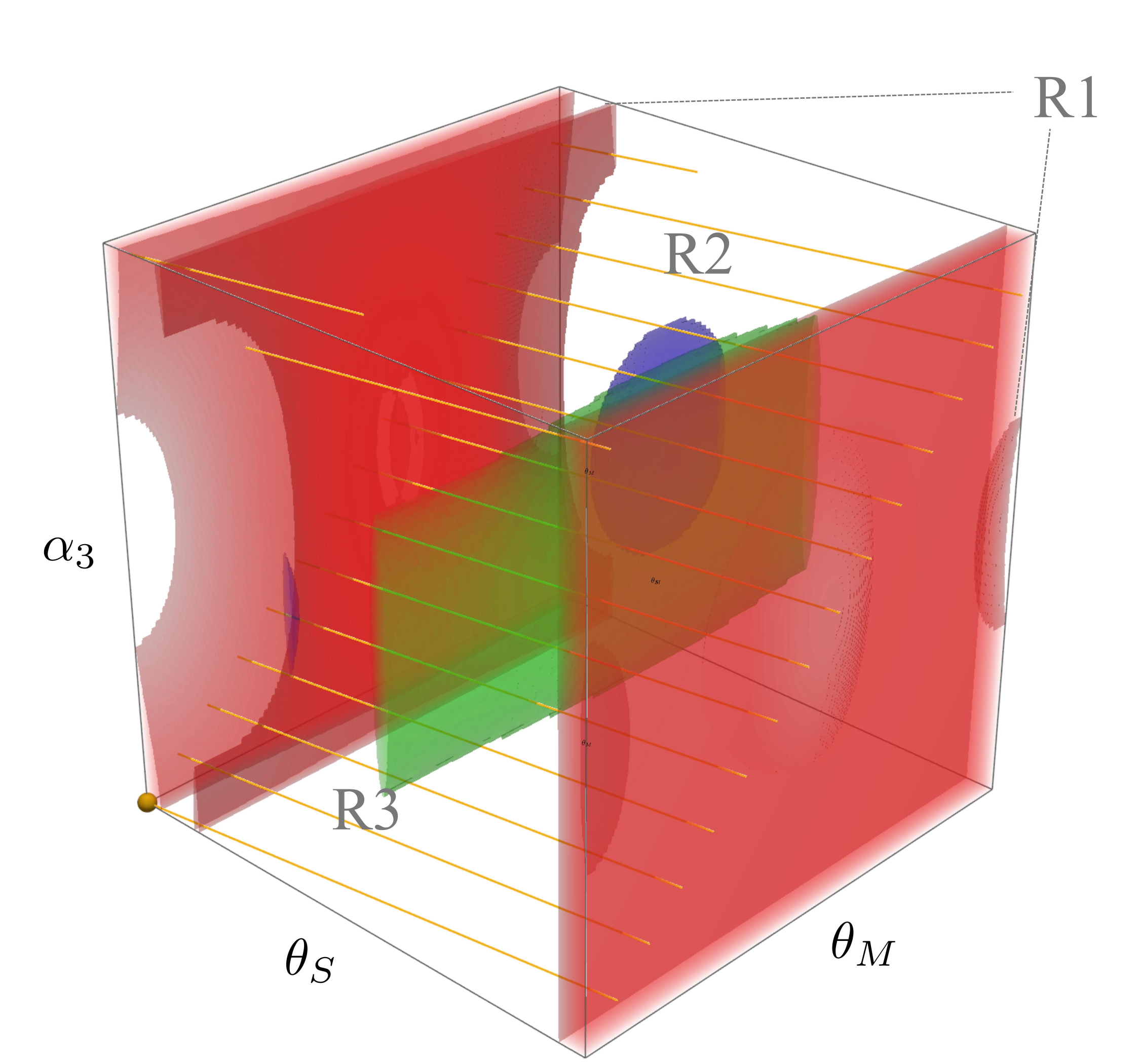}
        \caption{\acrshort{gdop} $> 6$ pillars decomposed into 3 connected components (R1-R3).}
        \label{fig:non_lsp_3d_red}
    \end{subfigure}
    \hfill
    \begin{subfigure}[b]{0.32\textwidth}
        \centering
        \includegraphics[width=\textwidth]{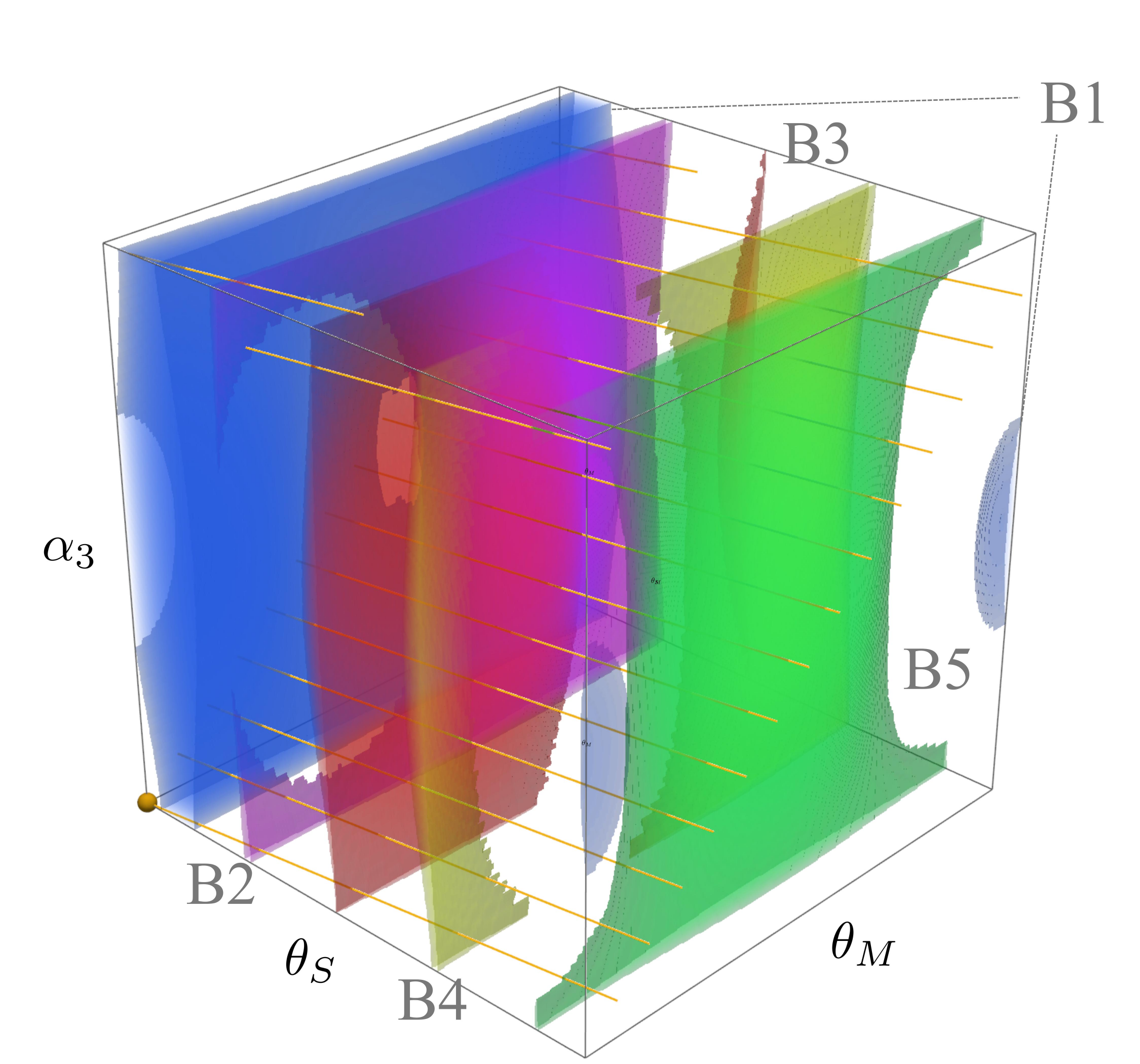}
        \caption{$n_{\mathrm{vis}} < 4$ pillars decomposed into 5 connected components (B1-B5).}
        \label{fig:non_lsp_3d_black}
    \end{subfigure}
    \caption{Three-dimensional angular space $(\theta_S, \theta_M, \alpha_3)$ for the non-\acrshort{lsp} user of Fig.~\ref{fig:non_lsp_snapshots}. (a)~Combined view. (b),(c)~Connected-component decomposition with periodic merging across $\theta_S$, $\theta_M$ boundaries. The gold curve in each panel traces the time evolution from the origin with slope $(\nu_S, \nu_M, \omega_3)$ over one draconic period.}
    \label{fig:non_lsp_3d}
\end{figure}

\begin{figure}[h!]
    \centering
    \begin{subfigure}[b]{0.48\textwidth}
        \centering
        \includegraphics[width=\textwidth]{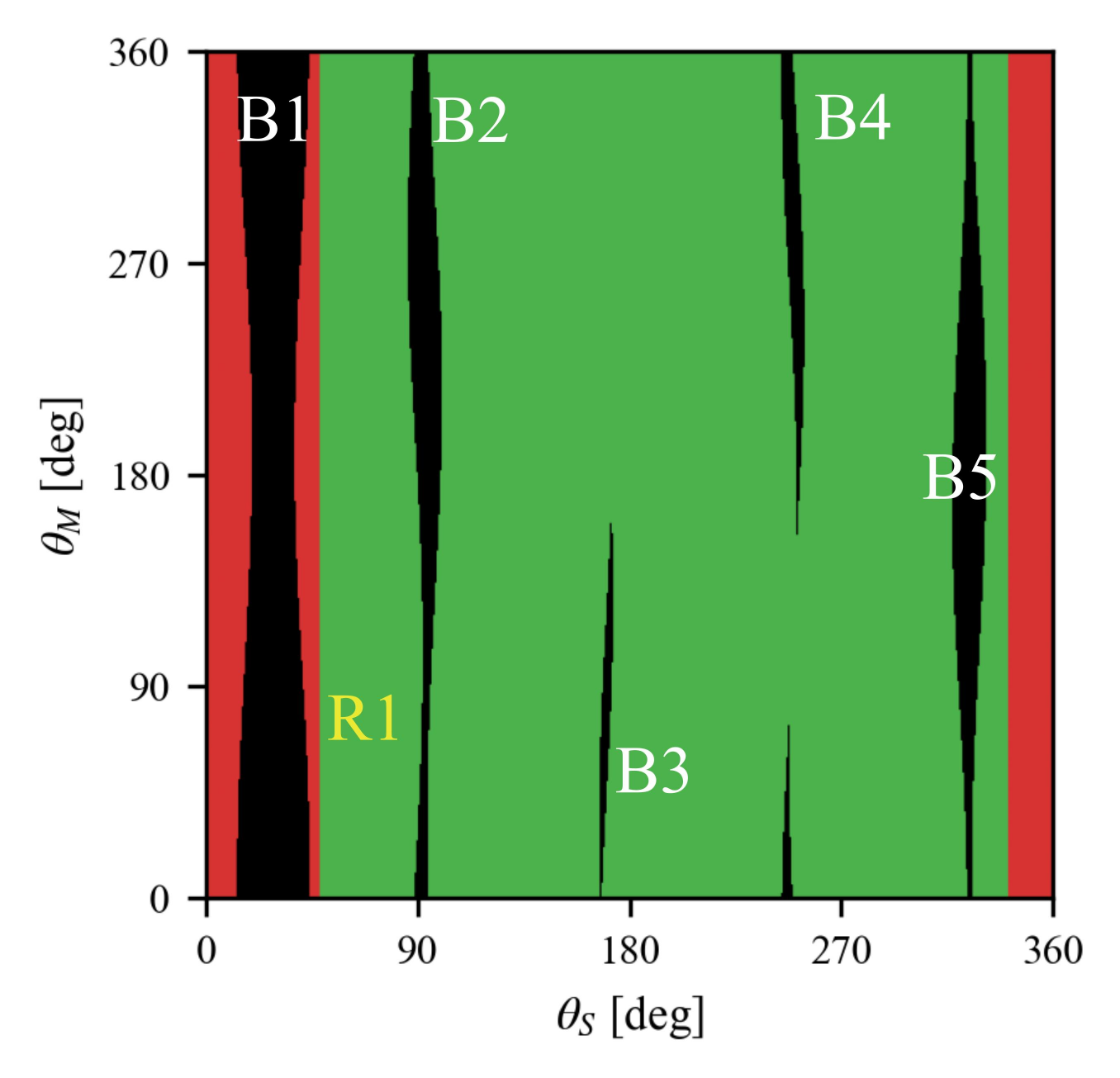}
        \caption{$\alpha_3 = 0^\circ$, $\varphi = -86.68^\circ$.}
        \label{fig:non_lsp_snap_a}
    \end{subfigure}
    \hfill
    \begin{subfigure}[b]{0.48\textwidth}
        \centering
        \includegraphics[width=\textwidth]{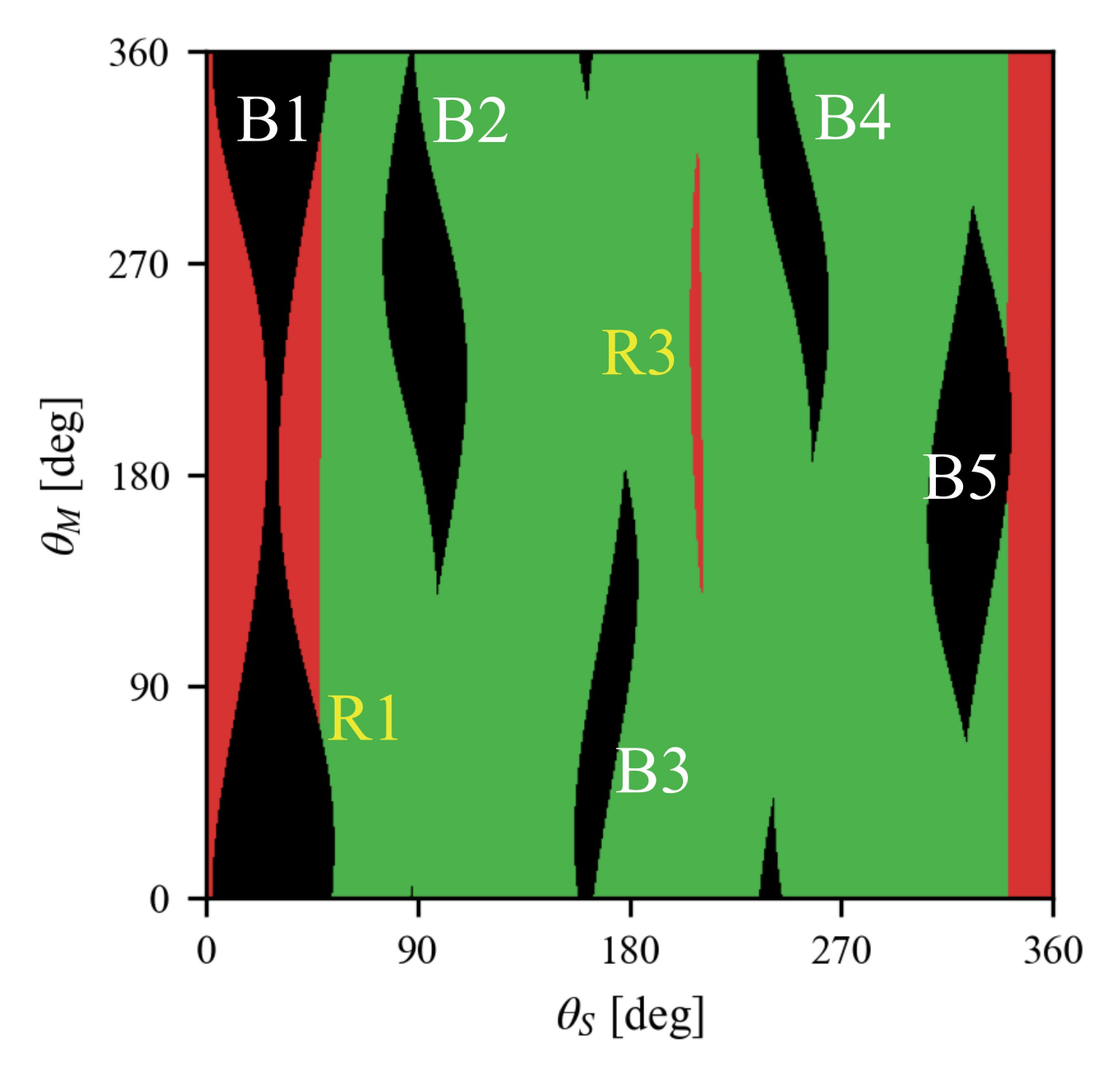}
        \caption{$\alpha_3 = 90^\circ$, $\varphi = -77.99^\circ$.}
        \label{fig:non_lsp_snap_b}
    \end{subfigure}
    \\[6pt]
    \begin{subfigure}[b]{0.48\textwidth}
        \centering
        \includegraphics[width=\textwidth]{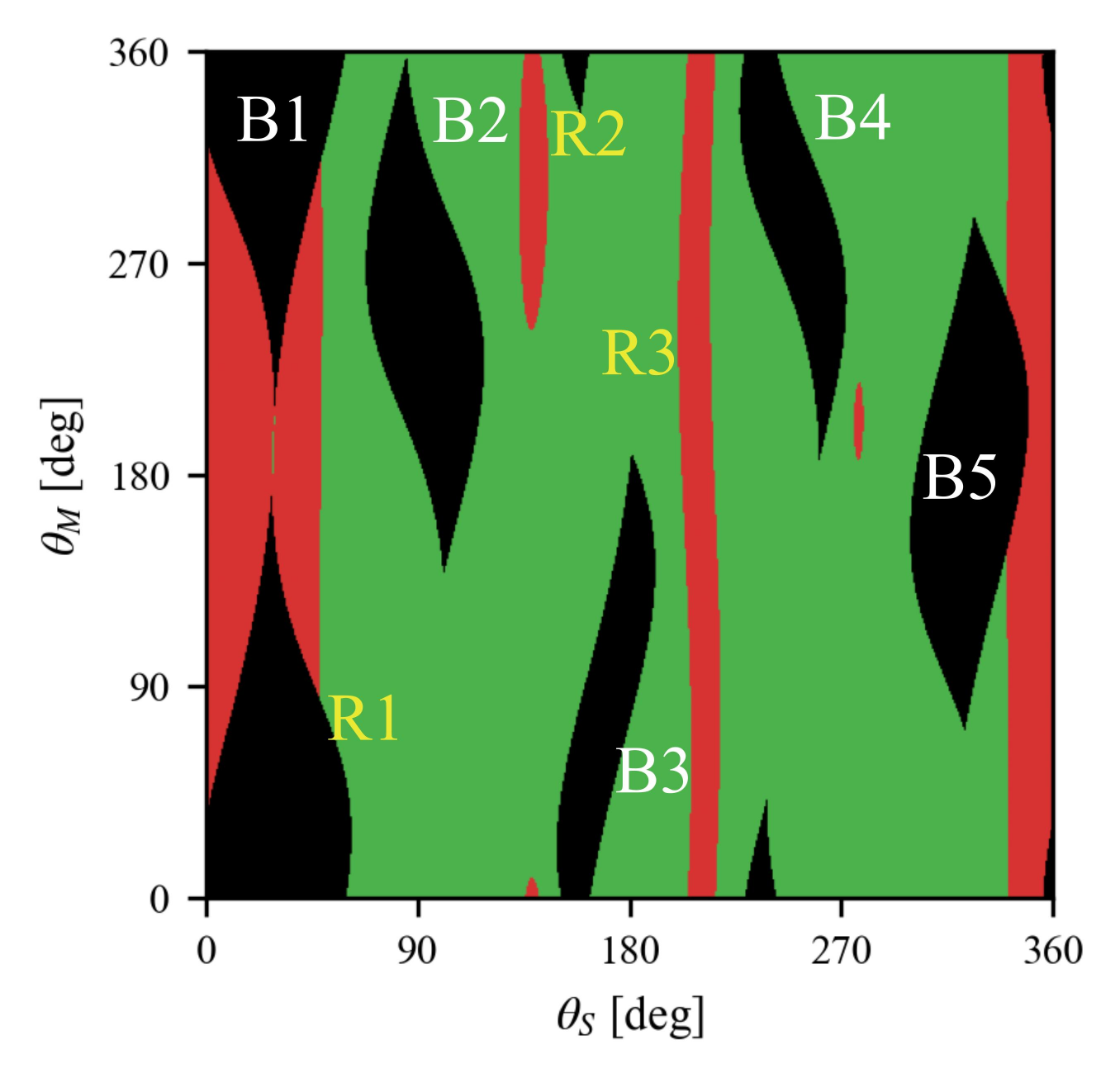}
        \caption{$\alpha_3 = 180^\circ$, $\varphi = -73.32^\circ$.}
        \label{fig:non_lsp_snap_c}
    \end{subfigure}
    \hfill
    \begin{subfigure}[b]{0.48\textwidth}
        \centering
        \includegraphics[width=\textwidth]{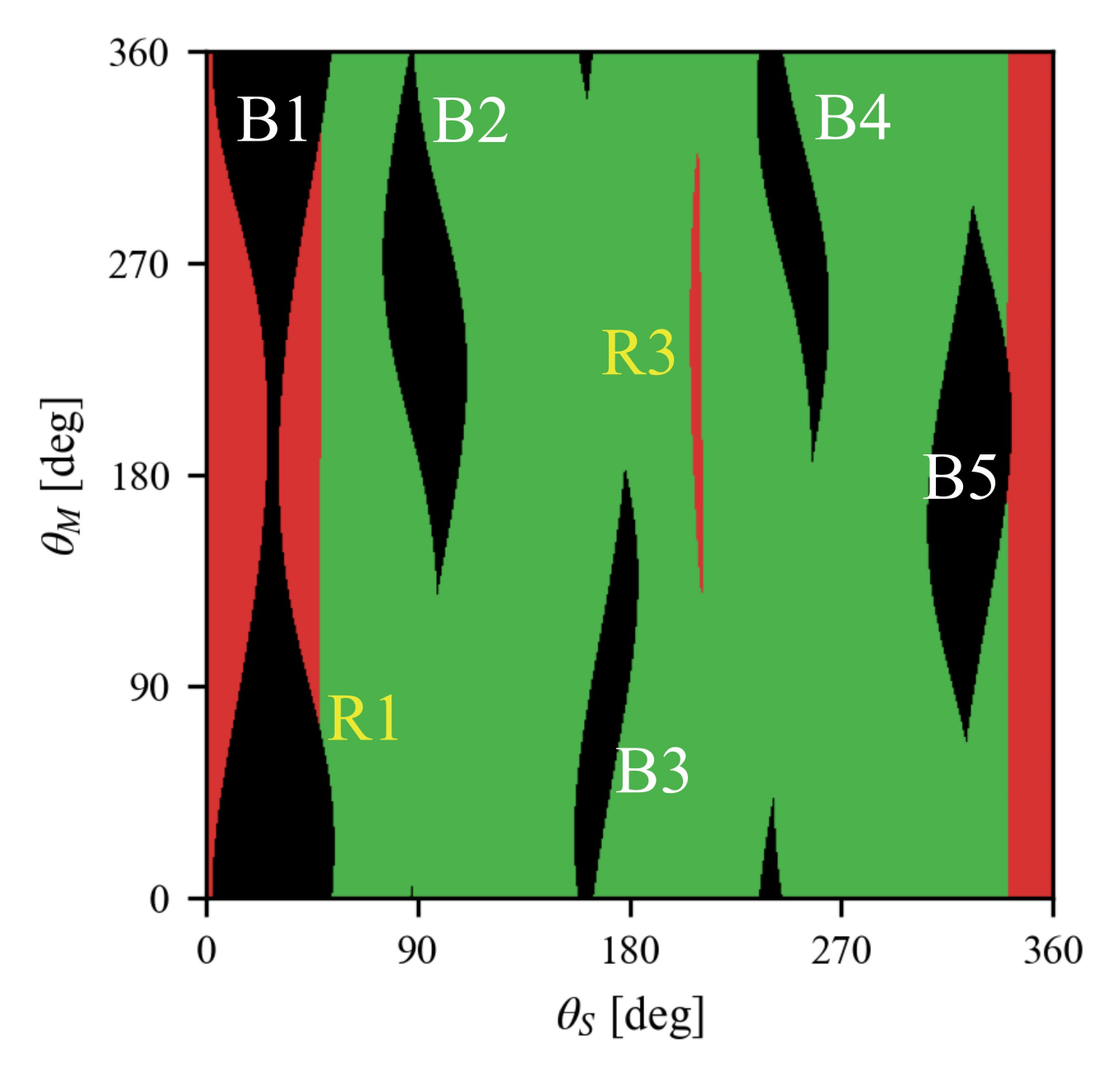}
        \caption{$\alpha_3 = 270^\circ$, $\varphi = -77.99^\circ$.}
        \label{fig:non_lsp_snap_d}
    \end{subfigure}
    \caption{\acrshort{gdop} torus snapshots at four draconic phases for a non-\acrshort{lsp} user at $\varphi_\mathrm{PA} = -80^\circ$, with the reference constellation of Fig.~\ref{fig:sweep_obl_phases}. Color convention follows Fig.~\ref{fig:sweep_obl_torus}.}
    \label{fig:non_lsp_snapshots}
\end{figure}

These $\omega_3$-induced oscillations within the \acrshort{mrf} result in a fundamental increase in the dimensionality. As such, for any non-\acrshort{lsp} locations within the \acrshort{dadm}, sampling in three-dimensional angular space is generally required to represent a constellation behavior over an infinite time horizon. Consider the optimal phasing from Fig.~\ref{fig:sweep_obl_phases}, now evaluated at the non-\acrshort{lsp} location $(\varphi_\mathrm{PA}, \lambda_{\mathrm{PA}}) = (-80^\circ, 0^\circ)$. The \acrshort{gdop} is a three-dimensional quasi-periodic function on the angular torus $(\theta_S, \theta_M, \alpha_3) \in [0, 2\pi)^3$, where $\alpha_3 = \omega_3 t$ denotes the added angular dimension parametrizing the user's draconic phase. Figure~\ref{fig:non_lsp_3d} renders this three-dimensional torus within a $500 \times 500 \times 100$ grid constructed on $(\theta_S, \theta_M, \alpha_3)$ with the time-evolution line overlaid. The time evolution traces a straight line in this three-dimensional space with slope $(\nu_S, \nu_M, \omega_3)$. Without loss of generality, the line starts at the origin at the initial time. The red and black volumes in Fig.~\ref{fig:non_lsp_3d_all} correspond to \acrshort{gdop} $> 6$ and $n_{\mathrm{vis}} < 4$, respectively, while the uncolored volume satisfies \acrshort{gdop} $\leq 6$. The bad-region volumes are decomposed into connected ``pillars'' in Figs.~\ref{fig:non_lsp_3d_red} and~\ref{fig:non_lsp_3d_black}, with periodic merging across the $(\theta_S, \theta_M)$ boundaries. The red and black pillars are labelled as R1-R3 and B1-B5 within Figs.~\ref{fig:non_lsp_3d_red} and~\ref{fig:non_lsp_3d_black}, respectively. This angular domain construction supplies an alternative, intuitive understanding for the time-domain behavior for \acrshort{gdop} as the time-lines (gold in Fig.~\ref{fig:non_lsp_3d}) encounter different pillars. Figure~\ref{fig:non_lsp_snapshots} illustrates four representative slices at $\alpha_3 \in \{0, \pi/2, \pi, 3\pi/2\}$, each evaluated as a two-dimensional $(\theta_S, \theta_M)$ torus at the corresponding instantaneous \acrshort{mrf} latitude from Eq.~\eqref{eq:lat-mrf}; the longitudinal symmetry of the \acrshort{dadm} trivially absorbs the longitude oscillation of Eq.~\eqref{eq:lon-mrf}, so each slice depends only on the instantaneous latitude. The red and black patterns encountered at each $\alpha_3$ are slices through the pillars identified in Fig.~\ref{fig:non_lsp_3d} and annotated (R1-R3, B1-B5). The coverage varies from $75.4\%$ at the pole-side extremum ($\varphi = -86.68^\circ$, $\alpha_3 = 0^\circ$) down to $61.1\%$ at the equator-side extremum ($\varphi = -73.32^\circ$, $\alpha_3 = 180^\circ$), with the two intermediate slices at $\alpha_3 = 90^\circ$ and $270^\circ$ coinciding by symmetry of Eq.~\eqref{eq:lat-mrf}. 

Although the full 3D angular search is conceptually straightforward, it requires more sampling and may be costly for a preliminary analysis. Three reduced routes are therefore natural.

\begin{enumerate}[label=(\arabic*)]
    \item \ul{Time-domain analysis}: A direct evaluation along the time-lines (gold in Fig.~\ref{fig:non_lsp_3d}) captures the configuration but requires care; the beat frequency $\beta = \nu_M - \omega_3$ is on the order of one year, so sampling over this ``long-enough'' horizon may be necessary to reduce biasing. Sampling at selected epochs to guarantee a sufficient distribution in $\beta$-angles may serve as a proxy.
    \item \ul{Time-averaged latitude}: Since the apparent latitude is a periodic function of $\alpha_3$ (Eq.~\eqref{eq:lat-mrf}), it may be approximated by its average within the rotating frame. Fixing this averaged latitude re-activates the longitudinal invariance within the \acrshort{dadm} and admits a long-term averaged performance estimate, particularly relevant when an averaged figure of merit, rather than a specific window-based requirement (e.g., daily window), is enforced.
    \item \ul{Latitude factorization}: This route is particularly suited to evaluating coverage across many fixed user locations. It is possible to exploit the \acrshort{dadm} longitudinal symmetry: the coverage at any $(\theta_S, \theta_M)$ depends on the user only through the instantaneous latitude ($\varphi$), so the family of bad-region structures across $\alpha_3$ collapses into a one-dimensional coverage function $\Gamma_{2D}(\varphi)$, defined as the fraction of the $(\theta_S, \theta_M)$ torus satisfying the \acrshort{gdop} criterion at an instantaneous user latitude $\varphi$. Each user's latitude $\varphi_{\mathrm{PA}}$ then enters only through the evolution induced by Eq.~\eqref{eq:lat-mrf}. The four snapshots in Fig.~\ref{fig:non_lsp_snapshots} are evaluations of this universal $\Gamma_{2D}(\varphi)$ at the latitudes prescribed for $\varphi_\mathrm{PA} = -80^\circ$, and the averaged coverage reduces to
    \begin{equation}
    \label{eq:non-lsp-cov}
    \overline{\text{cov}}(\varphi_\mathrm{PA}) \;=\; \frac{1}{2\pi}\int_0^{2\pi} \Gamma_{2D}\big(\varphi(\alpha_3;\,\varphi_\mathrm{PA})\big)\,d\alpha_3.
    \end{equation}
    A single $\Gamma_{2D}$ sweep then characterizes the entire family of fixed surface users, at the cost of capturing only torus-averaged figures of merit; window-based requirements (e.g., worst-case daily coverage) are discarded by the integral, since the temporal ordering of latitude visits is removed.
\end{enumerate}
The non-\acrshort{lsp} user geometry enlarges the dimensionality but does not invalidate the frequency-domain framework. Rather, it motivates a hierarchy of reduced approximations matched to the mission metric. Within this perspective, the simplifying assumptions invoked elsewhere in this section are not merely computational shortcuts but defensible early-stage design strategies. The idealization of App.~\ref{app:onaxis} ($\varphi_\mathrm{PA} = -90^\circ$) removes the user's draconic dimension entirely and admits both \acrshort{dadm} and \acrshort{cr3bp} surveys; the time-averaged latitude approximation (avenue~(2)) collapses the instantaneous latitude oscillation to its mean within the \acrshort{dadm} and re-activates the longitudinal invariance for global $(a, i)$ surveys; and the population-averaged formulation (avenue~(3), Eq.~\eqref{eq:non-lsp-cov}) is appropriate when the user-population coverage averaged over the draconic cycle is itself the figure of merit. Each sacrifices a different aspect of strict optimality in exchange for a tractable design space and for rankings that are insensitive to the choice of evaluation horizon. The window-aware time-domain alternative (avenue~(1)) preserves every degree of freedom at a cost that scales unfavorably with both the number of user locations and the simulation horizon, but ultimately re-introduces the same fragility under a different name: any constellation labeled ``optimized'' is so labeled only under the metric and propagation horizon selected during evaluation, and configurations that maximize coverage at one window length need not maximize it at another. This horizon-dependence is not a defect of the present idealizations but a structural feature of threshold-based metrics on quasi-periodic systems with multiple incommensurate frequencies.

\section{Spacing in $\theta_S$}
\label{app:uniform-phase}

The optimized phase distributions in Sec.~\ref{sec:design-exploration} generally exhibit non-uniform short-period spacing $\theta_{S,k}$, particularly at higher inclinations (Fig.~\ref{fig:phase_spread_obl}). This appendix confirms that enforcing uniform $\theta_S$ spacing ($\Delta\theta_S = 72^\circ$ for five satellites) can be substantially suboptimal, especially at higher inclinations. Two reference geometries are examined, both at $a = 14{,}200$~km under the frozen-orbit condition ($\omega = 90^\circ$, $A_L = 0$): a higher-inclination case at $i = 57^\circ$ ($e \approx 0.71$) and a moderate case at $i = 50.5^\circ$ ($e \approx 0.57$). For each, the differential evolution algorithm optimizes the medium-period phases $\theta_{M,k}$ while holding $\theta_S$ at uniform spacing; the result is compared against the unconstrained (free) optimization of both $\theta_S$ and $\theta_M$. The user is placed at the obliquity-corrected \acrshort{lsp} (fixed \acrshort{mrf} latitude $\varphi = -83.32^\circ$ and longitude $\lambda = 0^\circ$); coverage is evaluated on the full two-dimensional $(\theta_S, \theta_M)$ torus on a $500 \times 500$ grid, similar to Fig.~\ref{fig:sweep_obl_torus}.

The summary appears in Table~\ref{tab:uniform_comparison}. For each geometry and optimization mode, the \acrshort{gdop} coverage fraction and the standard deviation for the sorted consecutive $\theta_{S,k}$ gaps (as defined for Fig.~\ref{fig:phase_spread_obl} in Sec.~\ref{sec:design-exploration}) are listed; the standard deviation is a measure of departure from uniform spacing. The final row records the coverage gain from releasing the uniformity constraint. The performance gap depends strongly on the orbital eccentricity. At $i = 57^\circ$ ($e \approx 0.71$), the high eccentricity and inclination confine each satellite's orbit geometry to a narrower region within the $xy$-plane in the \acrshort{mrf} (Fig.~\ref{fig:strobo_hfem}). With uniform $\theta_S$ spacing, the satellites that are simultaneously above the elevation mask tend to cluster in azimuth and elevation as viewed from the south pole, resulting in poor geometric diversity and consequently high \acrshort{gdop}. The optimizer compensates by adopting a non-uniform $\theta_S$ distribution that sacrifices visibility in some portions of the torus in exchange for superior geometric diversity where satellites are visible, yielding an improvement of roughly $32$~pp for this configuration. An optimal configuration discussed in \citet{brack2025} is in line with this observation, where $60^\circ$ spacing is utilized for the first four gaps, resulting in $\theta_{S,1} - \theta_{S, 5} = 120^\circ$. In contrast, at $i = 50.5^\circ$ ($e \approx 0.57$), the uniform configuration achieves $71.7\%$ coverage, only $2.4$~pp below the free optimum, and the near-uniform gap distribution of the free solution ($\mathrm{std} \approx 8.5^\circ$) confirms that uniform spacing is close to optimal at this eccentricity value. The precise boundary between these distinct spacing regimes likely depends on the adopted performance metric; however, the contrast between the two geometries demonstrates that solutions beyond the intuitive evenly spaced distribution can be optimal, reinforcing the value of a broad lower-fidelity search in the preliminary design stage.

\begin{table}[htb]
    \centering
    \caption{Coverage comparison: uniform vs.\ free $\theta_S$ spacing at $a = 14{,}200$~km in the \acrshort{dadm}, with the user at the obliquity-corrected \acrshort{lsp} (\acrshort{mrf} latitude $\varphi = -83.32^\circ$ and longitude $\lambda = 0^\circ$).}
    \label{tab:uniform_comparison}
    \begin{tabular}{l cc cc}
        \toprule
        & \multicolumn{2}{c}{$i = 57^\circ$ ($e \approx 0.71$)} & \multicolumn{2}{c}{$i = 50.5^\circ$ ($e \approx 0.57$)} \\
        \cmidrule(lr){2-3} \cmidrule(lr){4-5}
        & Coverage & $\theta_S$ gap std & Coverage & $\theta_S$ gap std \\
        \midrule
        Uniform $\theta_S$ & $30.9\%$ & $0.0^\circ$ & $71.7\%$ & $0.0^\circ$ \\
        Free $\theta_S$    & $62.5\%$ & $47.3^\circ$ & $74.1\%$ & $8.5^\circ$ \\
        \midrule
        $\Delta$ [pp] & $+31.6$ & & $+2.4$ & \\
        \bottomrule
    \end{tabular}
\end{table}

\section{Initial Conditions for Table~\ref{tab:vz_comparison}}
\label{app:initial-conditions}

The sample constellation analyzed in Sec.~\ref{sec:fddc} shares the mean elements $a = 14{,}200$~km, $e = 0.5707$, $i = 50.50^\circ$, $\omega = 90^\circ$ (\acrshort{eof}) across all five satellites, under the frozen condition $A_L = 0$. The inter-satellite phase variables $(\theta_{S,k}, \theta_{M,k})$ are defined at the initial epoch $t_0 = 0$ (Eq.~\eqref{eq:eof_mrf}, arbitrarily selected without loss of generality), and are mapped to the mean anomaly $M_k = \theta_{S,k}$ and ascending node $\Omega_k = -\theta_{M,k}$, consistent with the convention $\theta_M \approx -\Omega_R$ (Sec.~\ref{subsec:freq-modes}). The per-satellite phase values are $(\theta_{S,k}, \theta_{M,k}) = (0,\,0)$, $(55.65^\circ,\,200.06^\circ)$, $(131.61^\circ,\,121.21^\circ)$, $(209.69^\circ,\,235.29^\circ)$, $(288.03^\circ,\,178.73^\circ)$ for $k=1,\dots,5$.

The osculating Keplerian elements in the \acrshort{mci} are listed in Table~\ref{tab:ic_all} at the common reference epoch $JD_0 = 2460941.5$ (2025-09-23 00:00:00 UTC) for the four initialization methods compared in Sec.~\ref{subsec:fddc-demo}: the na\"{\i}ve placement of the \acrshort{dadm} design means; the analytical Von-Zeipel mean-to-osculating transformation from Nie and Gurfil~\cite{nie2018lunar}; the \acrshort{fddc} refinement that targets only $(\nu_S,\nu_M)$; and the \acrshort{fddc} refinement that additionally suppresses the long-period amplitude to $A_L^* = 10^{-3}$. All four sets are intended for direct propagation in the \acrshort{hfem}. Transformation between the \acrshort{eof} and \acrshort{mci} (applicable for the na\"{\i}ve and semi-analytical transformations) is achieved via Eq.~\eqref{eq:eof_mci_snapshot}. In the last column of Table~\ref{tab:ic_all}, $f$ is the true anomaly. Since the \acrshort{mci} is referred to the ecliptic (Sec.~\ref{sec:prelim}), the tabulated inclinations all remain within a few degrees of $50.5^\circ$, the spread being set by $i_{EOF} \approx 5.15^\circ$ together with the osculating variation, and thereby reflect the common design inclination shared by the five satellites; referred instead to the Earth mean equator, the same five orbits span roughly $24^\circ$ to $79^\circ$ and the common design value is not apparent.

\begin{table}[htb]
    \centering
    \caption{Osculating Keplerian elements in the \acrshort{mci} at $JD_0$ for the four transformation/refinement methods.}
    \label{tab:ic_all}
    \begin{tabular}{cl rrrrrr}
        \toprule
        Method & Sat & $a$ [km] & $e$ & $i$ [deg] & $\Omega$ [deg] & $\omega$ [deg] & $f$ [deg] \\
        \midrule
        \multirow{5}{*}{Na\"{\i}ve}
            & 1 & 14200.0000 & 0.570679 & 45.7760 & 195.2297 &  86.8812 &   0.0000 \\
            & 2 & 14200.0000 & 0.570679 & 55.7438 & 352.7425 &  90.5640 & 123.4185 \\
            & 3 & 14200.0000 & 0.570679 & 51.2496 &  67.6948 &  96.7169 & 163.2966 \\
            & 4 & 14200.0000 & 0.570679 & 55.0974 & 319.8041 &  86.7808 & 190.0143 \\
            & 5 & 14200.0000 & 0.570679 & 55.2492 &  12.6942 &  92.8475 & 224.6268 \\
        \midrule
        \multirow{5}{*}{Semi-analytical~\cite{nie2018lunar}}
            & 1 & 14216.3941 & 0.531376 & 46.8041 & 195.1529 &  86.9342 &   0.0000 \\
            & 2 & 14306.8150 & 0.541933 & 56.5486 & 354.3855 &  89.4142 & 120.1944 \\
            & 3 & 14136.3146 & 0.590531 & 50.8724 &  66.2241 & 100.5586 & 163.6315 \\
            & 4 & 14146.4107 & 0.586530 & 54.7243 & 321.3566 &  82.8010 & 190.0152 \\
            & 5 & 14266.8217 & 0.537741 & 56.2822 &  12.6487 &  91.9924 & 228.9599 \\
        \midrule
        \multirow{5}{*}{\acrshort{fddc}}
            & 1 & 14261.8542 & 0.558874 & 46.2887 & 194.9008 &  89.1108 &   5.4234 \\
            & 2 & 14342.3055 & 0.563861 & 56.4955 & 352.3936 &  93.0438 & 122.1989 \\
            & 3 & 14187.7904 & 0.583407 & 51.6583 &  69.4216 &  99.8083 & 163.7722 \\
            & 4 & 14188.4836 & 0.576141 & 56.6710 & 326.6200 &  80.5862 & 191.3290 \\
            & 5 & 14316.8148 & 0.567657 & 55.9803 &  12.3176 &  96.1957 & 224.8628 \\
        \midrule
        \multirow{5}{*}{\acrshort{fddc} ($A_L^* = 10^{-3}$)}
            & 1 & 14248.9702 & 0.502084 & 48.3980 & 194.7351 &  88.1332 &   3.0138 \\
            & 2 & 14344.8583 & 0.532336 & 57.8229 & 353.1834 &  87.8532 & 122.7181 \\
            & 3 & 14182.0667 & 0.575660 & 52.1616 &  68.0461 &  98.6967 & 163.5679 \\
            & 4 & 14200.2753 & 0.577647 & 56.0622 & 318.9773 &  83.5064 & 191.3009 \\
            & 5 & 14317.6758 & 0.526227 & 57.5052 &  13.3070 &  92.0420 & 232.4260 \\
        \bottomrule
    \end{tabular}
\end{table}

\section{Fourier Decomposition for an HFEM Trajectory}
\label{app:fourier}

The current appendix includes further details on the Fourier decomposition of Sec.~\ref{subsec:fourier-reopt}, intended to serve two purposes. First, the tables below supply the complete coefficient set required to reconstruct the surrogate employed there, rendering the reported coverage reproducible without repeating the \acrshort{hfem} propagation. Second, although the decomposition pertains to a single sample \acrshort{elfo}, the amplitude ordering is structurally informative: the leading peaks along every axis are integer combinations of $\nu_S$ and $\nu_M$ alone, with the ephemeris-driven terms ($p, q, s \neq 0$) entering an order of magnitude lower in amplitude. Since only $m_j$ and $n_j$ enter the inter-satellite shift in Eq.~\eqref{eq:fourier-shift}, this ordering is what renders the per-peak phase rotation effective.

The \acrshort{mrf} Cartesian position for sample satellite~1 (\acrshort{fddc} with $A_L^* = 10^{-3}$, see Appendix~\ref{app:initial-conditions}; $a \approx 14{,}200$~km, $e \approx 0.57$, $i \approx 50.5^\circ$) serves as a reference trajectory. The numerically integrated trajectory within the \acrshort{hfem} over ten years from the initial epoch $JD_0 = 2460941.5$ (2025-09-23 00:00:00 UTC), is decomposed as a finite Fourier series,
\begin{equation}
    l_*r_{\zeta}(t) = l_*\bar r_{\zeta} + \sum_{j=1}^{N_{\zeta}} A_{j,\zeta}\,\cos\!\left( \mathfrak{f}_{j,\zeta} t + \phi_{j,\zeta} \right),\qquad \zeta\in\{x,y,z\}.
    \label{eq:fourier-hfem}
\end{equation}
The peak frequencies $\mathfrak{f}_{j,\zeta}$ (distinct from the internal $\bm{\nu}$ and external $\bm{\omega}$ notation) are identified via iterative refinement following Laskar~\cite{laskar1999introduction}. Subsequently, these peak frequencies are labeled as integer combinations of the five candidate fundamental frequencies,
\begin{equation}
    \mathfrak{f}_{j,\zeta} \approx m_j \nu_S + n_j \nu_M + p_j \omega_1 + q_j \omega_2 + s_j \omega_3,
\end{equation}
where $\nu_S = 15.546$~rad/nd and $\nu_M = 1.085$~rad/nd are the short- and medium-period frequencies targeted by the \acrshort{fddc} (Sec.~\ref{sec:fddc}), and $\omega_1, \omega_2, \omega_3$ are the anomalistic, synodic, and draconic lunar frequencies introduced in Sec.~\ref{subsec:freq-across-models}. With the suppressed long-period motion ($A_L^* = 10^{-3}$), the long-period frequency $\nu_L$ does not appear in the dominant peaks for the signals. The iterative Laskar refinement is configured to extract up to $40$ components per axis; the procedure terminates earlier when a candidate peak falls within two spectral bins of peaks that are already extracted. This early termination occurs after $18$ components on the $z$-axis. Integer labels are then assigned by selecting, among all combinations within a tolerance of $8.2\cdot 10^{-3}$~rad/nd (approximately one Fourier bin width of the ten-year analysis window, set as the smallest distinguishable frequency separation), the one with the smallest sum of absolute integer coefficients. Tables~\ref{tab:fourier-x}-\ref{tab:fourier-z} report such linear combinations as well as the resulting residual, $\Delta \mathfrak{f}$ for the extracted peaks. Since $\omega_1$ and $\omega_3$ are very close, some higher-order labels may admit alternative decompositions of comparable quality. Nevertheless, all extracted peaks are matched within tolerance and yield consistent labels across independent extractions based on 10- and 30-year trajectories. The above choices (number of components, matching tolerance, basis cardinality) are sufficient for the Fourier surrogate constructed in Sec.~\ref{sec:fddc}; further refinement of the decomposition itself, together with its sensitivity to those choices, lies beyond the scope of the present analysis.

\begin{table}[h!]
    \centering
    \caption{Fourier coefficients for $l_* r_x$, \acrshort{hfem} satellite~1 (\acrshort{fddc} with $A_L^* = 10^{-3}$).}
    \label{tab:fourier-x}
    \begin{tabular}{r rrrrr r r r r}
        \toprule
        $j$ & $m$ & $n$ & $p$ & $q$ & $s$ & $\mathfrak{f}$ [rad/nd] & $A$ [km] & $\phi$ [rad] & $\Delta\mathfrak{f}$ [$10^{-4}$ rad/nd] \\
        \midrule
        1 & $ 1$ & $-1$ & $ 0$ & $ 0$ & $ 0$ & $  14.461$ & $   9771.8$ & $+1.665$ & $0.07$ \\
        2 & $ 0$ & $ 1$ & $ 0$ & $ 0$ & $ 0$ & $   1.085$ & $   7974.5$ & $+1.501$ & $0.51$ \\
        3 & $ 2$ & $-1$ & $ 0$ & $ 0$ & $ 0$ & $  30.008$ & $   2446.2$ & $+1.689$ & $0.64$ \\
        4 & $ 1$ & $ 1$ & $ 0$ & $ 0$ & $ 0$ & $  16.631$ & $   1986.9$ & $+1.525$ & $1.08$ \\
        5 & $ 3$ & $-1$ & $ 0$ & $ 0$ & $ 0$ & $  45.554$ & $    913.4$ & $+1.713$ & $1.21$ \\
        6 & $ 1$ & $-1$ & $-1$ & $ 0$ & $ 0$ & $  13.470$ & $    560.1$ & $-0.780$ & $0.10$ \\
        7 & $ 2$ & $ 1$ & $ 0$ & $ 0$ & $ 0$ & $  32.178$ & $    556.1$ & $+1.548$ & $1.65$ \\
        8 & $ 1$ & $-1$ & $ 1$ & $ 0$ & $ 0$ & $  15.453$ & $    549.7$ & $+0.968$ & $0.02$ \\
        9 & $ 0$ & $ 1$ & $ 1$ & $ 0$ & $ 0$ & $   2.076$ & $    515.8$ & $-2.338$ & $0.47$ \\
        10 & $ 0$ & $ 1$ & $-1$ & $ 0$ & $ 0$ & $   0.093$ & $    445.4$ & $+2.199$ & $0.56$ \\
        11 & $ 4$ & $-1$ & $ 0$ & $ 0$ & $ 0$ & $  61.100$ & $    403.6$ & $+1.736$ & $1.79$ \\
        12 & $ 3$ & $ 1$ & $ 0$ & $ 0$ & $ 0$ & $  47.724$ & $    222.4$ & $+1.572$ & $2.22$ \\
        13 & $ 5$ & $-1$ & $ 0$ & $ 0$ & $ 0$ & $  76.647$ & $    195.6$ & $+1.760$ & $2.36$ \\
        14 & $ 2$ & $-3$ & $ 0$ & $ 0$ & $ 0$ & $  27.838$ & $    176.5$ & $-1.313$ & $0.38$ \\
        15 & $ 1$ & $ 1$ & $ 1$ & $ 0$ & $ 0$ & $  17.623$ & $    140.4$ & $-2.314$ & $1.04$ \\
        16 & $ 3$ & $-3$ & $ 0$ & $ 0$ & $ 0$ & $  43.384$ & $    132.3$ & $-1.289$ & $0.19$ \\
        17 & $ 2$ & $-1$ & $ 1$ & $ 0$ & $ 0$ & $  30.999$ & $    130.3$ & $+0.992$ & $0.60$ \\
        18 & $ 2$ & $-1$ & $-1$ & $ 0$ & $ 0$ & $  29.016$ & $    126.5$ & $-0.757$ & $0.67$ \\
        19 & $ 1$ & $ 1$ & $-1$ & $ 0$ & $ 0$ & $  15.640$ & $    123.9$ & $+2.224$ & $1.13$ \\
        20 & $ 1$ & $-1$ & $-1$ & $ 2$ & $ 0$ & $  15.320$ & $    115.6$ & $+2.648$ & $0.10$ \\
        21 & $ 1$ & $-1$ & $ 1$ & $-2$ & $ 0$ & $  13.603$ & $    110.9$ & $-2.458$ & $0.06$ \\
        22 & $ 4$ & $ 1$ & $ 0$ & $ 0$ & $ 0$ & $  63.270$ & $    103.6$ & $+1.596$ & $2.79$ \\
        23 & $ 0$ & $ 1$ & $-1$ & $ 2$ & $ 0$ & $   1.944$ & $    102.1$ & $-0.660$ & $0.51$ \\
        24 & $ 6$ & $-1$ & $ 0$ & $ 0$ & $ 0$ & $  92.193$ & $    100.6$ & $+1.784$ & $2.92$ \\
        25 & $ 0$ & $ 1$ & $ 1$ & $-2$ & $ 0$ & $   0.226$ & $     95.8$ & $+0.517$ & $0.46$ \\
        26 & $ 4$ & $-3$ & $ 0$ & $ 0$ & $ 0$ & $  58.931$ & $     88.3$ & $-1.265$ & $0.77$ \\
        27 & $ 0$ & $ 0$ & $ 0$ & $ 2$ & $-1$ & $   0.846$ & $     70.5$ & $-1.704$ & $0.03$ \\
        28 & $ 1$ & $-1$ & $ 0$ & $-2$ & $ 0$ & $  12.611$ & $     66.4$ & $+1.377$ & $0.04$ \\
        29 & $ 0$ & $ 1$ & $ 0$ & $ 2$ & $ 0$ & $   2.935$ & $     60.8$ & $+1.786$ & $0.47$ \\
        30 & $ 5$ & $-3$ & $ 0$ & $ 0$ & $ 0$ & $  74.477$ & $     57.6$ & $-1.241$ & $1.34$ \\
        31 & $ 7$ & $-1$ & $ 0$ & $ 0$ & $ 0$ & $ 107.740$ & $     53.8$ & $+1.807$ & $3.49$ \\
        32 & $ 0$ & $ 2$ & $ 0$ & $ 0$ & $-1$ & $   1.167$ & $     53.7$ & $-2.333$ & $6.17$ \\
        33 & $ 5$ & $ 1$ & $ 0$ & $ 0$ & $ 0$ & $  78.817$ & $     52.5$ & $+1.619$ & $3.35$ \\
        34 & $ 1$ & $ 0$ & $ 0$ & $ 0$ & $ 0$ & $  15.546$ & $     50.6$ & $+0.077$ & $1.19$ \\
        35 & $ 3$ & $-1$ & $ 1$ & $ 0$ & $ 0$ & $  46.546$ & $     46.3$ & $+1.016$ & $1.18$ \\
        36 & $ 0$ & $ 2$ & $ 0$ & $ 0$ & $ 0$ & $   2.170$ & $     44.6$ & $-0.133$ & $0.55$ \\
        37 & $ 2$ & $ 1$ & $ 1$ & $ 0$ & $ 0$ & $  33.169$ & $     43.2$ & $-2.291$ & $1.60$ \\
        38 & $ 3$ & $-1$ & $-1$ & $ 0$ & $ 0$ & $  44.563$ & $     41.6$ & $-0.735$ & $1.23$ \\
        39 & $ 0$ & $ 0$ & $ 0$ & $ 0$ & $ 1$ & $   1.006$ & $     41.2$ & $-1.656$ & $15.09$ \\
        40 & $ 1$ & $ 0$ & $ 0$ & $ 0$ & $ 1$ & $  16.550$ & $     38.5$ & $+2.229$ & $4.97$ \\
        \bottomrule
    \end{tabular}
\end{table}

\begin{table}[h!]
    \centering
    \caption{Fourier coefficients for $l_* r_y$, \acrshort{hfem} satellite~1 (\acrshort{fddc} with $A_L^* = 10^{-3}$).}
    \label{tab:fourier-y}
    \begin{tabular}{r rrrrr r r r r}
        \toprule
        $j$ & $m$ & $n$ & $p$ & $q$ & $s$ & $\mathfrak{f}$ [rad/nd] & $A$ [km] & $\phi$ [rad] & $\Delta\mathfrak{f}$ [$10^{-4}$ rad/nd] \\
        \midrule
        1 & $ 1$ & $-1$ & $ 0$ & $ 0$ & $ 0$ & $  14.461$ & $   9380.2$ & $+0.094$ & $0.07$ \\
        2 & $ 0$ & $ 1$ & $ 0$ & $ 0$ & $ 0$ & $   1.085$ & $   6510.7$ & $+3.072$ & $0.51$ \\
        3 & $ 2$ & $-1$ & $ 0$ & $ 0$ & $ 0$ & $  30.008$ & $   2412.4$ & $+0.118$ & $0.64$ \\
        4 & $ 1$ & $ 1$ & $ 0$ & $ 0$ & $ 0$ & $  16.631$ & $   1606.9$ & $+3.095$ & $1.07$ \\
        5 & $ 3$ & $-1$ & $ 0$ & $ 0$ & $ 0$ & $  45.554$ & $    931.6$ & $+0.142$ & $1.21$ \\
        6 & $ 1$ & $-1$ & $ 1$ & $ 0$ & $ 0$ & $  15.453$ & $    548.7$ & $-0.602$ & $0.03$ \\
        7 & $ 1$ & $-1$ & $-1$ & $ 0$ & $ 0$ & $  13.470$ & $    531.5$ & $-2.351$ & $0.10$ \\
        8 & $ 0$ & $ 1$ & $-1$ & $ 0$ & $ 0$ & $   0.093$ & $    442.2$ & $-2.515$ & $0.54$ \\
        9 & $ 2$ & $ 1$ & $ 0$ & $ 0$ & $ 0$ & $  32.178$ & $    431.9$ & $+3.119$ & $1.65$ \\
        10 & $ 4$ & $-1$ & $ 0$ & $ 0$ & $ 0$ & $  61.100$ & $    425.9$ & $+0.166$ & $1.79$ \\
        11 & $ 0$ & $ 1$ & $ 1$ & $ 0$ & $ 0$ & $   2.076$ & $    373.5$ & $-0.766$ & $0.47$ \\
        12 & $ 5$ & $-1$ & $ 0$ & $ 0$ & $ 0$ & $  76.647$ & $    213.7$ & $+0.189$ & $2.35$ \\
        13 & $ 3$ & $ 1$ & $ 0$ & $ 0$ & $ 0$ & $  47.724$ & $    167.8$ & $-3.140$ & $2.22$ \\
        14 & $ 2$ & $-3$ & $ 0$ & $ 0$ & $ 0$ & $  27.838$ & $    167.4$ & $-2.884$ & $0.38$ \\
        15 & $ 2$ & $-1$ & $ 1$ & $ 0$ & $ 0$ & $  30.999$ & $    129.9$ & $-0.578$ & $0.61$ \\
        16 & $ 3$ & $-3$ & $ 0$ & $ 0$ & $ 0$ & $  43.384$ & $    128.2$ & $-2.860$ & $0.19$ \\
        17 & $ 2$ & $-1$ & $-1$ & $ 0$ & $ 0$ & $  29.016$ & $    125.0$ & $-2.328$ & $0.67$ \\
        18 & $ 1$ & $ 1$ & $-1$ & $ 0$ & $ 0$ & $  15.640$ & $    124.0$ & $-2.493$ & $1.09$ \\
        19 & $ 6$ & $-1$ & $ 0$ & $ 0$ & $ 0$ & $  92.193$ & $    113.7$ & $+0.213$ & $2.92$ \\
        20 & $ 1$ & $-1$ & $-1$ & $ 2$ & $ 0$ & $  15.320$ & $    113.2$ & $+1.076$ & $0.07$ \\
        21 & $ 1$ & $-1$ & $ 1$ & $-2$ & $ 0$ & $  13.603$ & $    105.0$ & $+2.254$ & $0.06$ \\
        22 & $ 1$ & $ 1$ & $ 1$ & $ 0$ & $ 0$ & $  17.623$ & $     99.0$ & $-0.742$ & $1.04$ \\
        23 & $ 0$ & $ 1$ & $ 1$ & $-2$ & $ 0$ & $   0.226$ & $     89.4$ & $+2.092$ & $0.55$ \\
        24 & $ 4$ & $-3$ & $ 0$ & $ 0$ & $ 0$ & $  58.931$ & $     87.2$ & $-2.836$ & $0.77$ \\
        25 & $ 0$ & $-1$ & $ 0$ & $ 2$ & $ 0$ & $   0.765$ & $     82.5$ & $+0.377$ & $0.21$ \\
        26 & $ 4$ & $ 1$ & $ 0$ & $ 0$ & $ 0$ & $  63.270$ & $     76.1$ & $-3.117$ & $2.78$ \\
        27 & $ 0$ & $ 1$ & $-1$ & $ 2$ & $ 0$ & $   1.944$ & $     73.7$ & $+0.914$ & $0.53$ \\
        28 & $ 1$ & $-1$ & $ 0$ & $ 2$ & $ 0$ & $  16.312$ & $     73.3$ & $-2.748$ & $0.50$ \\
        29 & $ 1$ & $-1$ & $ 0$ & $-2$ & $ 0$ & $  12.611$ & $     63.2$ & $-0.194$ & $0.04$ \\
        30 & $ 7$ & $-1$ & $ 0$ & $ 0$ & $ 0$ & $ 107.740$ & $     63.0$ & $+0.236$ & $3.48$ \\
        31 & $ 5$ & $-3$ & $ 0$ & $ 0$ & $ 0$ & $  74.477$ & $     57.9$ & $-2.812$ & $1.34$ \\
        32 & $ 0$ & $ 0$ & $ 0$ & $ 2$ & $-1$ & $   0.846$ & $     56.1$ & $-0.118$ & $0.44$ \\
        33 & $ 1$ & $ 0$ & $ 0$ & $ 0$ & $ 1$ & $  16.550$ & $     48.5$ & $-2.590$ & $3.06$ \\
        34 & $ 3$ & $-1$ & $ 1$ & $ 0$ & $ 0$ & $  46.546$ & $     45.7$ & $-0.554$ & $1.18$ \\
        35 & $ 0$ & $ 2$ & $ 0$ & $ 0$ & $-1$ & $   1.167$ & $     45.5$ & $-0.749$ & $5.84$ \\
        36 & $ 3$ & $-1$ & $-1$ & $ 0$ & $ 0$ & $  44.563$ & $     44.0$ & $-2.306$ & $1.23$ \\
        37 & $ 0$ & $ 1$ & $ 0$ & $ 2$ & $ 0$ & $   2.935$ & $     42.2$ & $-2.925$ & $0.48$ \\
        38 & $ 0$ & $ 0$ & $ 0$ & $ 0$ & $ 1$ & $   1.005$ & $     41.9$ & $+0.043$ & $11.56$ \\
        39 & $ 6$ & $-3$ & $ 0$ & $ 0$ & $ 0$ & $  90.023$ & $     38.3$ & $-2.788$ & $1.91$ \\
        40 & $ 2$ & $ 1$ & $-1$ & $ 0$ & $ 0$ & $  31.186$ & $     36.2$ & $-2.471$ & $1.66$ \\
        \bottomrule
    \end{tabular}
\end{table}

\begin{table}[h!]
    \centering
    \caption{Fourier coefficients for $l_* r_z$, \acrshort{hfem} satellite~1 (\acrshort{fddc} with $A_L^* = 10^{-3}$).}
    \label{tab:fourier-z}
    \begin{tabular}{r rrrrr r r r r}
        \toprule
        $j$ & $m$ & $n$ & $p$ & $q$ & $s$ & $\mathfrak{f}$ [rad/nd] & $A$ [km] & $\phi$ [rad] & $\Delta\mathfrak{f}$ [$10^{-4}$ rad/nd] \\
        \midrule
        1 & $ 1$ & $ 0$ & $ 0$ & $ 0$ & $ 0$ & $  15.546$ & $   9974.3$ & $+0.024$ & $0.57$ \\
        2 & $ 2$ & $ 0$ & $ 0$ & $ 0$ & $ 0$ & $  31.093$ & $   2484.6$ & $+0.048$ & $1.15$ \\
        3 & $ 3$ & $ 0$ & $ 0$ & $ 0$ & $ 0$ & $  46.639$ & $    935.3$ & $+0.072$ & $1.72$ \\
        4 & $ 0$ & $ 2$ & $ 0$ & $ 0$ & $ 0$ & $   2.170$ & $    535.1$ & $-0.140$ & $1.02$ \\
        5 & $ 4$ & $ 0$ & $ 0$ & $ 0$ & $ 0$ & $  62.185$ & $    418.3$ & $+0.095$ & $2.29$ \\
        6 & $ 1$ & $-2$ & $ 0$ & $ 0$ & $ 0$ & $  13.376$ & $    275.8$ & $+0.164$ & $0.43$ \\
        7 & $ 5$ & $ 0$ & $ 0$ & $ 0$ & $ 0$ & $  77.732$ & $    205.7$ & $+0.119$ & $2.85$ \\
        8 & $ 2$ & $-2$ & $ 0$ & $ 0$ & $ 0$ & $  28.923$ & $    111.9$ & $-2.954$ & $0.11$ \\
        9 & $ 0$ & $ 1$ & $ 0$ & $ 0$ & $-1$ & $   0.081$ & $    108.3$ & $-0.435$ & $0.47$ \\
        10 & $ 3$ & $-2$ & $ 0$ & $ 0$ & $ 0$ & $  44.469$ & $    108.0$ & $-2.930$ & $0.69$ \\
        11 & $ 6$ & $ 0$ & $ 0$ & $ 0$ & $ 0$ & $  93.278$ & $    107.4$ & $+0.142$ & $3.42$ \\
        12 & $ 1$ & $ 2$ & $ 0$ & $ 0$ & $ 0$ & $  17.716$ & $     89.7$ & $-0.116$ & $1.60$ \\
        13 & $ 1$ & $-1$ & $ 0$ & $ 0$ & $ 1$ & $  15.465$ & $     86.5$ & $+0.637$ & $1.62$ \\
        14 & $ 0$ & $ 2$ & $ 1$ & $ 0$ & $ 0$ & $   3.161$ & $     84.4$ & $+2.304$ & $0.97$ \\
        15 & $ 4$ & $-2$ & $ 0$ & $ 0$ & $ 0$ & $  60.016$ & $     77.8$ & $-2.906$ & $1.27$ \\
        16 & $ 7$ & $ 0$ & $ 0$ & $ 0$ & $ 0$ & $ 108.824$ & $     58.5$ & $+0.166$ & $3.98$ \\
        17 & $ 1$ & $ 0$ & $ 2$ & $ 0$ & $-2$ & $  15.522$ & $     54.2$ & $+0.804$ & $2.16$ \\
        18 & $ 1$ & $ 0$ & $-2$ & $ 0$ & $ 2$ & $  15.571$ & $     39.6$ & $-1.242$ & $0.28$ \\
        \bottomrule
    \end{tabular}
\end{table}

\end{document}